\documentclass[12pt]{book}

\newcommand{\dem}{{\bf Preuve\ }}
\newcommand{\fin}{\rule{0.5em}{0.5em}}
\newcommand{\semi}{$\left\{T(t)\right\}_{t\geq 0}\in{\cal SG}(M,\omega)\;\;$}
\newcommand{\dif}{$\left\{T(t)\right\}_{t\geq 0}\in{\cal SGD}(M,\omega)\;\;$}

\begin{document}


\begin{titlepage}
\begin{center}
\large{M\'EMOIRE DE RECHERCHE\\
de Math\'ematiques Pures\\
intitul\'e}\\
\vspace{4cm}
\large{\bf ``LA FORMULE DE LIE - TROTTER\\
\bf POUR LES\\
\bf SEMI-GROUPES FORTEMENT CONTINUS''}\\
\vspace{1cm}
\large{par  {\bf Ludovic Dan LEMLE}\\
sous la direction de  {\bf Gilles CASSIER}}\\
\vspace{3.5cm}
\large{SOUTENU \`A L'UNIVERSIT\'E CLAUDE BERNARD LYON 1}\\
\vspace{5.5cm}
\large{Le 4 juillet 2001}
\end{center}
\newpage

{\bf Remerciements.}

Tout d'abord, je veux profiter de cette occasion pour pr\'esenter mes 
remerciements \`a Monsieur {\bf Gilles Cassier} de l'Universit\'e Claude Benard Lyon 1 pour le choix du sujet de ce m\'emoire, ses suggestions et son aide constante. J'ai \'et\'e tr\`es enchant\'e de cette collaboration avec Monsieur Gilles Cassier.

De m\^eme, je veux exprimer ma reconnaissance \`a Monsieur {\bf Dan Timotin} de l'Institut de Math\'ematiques de l'Acad\'emie Roumaine de Bucharest pour ses excellents conseils pendant son sejour \`a Lyon.

En m\^eme temps, il convient d'exprimer ma gratitude \`a Monsieur {\bf Dumitru Ga\c spar} et \`a Monsieur {\bf Nicolae Suciu} de l'Universit\'e de l'Ouest de Timi\c soara, qui m'ont donn\'e l'occasion d'\'etudier dans une tr\`es grande universit\'e europ\'eenne. 

Finalement, je veux exprimer toute ma reconnaissance pour l'hospitalit\'e et l'amabilit\'e avec lesquelles j'ai \'et\'e accueilli par toutes les personnes que j'ai eu le plaisir de conna\^{\i}tre pendant mon stage \`a l'Universit\'e Claude Bernard Lyon 1.
\vspace{1cm}
\begin{flushright}
Ludovic Dan LEMLE\\
Facultatea de Inginerie\\
Str. Revolu\c tiei$\quad$ Nr. 5\\
COD 2750$\quad$ Hunedoara\\
ROM\^ANIA\\
tel: +4 02 54 20 75 43\\
e-mail: {\sf lemledan@fih.utt.ro}
\end{flushright}

\end{titlepage}

\tableofcontents

\chapter{Introduction}

\section{Pr\'eliminaires}

\hspace{1cm}Dans la suite, nous noterons par $\cal E$ un espace de Banach
sur le corps des nombres complexes ${\bf C}$ et
par $\cal B(E)$ l'alg\`ebre de Banach des op\'erateurs lin\'eaires
born\'es dans $\cal E$. Nous d\'esignerons par $I$ l'unit\'e de $B(E)$.\\
Pour un op\'erateur lin\'eaire $A:{\cal D}(A)\subset{\cal
E}\longrightarrow{\cal E}$ nous noterons par:
$$
{\cal I}\mbox{\em m }A=\{Ax\left|x\in{\cal D}(A)\right.\}
$$
l'image de $A$ et par:
$$
{\cal K}\mbox{\em er }A=\{x\in{\cal D}(A)\left|Ax=0\right.\}
$$
le noyau de $A$.\\
L'op\'erateur $A:{\cal D}(A)\subset{\cal E}\longrightarrow{\cal I}\mbox{\em m }A$ est
surjectif. Si ${\cal K}\mbox{\em er }A=\{0\}$, alors $A$ est injectif. Pour
un op\'erateur bijectif, on peut d\'efinir l'op\'erateur inverse:
$$
A^{-1}:{\cal D}\left(A^{-1}\right)\subset{\cal E}\longrightarrow{\cal
E}
$$
par $A^{-1}y=x$ si $Ax=y$. Evidemment ${\cal
D}\left(A^{-1}\right)={\cal I}\mbox{\em m }A$. Dans la suite nous noterons par $\cal GL(E)$ l'ensemble des \'el\'ements
inversibles de $\cal B(E)$. L'ensemble ${\cal GL(E)}$ est un ensemble ouvert dans ${\cal
B(E)}$ (\cite[Theorem 4.1.13, pag. 145]{istratescu}).\\
Soit $A:{\cal D}(A)\subset{\cal E}\longrightarrow{\cal E}$ un
op\'erateur lin\'eaire. Pour tout $n\in{\bf N}$, nous d\'efinissons:
$$
A^n:{\cal D}(A^n)\longrightarrow {\cal E}
$$
par:
$$
A^0=I\:\: ,\:\: A^1=A\:\: ,...,\:\:A^n=A\left(A^{n-1}\right)\quad,
$$
o\`u:
$$
{\cal D}(A^n)=\left\{x\in{\cal D}(A^{n-1})\left|A^{n-1}x\in{\cal
D}(A)\right.\right\}
$$
quel que soit $n\in{\bf N}$.

\newtheorem{lema}{Lemme}[section]
\begin{lema}\label{num4}
Soit $f:[a,b]\rightarrow \cal E$ une fonction continue. Alors:
$$
\lim_{t\rightarrow 0}\frac{1}{t}\int\limits_a^{a+t}\!\!f(s)\:ds=f(a)\quad.
$$
\end{lema}
\dem
Nous avons:
$$
\left\|\frac{1}{t}\int\limits_{a}^{a+t}\!\!f(s)\:ds-f(a)\right\|=\left\|\frac{1}{t}\int\limits_{a}^{a+t}\![f(s)-f(a)]\:ds\right\|\leq\sup_{s\in[a,a+t]}\left\|f(s)-f(a)\right\|\quad.
$$
L'\'egalit\'e de l'\'enonc\'e r\'esulte de la continuit\'e de l'application $f$.\fin

\begin{lema}\label{num2}
Si $A\in \cal B(E)$ et $\|A\|< 1$, alors $I-A\in \cal GL(E)$ et:
$$\left(I-A\right)^{-1}=\sum\limits_{n=0}^{\infty}A^n\quad.$$
\end{lema}
\dem
Soit $Y_n=I+A+A^2+\ldots+A^n$. Alors:
$$
\left\|Y_{n+p}-Y_n\right\|\leq\frac{\left\|A\right\|^{n+1}}{1-\|A\|}\longrightarrow
0\quad\mbox{pour}\quad n\rightarrow\infty.
$$
Par cons\'equent, $\left(Y_n\right)_{n\in{\bf N}}$ est une suite
  Cauchy. Mais $\cal B(E)$ est une alg\`ebre de Banach. La suite
  $\left(Y_n\right)_{n\in{\bf N}}$ est donc convergente.
Notons  $Y\in \cal B(E)$ sa limite.
De l'\'egalit\'e $(I-A)Y_n=I-A^{n+1}$, il r\'esulte que
  $\lim_{n\rightarrow\infty}(I-A)Y_n=I$, d'o\`u $(I-A)Y=I.$\\
Nous obtenons $Y(I-A)=I$ de fa\c con analogue.\\
Finalement, on voit que $I-A\in\cal GL(E)$ et que
  $\left(I-A\right)^{-1}=\sum\limits_{n=0}^{\infty}A^n.$\fin

\newtheorem{obs}[lema]{Remarque}
\begin{obs}\label{num5}
\em
Si $\|I-A\|<1$, alors $A\in\cal GL(E)$ et
$A^{-1}=\sum\limits_{n=0}^{\infty}(I-A)^n.$
\end{obs}

\newtheorem{definitie}[lema]{D\'efinition}
\begin{definitie}
L'ensemble:
$$\rho (A)=\left \{ \lambda \in {\cal {\bf C}}\left| \left(\lambda I-A \right)^{-1}\mbox {est inversible dans} \in {\cal B(E)}\right.\right\}$$
s'appelle l'ensemble r\'esolvant de $A\in\cal B(E).$
\end{definitie}

\newtheorem{prop}[lema]{Proposition}
\begin{prop}
Soit $A\in \cal B(E)$. Alors $\rho (A)$ est un ensemble ouvert.
\end{prop}
\dem
D\'efinissons l'application:
$$
\phi:\cal {\bf C}\longrightarrow\cal B(E)
$$
par:
$$
\phi(\lambda)=\lambda I-A\quad.
$$
Evidemment, $\phi$ est continue. Si $\lambda \in \rho(A)$, alors
$\lambda I-A\in \cal GL(E)$ et par suite $\rho (A)=\phi^{-1}\left(\cal
GL(E)\right)$. Comme $\cal GL(E)$ est un ensemble ouvert, on voit que $\rho (A)$ est ouvert.\fin

\begin{definitie}
L'application:
$$
R(\:.\:;A):\rho (A)\longrightarrow\cal B(E)
$$
$$
R(\lambda;A)=\left(\lambda I-A\right)^{-1}
$$
s'appelle la r\'esolvante de A.
\end{definitie}

\begin{prop}\label{num3}
La r\'esolvante d'un op\'erateur lin\'eaire $A\in \cal B(E)$, a les propri\'et\'es
suivantes:\\
i) si $\lambda,\mu \in \rho (A)$, alors:
$$
R(\lambda;A)-R(\mu;A)=(\mu - \lambda)R(\lambda;A)R(\mu;A)\quad;
$$
ii) $R(\:.\:;A)$ est une application analytique sur $\rho (A)$;\\
iii) si $\lambda \in \cal{\bf C}$ et $|\lambda|>\|A\|$, alors $\lambda\in
\rho (A)$ et nous avons:
$$
R(\lambda;A)=\sum\limits_{n=0}^{\infty}\frac{A^n}{\lambda ^{n+1}}\quad;
$$
iv) Nous avons:
$$
\frac{d^n}{d\lambda ^n}R(\lambda;A)=(-1)^nn!{R(\lambda;A)}^{n+1}
$$
quels que soient $n\in{\bf N}^*$ et $\lambda \in \rho(A).$
\end{prop}
\dem
i) Nous avons successivement:
\begin{eqnarray*}
& &R(\lambda;A)-R(\mu;A)=\left(\lambda I-A\right)^{-1}-\left(\mu
I-A\right)^{-1}=\\
&=&\left(\lambda I-A\right)^{-1}\left(\mu I-A-\lambda
I+A\right)\left(\mu I-A\right)^{-1}=\\
&=&\left(\mu -\lambda\right)R(\lambda;A)R(\mu;A)
\end{eqnarray*}
quels que soient $\lambda,\mu\in\rho(A).$\\
ii) Soit $\lambda_0 \in \rho(A)$. Notons
$D\left(\lambda_0;\frac{1}{\left\|R(\lambda_0;A)\right\|}\right)$ le
disque ouvert de centre $\lambda_0$ et de rayon $\frac{1}{\|R(\lambda;A)\|}.$ Alors,
pour $\lambda \in
D\left(\lambda_0;\frac{1}{\left\|R(\lambda_0;A)\right\|}\right)$, nous
avons:
$$
\lambda
I-A=\left[I-(\lambda_0-\lambda)R(\lambda_0;A)\right](\lambda_0I-A)\quad.
$$
Mais:
$$
\left\|(\lambda_0-\lambda)R(\lambda_0;A)\right\|=|\lambda_0-\lambda|\|R(\lambda_0;A)\|<
1\quad.
$$
Compte tenu du lemme \ref{num2}, il r\'esulte que:
$$
I-(\lambda_0-\lambda)R(\lambda_0;A)\in \cal GL(E)\quad,
$$
d'o\`u $\lambda I-A\in \cal GL(E)$ et:
\begin{eqnarray*}
{\left(\lambda I-A\right)}^{-1}&=&{\left(\lambda_0I-A\right)}^{-1}{\left[I-(\lambda_0-\lambda)R(\lambda_0;A)\right]}^{-1}=\\
&=&R(\lambda_0;A)\sum\limits_{n=0}^{\infty}{\left(\lambda_0-\lambda\right)}^n{R(\lambda_0;A)}^n=\\&=&\sum\limits_{n=0}^{\infty}(-1)^n{\left(\lambda-\lambda_0\right)}^n{R(\lambda_0;A)}^{n+1}\quad.
\end{eqnarray*}
Donc $R(\:.\:;A)$ est analytique sur $\rho(A)$.\\
iii) Soit $\lambda \in \cal{\bf C}$ tel que $|\lambda|>\|A\|.$ Alors
$\|\lambda^{-1}A\|<1$, d'o\`u ${I-\lambda^{-1}A}\in\cal GL(E)$. De plus:
$$
{\left(I-\lambda^{-1}A\right)}^{-1}=\sum\limits_{n=0}^{\infty}{\left(\lambda^{-1}A\right)}^n=\sum\limits_{n=0}^{\infty}\frac{A^n}{\lambda^{n}}\quad.
$$
Par cons\'equent:
$$
R(\lambda;A)=\left(\lambda I-A\right)^{-1}=\lambda^{-1}\left(I-\lambda^{-1}A\right)^{-1}=\sum\limits_{n=0}^{\infty}\frac{A^n}{\lambda^{n+1}}\quad.
$$
L'assertion (iv) s'obtient par r\'ecurrence. Pour $n=1$, nous avons:
$$
\frac{d}{d\lambda}R(\lambda;A)=\frac{d}{d\lambda}(\lambda
I-A)^{-1}=-(\lambda I-A)^{-2}=R(\lambda;A)^2\quad.
$$
Supposons que pour $k\in{\bf N}$, on ait:
$$
\frac{d^k}{d\lambda^k}R(\lambda;A)=(-1)^kk!R(\lambda;A)^{k+1}\quad.
$$
Montrons que:
$$
\frac{d^{k+1}}{d\lambda^{k+1}}R(\lambda;A)=(-1)^{k+1}(k+1)!R(\lambda;A)^{k+2}\quad.
$$
Nous avons:
\begin{eqnarray*}
&
&\frac{d^{k+1}}{d\lambda^{k+1}}R(\lambda;A)=\frac{d}{d\lambda}\left(\frac{d^k}{d\lambda^k}R(\lambda;A)\right)=\\
&=&\frac{d}{d\lambda}\left[(-1)^kk!(\lambda I-A)^{-k-1}\right]=\\
&=&(-1)^kk!(-k-1)(\lambda I-A)^{-k-2}=(-1)^{k+1}(k+1)!R(\lambda;A)^{k+2}
\end{eqnarray*}
et par cons\'equent:
$$
\frac{d^n}{d\lambda^n}R(\lambda;A)=(-1)^nn!R(\lambda;A)^{n+1}\quad,\quad(\forall)n\in{\bf
N}^*.\fin
$$

\begin{obs}
\em Compte tenu de la proposition \ref{num3} (iii), il r\'esulte que:
$$
\left\{\lambda\in{\bf C}\left||\lambda|>\|A\|\right.\right\}\subset\rho(A).
$$
\end{obs}

\begin{definitie}
L'ensemble $\sigma(A)={\bf C}-\rho(A)$ s'appelle le spectre de
$A\in\cal B(E).$
\end{definitie}

\begin{prop}
Soit $A\in\cal B(E)$. Alors:\\
i) $\sigma(A)\neq\emptyset;$\\
ii) $\sigma(A)$ est un ensemble compact.
\end{prop}
\dem
i) Supposons que $\sigma(A)=\emptyset$. Alors $\rho(A)={\bf
C}$. Par cons\'equent, l'application $\lambda\longmapsto(\lambda I-A)^{-1}$ est
d\'efinie sur ${\bf C}$. De plus, pour $|\lambda|>\|A\|$, nous avons:
$$
R(\lambda;A)=\sum\limits_{n=0}^{\infty}\frac{A^n}{\lambda^{n+1}}\quad,\quad(\forall)\lambda\in\rho(A).
$$
Il s'ensuit que:
$$
\lim_{|\lambda|\rightarrow\infty}R(\lambda;A)=0.
$$
Donc il existe $M>0$ tel que $\|R(\lambda;A)\|<M$, $(\forall)\lambda\in\cal
{\bf C}$. Le th\'eor\`eme de Liouville (\cite[pag. 231]{dunford-schwartz}) implique que
$R(\:.\:;A)$ est constante sur $\cal {\bf C}$ et que cette constante ne
peut \^etre que $0$. Donc ${\left(\lambda I-A\right)}^{-1}=0$ pour
tout $\lambda\in{\bf C}$, ce qui est absurde.
Par cons\'equent $\sigma(A)\neq\emptyset.$\\
ii) Compte tenu de la proposition \ref{num3} (iii), nous obtenons que:
$$
\sigma(A)\subset\left\{\lambda\in{\bf C}\left||\lambda|\leq\|A\|\right.\right\}.
$$
L'ensemble $\sigma(A)$ est donc born\'e. Comme nous avons vu que $\sigma(A)$
est un ensemble ferm\'e, il est donc compact.\fin

\begin{definitie}
Pour un op\'erateur lin\'eaire $A\in{\cal B(E)}$, le nombre
$$
r(A)=\sup_{\lambda\in{\sigma(A)}}|\lambda|
$$
s'appelle le rayon spectral de $A$.
\end{definitie}

\begin{obs}
\em
Evidemment, pour un op\'erateur $A\in{\cal B(E)}$, $\sigma(A)$ est contenu dans l'int\'erieur du cercle de centre $O$ et de rayon $r(A)$. De plus, on peut montrer que
$$
r(A)=\lim_{n\rightarrow\infty}\left\|A^n\right\|^{\frac{1}{n}}
$$
et on voit que $r(A)\leq\|A\|$.
\end{obs}

Par la suite, nous pr\'esenterons quelques probl\`emes concernant la
th\'eorie spectrale pour un op\'erateur lin\'eaire ferm\'e $A:{\cal
D}(A)\subset{\cal E}\longrightarrow {\cal E}.$

\begin{definitie}
L'ensemble:
$$\rho(A)=\left\{\lambda\in{\bf C}\left|\lambda I-A:{\cal
D}(A)\longrightarrow{\cal E}\mbox{ est op\'erateur bijectif}\right.\right\}
$$
s'appele l'ensemble r\'esolvant de $A$.
\end{definitie}

\begin{obs}
\em Il r\'esulte du th\'eor\`eme du graphe ferm\'e (\cite[Theorem II.2.4, pag. 57]
{dunford-schwartz}) que l'op\'erateur:
$$
{\left(\lambda I-A\right)}^{-1}:{\cal E}\longrightarrow {\cal E}
$$
est continu dans $\cal E$.
\end{obs}

\begin{definitie}
L'application:
$$
R(\:.\:;A):\rho(A)\longrightarrow \cal B(E)
$$
$$
R(\lambda;A)={\left(\lambda
I-A\right)}^{-1}\quad,\quad(\forall)\lambda\in\rho(A)
$$
s'appelle la r\'esolvante de $A$.
\end{definitie}


\begin{prop}\label{num14}
Soit $A:{\cal D}(A)\subset{\cal E}\longrightarrow {\cal E}$,
un op\'erateur lin\'eaire ferm\'e. Alors:\\
i) $\rho(A)$ est un ensemble ouvert et $R(\:.\:;A)$ est une application
analytique sur $\rho(A);$\\
ii) si $\lambda,\mu\in\rho(A)$, alors:
$$
R(\lambda;A)-R(\mu;A)=(\mu -\lambda)R(\lambda;A)R(\mu;A)\quad;
$$
iii) Nous avons:
$$
\frac{d^n}{d\lambda^n}R(\lambda;A)=(-1)^nn!{R(\lambda;A)}^{n+1}
$$
quels que soient $n\in{\bf N}$ et $\lambda\in\rho(A).$
\end{prop}
\dem
Elle est analogue \`a celle de la proposition \ref{num3}.\fin

\begin{definitie}
L'ensemble $\sigma(A)={\bf C}-\rho(A)$ s'appelle le spectre de
$A$.
\end{definitie}

\begin{obs}
\em $\sigma(A)$ est un ensemble ferm\'e.
\end{obs}

\begin{obs}
\em Il existe des op\'erateurs ferm\'es qui ont un spectre non
 born\'e.
\end{obs}

\newtheorem{exemplu}[lema]{Exemple}
\begin{exemplu}
\em
Prenons ${\cal E}={\cal C}_{[0,1]}$ et consid\'erons l'op\'erateur:
$$
D:{\cal C}^1_{[0,1]}\longrightarrow \cal E
$$
d\'efini par:
$$
Df=f^{'}
$$
Dans ce cas, nous avons $\sigma(D)={\bf C}$.
\end{exemplu}

\newtheorem{teo}[lema]{Th\'eor\`eme}

\begin{definitie}
Soit ${\bf D}\subset{\bf C}$ un ensemble ouvert. Une application
analytique:
$$
{\bf D}\ni\lambda\longmapsto R_\lambda\in{\cal B(E)}
$$
qui v\'erifie la propri\'et\'e:
$$
R_\lambda-R_\mu=(\mu-\lambda)R_\lambda
R_\mu\quad,\quad(\forall)\lambda,\mu\in{\bf D},
$$
s'appelle une pseudo-r\'esolvante.
\end{definitie}

\begin{teo}\label{num19}
Soit ${\bf D}\ni\lambda\longmapsto R_\lambda\in{\cal B(E)}$
 une pseudo-r\'esolvante. Alors:\\
i) $R_\lambda R_\mu=R_\mu R_\lambda$, $(\forall)\lambda,\mu\in {\bf D}$;\\
ii) ${\cal K}er R_\lambda$ et ${\cal I}m R_\lambda$ ne d\'ependent pas de
$\lambda\in {\bf D}$;\\
iii) $R_\lambda$ est la r\'esolvante d'un op\'erateur lin\'eaire $A$
ferm\'e et d\'efini sur un sous espace dense si et seulement si ${\cal
K}er R_\lambda=\{0\}$ et
$\overline{{\cal I}{m R_\lambda}}={\cal E}$.
\end{teo}
\dem
i) Soient $\lambda,\mu\in{\bf D}$. Alors, nous avons:
$$
R_\lambda-R_\mu=(\mu-\lambda)R_\lambda R_\mu
$$
et:
$$
R_\mu-R_\lambda=(\lambda-\mu)R_\mu R_\lambda\quad,
$$
d'o\`u:
$$
0=(\mu-\lambda)R_\lambda R_\mu+(\lambda-\mu)R_\mu R_\lambda\quad.
$$
Par suite, on a $R_\lambda R_\mu=R_\mu R_\lambda$.\\
ii) Soient $\mu\in{\bf D}$ et $x\in{\cal K}\mbox{\em er }R_\mu$. Alors $R_\mu
x=0$. Si $\lambda\in{\bf D}$, on a:
$$
R_\lambda x-R_\mu x=(\mu-\lambda)R_\lambda R_\mu x\quad.
$$
Donc $R_\lambda x=0$. Par cons\'equent $x\in{\cal K}\mbox{\em er }R_\lambda$. Il
s'ensuit  que ${\cal K}\mbox{\em er }R_\lambda$ ne d\'epend pas de
$\lambda\in{\bf D}$.\\
Soient $\mu\in{\bf D}$ et $y\in{\cal I}\mbox{\em m }R_\mu$. Alors il existe
$x\in{\cal E}$ tel que $R_\mu x=y$. Si $\lambda\in{\bf D}$, nous
avons:
$$
R_\lambda x-R_\mu x=(\mu-\lambda)R_\lambda R_\mu x\quad.
$$
Donc:
$$
R_\lambda x-y=(\mu-\lambda)R_\lambda y\quad,
$$
ou bien:
$$
y=R_\lambda\left(x+(\lambda-\mu)y\right)\quad.
$$
Donc il existe $u=x+(\lambda-\mu)y\in{\cal E}$ tel que $y=R_\lambda
u$. Par cons\'equent $y\in{\cal I}\mbox{\em m }R_\lambda$. Il s'ensuit
que ${\cal I}\mbox{\em m }R_\lambda$ ne d\'epend pas de
$\lambda\in{\bf D}$.\\
iii) $\Longrightarrow$ Si $R_\lambda$ est une r\'esolvante pour
un op\'erateur lin\'eaire $A$ ferm\'e et d\'efini sur un sous espace dense, alors
$R_\lambda$ est une application bijective, d'o\`u ${\cal K}\mbox{\em er
}R_\lambda=\{0\}$ et $R_\lambda=\left(\lambda
I-A\right)^{-1}$. Par suite, ${R_\lambda}^{-1}=\lambda I-A$ et $\overline{{\cal
D}\left({R_\lambda}^{-1}\right)}=\overline{{\cal D}(A)}={\cal E}$. Par
cons\'equent $\overline{{\cal I}\mbox{\em m }R_\lambda}=\overline{{\cal
D}\left({R_\lambda}^{-1}\right)}={\cal E}$.\\
$\Longleftarrow$ Soient ${\bf D}\ni\lambda\longmapsto R_\lambda\in{\cal B(E)}
$ une
pseudo-r\'esolvante et $\lambda\in{\bf D}$ tel que ${\cal K}\mbox{\em er
}R_\lambda=\{0\}$. Alors pour $y\in{\cal I}\mbox{\em m }R_\lambda$, il
existe un seul $x_\lambda\in{\cal E}$ tel que $y=R_\lambda
x_\lambda$. Mais pour $\lambda,\mu\in{\bf D}$, on a:
$$
R_\lambda y-R_\mu y=(\mu-\lambda)R_\lambda R_\mu y\quad.
$$
D'autre part:
\begin{eqnarray*}
& &R_\lambda y-R_\mu y=R_\lambda R_\mu x_\mu-R_\mu R_\lambda
x_\lambda=\\
&=&R_\lambda R_\mu x_\mu-R_\lambda R_\mu x_\lambda=R_\lambda
R_\mu\left(x_\mu-x_\lambda\right)\quad.
\end{eqnarray*}
Donc $x_\mu-x_\lambda=(\mu-\lambda)y$, d'o\`u $\lambda y-x_\lambda=\mu
y-x_\mu$.
Par cons\'equent, l'op\'erateur:
$$
A:{\cal I}\mbox{\em m }R_\lambda\longrightarrow{\cal E}
$$
$$
Ay=\lambda y-x_\lambda=\lambda y-{R_\lambda}^{-1}y
$$
est correctement d\'efini (valeur ind\'ependante de $\lambda$). De
m\^eme $\overline{{\cal D}(A)}=\overline{{\cal I}\mbox{\em m
}R_\lambda}={\cal E}$. Puis que $R_\lambda\in{\cal B(E)}$, il
r\'esulte du th\'eor\`eme du graphe ferm\'e
(\cite[Theorem II.2.4, pag. 57]{dunford-schwartz}) que $R_\lambda^{-1}$ est un op\'erateur
ferm\'e. Donc $A=\lambda I-{R_\lambda}^{-1}$ est un op\'erateur
ferm\'e. De plus, on a:
$$
{R_\lambda}^{-1}y=x_\lambda=\lambda y-Ay=(\lambda I-A)y\quad.
$$
Par cons\'equent $R_\lambda=(\lambda I-A)^{-1}$ est la r\'esolvante de $A$.\fin
\vspace{2cm}

  \section{Les op\'erateurs dissipatifs}

\hspace{1cm}Dans la suite, nous notons par ${\cal E}^*$ l'espace dual du ${\cal
E}$ et par $\|\:.\:\|_*$ sa norme. Pour tout $x\in{\cal E}$, nous
d\'esignerons par ${\cal J}(x)$ l'ensemble:
$$
\left\{x^*\in{\cal
E}^*\left|\langle x,x^*\rangle=\|x\|^2=\|x^*\|_*^2\right.\right\}\quad.
$$

\begin{definitie}
On dit que l'op\'erateur lin\'eaire $A:{\cal D}(A)\subset{\cal
E}\longrightarrow{\cal E}$ est dissipatif si pour tout $x\in{\cal
D}(A)$, il existe $x^*\in{\cal J}(x)$ tel que
$\mbox{Re}\langle Ax,x^*\rangle\leq 0$.
\end{definitie}
Dans la proposition suivante nous pr\'esentons une caract\'erisation
tr\`es utile pour les op\'erateurs dissipatifs.

\begin{prop}\label{num34}
Un op\'erateur lin\'eaire $A:{\cal D}(A)\subset{\cal
E}\longrightarrow{\cal E}$ est dissipatif si et seulement si pour tout
$\alpha>0$ on a:
$$
\left\|(\alpha
I-A)x\right\|\geq\alpha\|x\|\quad,\quad(\forall)x\in{\cal D}(A).
$$
\end{prop}
\dem
$\Longrightarrow$ Supposons que $A:{\cal D}(A)\subset{\cal
E}\longrightarrow{\cal E}$ est un op\'erateur dissipatif. Pour tout
$x\in{\cal D}(A)$, il existe $x^*\in{\cal J}(x)$ tel que $\mbox{\em
Re}\langle Ax,x^*\rangle\leq 0$. Si $\alpha>0$, alors nous avons:
\begin{eqnarray*}
& &\left\|(\alpha I-A)x\right\|\|x\|=\left\|(\alpha
I-A)x\right\|\left\|x^*\right\|_*\geq\\
&\geq&\left|\langle(\alpha I-A)x,x^*\rangle\right|\geq\mbox{\em
Re}\langle(\alpha I-A)x,x^*\rangle=\\
&=&\mbox{\em Re}\langle\alpha x,x^*\rangle-\mbox{\em Re}\langle
Ax,x^*\rangle\geq\alpha\|x\|^2\quad,
\end{eqnarray*}
d'o\`u il r\'esulte l'in\'egalit\'e de l'\'enonc\'e.\\
$\Longleftarrow$ Soit $A:{\cal D}(A)\subset{\cal
E}\longrightarrow{\cal E}$ tel que pour tout $\alpha>0$ et $x\in{\cal
D}(A)$ on ait:
$$
\left\|(\alpha
I-A)x\right\|\geq\alpha\|x\|\quad.
$$
Soit $y_\alpha^*\in{\cal J}\left((\alpha I-A)x\right)$. On a donc:
$$
\langle(\alpha I-A)x,y_\alpha^*\rangle=\left\|(\alpha I-A)x\right\|^2=\left\|y_\alpha^*\right\|_*^2\quad,
$$
d'o\`u:
$$
\left\|y_\alpha^*\right\|_*=\left\|(\alpha
I-A)x\right\|\geq\alpha\|x\|\quad.
$$
Nous d\'efinissons:
$$
z_\alpha^*=\frac{y_\alpha^*}{\left\|y_\alpha^*\right\|_*}
$$
et d\'esignons par ${\cal B}_1({\cal E}^*)$ la boule unit\'e de
${\cal E}^*$ et par $\partial{\cal B}_1({\cal E}^*)$ sa
fronti\`ere. Il est \'evident que $z_\alpha^*\in\partial{\cal
B}_1({\cal E}^*)$. De plus:
\begin{eqnarray*}
\alpha\|x\|&\leq&\left\|(\alpha
I-A)x\right\|=\frac{1}{\left\|y_\alpha^*\right\|_*}\langle(\alpha
I-A)x,y_\alpha^*\rangle=\\
&=&\langle(\alpha I-A)x,z_\alpha^*\rangle
\end{eqnarray*}
et par cons\'equent:
\begin{eqnarray*}
& &\alpha\|x\|\leq\mbox{\em Re}\langle(\alpha
I-A)x,z_\alpha^*\rangle=\mbox{\em Re}\langle\alpha x,z_\alpha^*\rangle-\mbox{\em Re}\langle
Ax,z_\alpha^*\rangle\leq\\
&\leq&\alpha\left|\langle x,z_\alpha^*\rangle\right|-\mbox{\em Re}\langle
Ax,z_\alpha^*\rangle\leq\alpha\|x\|\left\|z_\alpha^*\right\|_*-\mbox{\em Re}\langle
Ax,z_\alpha^*\rangle=\\
&=&\alpha\|x\|-\mbox{\em Re}\langle
Ax,z_\alpha^*\rangle\quad.
\end{eqnarray*}
Il s'ensuit que:
$$
\mbox{\em Re}\langle Ax,z_\alpha^*\rangle\leq 0\quad,
$$
d'o\`u:
$$
-\mbox{\em Re}\langle
Ax,z_\alpha^*\rangle\leq\left|\langle Ax,z_\alpha^*\rangle\right|\leq\|Ax\|\left\|z_\alpha^*\right\|_*=\|Ax\|
$$
et par cons\'equent:
$$
\alpha\|x\|\leq\alpha\mbox{\em Re}\langle
x,z_\alpha^*\rangle+\|Ax\|\quad.
$$
Donc:
$$
\mbox{\em Re}\langle
x,z_\alpha^*\rangle\geq\|x\|-\frac{1}{\alpha}\|Ax\|\quad.
$$
D'autre part, en appliquant le th\'eor\`eme d'Alaoglu
(\cite[Theorem V.4.2, pag. 424]{dunford-schwartz}), on voit que la
boule unit\'e ${\cal B}_1({\cal E}^*)$ est faiblement compacte. Par
cons\'equent, il existe une sous suite
$\left(z_\beta^*\right)_{\beta>0}\subset\left(z_\alpha^*\right)_{\alpha>0}$
et il existe $z^*\in{\cal B}_1({\cal E}^*)$ tel que:
$$
z_\beta^*\longrightarrow z^*\quad\mbox{si}\quad\beta\rightarrow\infty
$$
pour la topologie faible. Comme on a
$$
\mbox{\em Re}\langle
Ax,z_\beta^*\rangle\leq 0
$$
et
$$
\mbox{\em Re}\langle
x,z_\beta^*\rangle\geq\|x\|-\frac{1}{\beta}\|Ax\|\quad,
$$
on obtient par passage \`a limite en $\beta\rightarrow\infty$:
$$
\mbox{\em Re}\langle
Ax,z^*\rangle\leq 0
$$
et:
$$
\mbox{\em Re}\langle x,z^*\rangle\geq\|x\|\quad.
$$
Mais comme:
$$
\mbox{\em Re}\langle x,z^*\rangle\leq\left|\langle
x,z^*\rangle\right|\leq\|x\|\|z^*\|_*\leq\|x\|\quad,
$$
il s'ensuit que:
$$
\langle x,z^*\rangle=\|x\|\quad.
$$
Si nous prenons $x^*=\|x\|z^*$, il vient:
$$
\langle x,x^*\rangle=\langle x,\|x\|z^*\rangle=\|x\|\langle
x,z^*\rangle=\|x\|^2\quad.
$$
Il en r\'esulte que $x^*\in{\cal J}(x)$. Finalement, on voit que $\mbox{\em Re}\langle
Ax,x^*\rangle\leq 0$, d'o\`u l'on tire que l'op\'erateur $A$ est
dissipatif.\fin

\begin{prop}\label{num35}
Soit $A:{\cal D}(A)\subset{\cal
E}\longrightarrow{\cal E}$ un op\'erateur dissipatif. S'il existe
$\alpha_0>0$ tel que ${\cal I}\mbox{m }(\alpha_0I-A)={\cal E}$, alors
pour tout $\alpha>0$ on a ${\cal I}\mbox{m }(\alpha I-A)={\cal E}$.
\end{prop}
\dem
Soient $A:{\cal D}(A)\subset{\cal
E}\longrightarrow{\cal E}$ un op\'erateur dissipatif et $\alpha_0>0$
tel que ${\cal I}\mbox{\em m }(\alpha_0I-A)={\cal E}$. Compte tenu de
la proposition \ref{num34}, on voit que:
$$
\left\|(\alpha_0I-A)x\right\|\geq\alpha_0\|x\|\quad,\quad(\forall)x\in{\cal
D}(A)
$$
et comme ${\cal I}\mbox{\em m }(\alpha_0I-A)={\cal E}$, il en r\'esulte
que $\alpha_0I-A\in{\cal GL(E)}$ et $\alpha_0$ appartient donc bien \`a $\rho(A)$. 
Soit $(x_n)_{n\in{\bf N}}\subset{{\cal D}(A)}$ tel que $x_n\longrightarrow x$
et $Ax_n\longrightarrow y$ si $n\rightarrow\infty$. Il est clair que:
$$
(\alpha_0I-A)x_n\longrightarrow\alpha_0x-y\quad\mbox{si}\quad n\rightarrow\infty
$$
et par cons\'equent:
$$
x_n=R(\alpha_0;A)(\alpha_0I-A)x_n\longrightarrow
R(\alpha_0;A)(\alpha_0x-y)\quad\mbox{si}\quad n\rightarrow\infty.
$$
Par suite, nous obtenons:
$$
R(\alpha_0;A)(\alpha_0x-y)=x\quad.
$$
Comme ${\cal I}\mbox{\em m }R(\alpha_0;A)\subset{{\cal D}(A)}$, on voit que $x\in{\cal D}(A)$. De plus:
$$
(\alpha_0I-A)x=\alpha_0x-y\quad,
$$
d'o\`u il r\'esulte que $Ax=y$. Par cons\'equent, $A$ est un op\'erateur
ferm\'e.\\
Nous d\'esignerons par ${\cal A}$ l'ensemble:
$$
\left\{\alpha\in]0,\infty)\left|{\cal I}\mbox{\em m}(\alpha I-A)={\cal
E}\right.\right\}\quad.
$$
Soit $\alpha\in{\cal A}$. Comme $A$ est un op\'erateur dissipatif, on voit
que:
$$
\left\|(\alpha
I-A)x\right\|\geq\alpha\|x\|\quad,\quad(\forall)x\in{\cal D}(A),
$$
d'o\`u il r\'esulte que $\alpha\in\rho(A)$. Puisque $\rho(A)$ est un
ensemble ouvert, il existe un voisinage ${\cal V}$ de $\alpha$
contenu dans $\rho(A)$. Comme ${\cal V}\cap]0,\infty)\subset{\cal A}$, on voit que ${\cal A}$ est un ensemble ouvert.\\
Soit $(\alpha_n)_{n\in{\bf
N}}\subset{\cal A}$ tel que $\alpha_n\longrightarrow\alpha$ si
$n\rightarrow\infty$. Comme ${\cal I}\mbox{\em m }(\alpha_nI-A)={\cal
E}$, $(\forall)n\in{\bf N}$, on observe que pour tout $y\in{\cal E}$, il
existe $x_n\in{\cal D}(A)$ tel que:
$$
(\alpha_nI-A)x_n=y\quad,\quad(\forall)n\in{\bf N},
$$
et par suite, il existe $C>0$ tel que:
$$
\|x_n\|\leq\frac{1}{\alpha_n}\|y\|\leq C\quad,\quad(\forall)n\in{\bf
N}.
$$
Par cons\'equent:
\begin{eqnarray*}
&
&\alpha_n\|x_n-x_m\|\leq\|(\alpha_mI-A)(x_n-x_m)\|=\\
&=&\|(\alpha_mI-A)x_n-(\alpha_mI-A)x_m\|=\|\alpha_mx_n-Ax_n-y\|=\\
&=&\|\alpha_mx_n-\alpha_nx_n+\alpha_nx_n-Ax_n-y\|=\\
&=&\|(\alpha_m-\alpha_n)x_n+y-y\|=|\alpha_m-\alpha_n|\|x_n\|\leq
C|\alpha_m-\alpha_n|\quad,
\end{eqnarray*}
d'o\`u il r\'esulte que $(x_n)_{n\in{\bf N}}$ est une suite de
Cauchy. Puisque ${\cal E}$ est un espace de Banach, il s'ensuit que
$(x_n)_{n\in{\bf N}}$ converge vers un point $x\in{\cal E}$. Alors, on en d\'eduit que:
$$
Ax_n\longrightarrow\alpha x-y\quad\mbox{si}\quad n\rightarrow\infty
$$
et comme $A$ est un op\'erateur ferm\'e, on obtient $x\in{\cal D}(A)$
et $\alpha x-Ax=y$. Par suite, ${\cal I}\mbox{\em m }(\alpha
I-A)={\cal E}$ et $\alpha\in{\cal A}$. Donc ${\cal A}$ est ferm\'e dans
$]0,\infty)$ et comme il existe $\alpha_0\in{\cal A}$, nous
d\'eduisons que ${\cal A}=]0,\infty)$.\fin
\vspace{2cm}

  \section{Semi-groupes uniform\'ement continus}

\hspace{1cm}Dans la suite nous pr\'esenterons quelques probl\`emes concernant les
  semi-groupes uniform\'ement continus d'op\'erateurs lin\'eaires born\'es
  sur un espace de Banach $\cal E$.

\begin{definitie}
On appelle semi-groupe uniform\'ement continu d'op\'erateurs li-\\n\'eaires
born\'es sur $\cal E$ une famille
$\left\{T(t)\right\}_{t\geq0}\subset\cal B(E)$ v\'erifiant les propri\'et\'es
sui-\\vantes:\\
i) $T(0)=I;$\\
ii) $T(t+s)=T(t)T(s)\quad,\quad(\forall) t,s\geq 0;$\\
iii) $\lim_{t\searrow 0}\left\|T(t)-I\right\|=0.$
\end{definitie}

\begin{definitie}
On appelle g\'en\'erateur infinit\'esimal du semi-groupe uniform\'ement
continu $\left\{T(t)\right\}_{t\geq 0}$ l'op\'erateur lin\'eaire:
$$
A:{\cal E}\longrightarrow{\cal E}\quad,
$$
$$A=\lim_{t\searrow 0}\frac{T(t)-I}{t}\quad.
$$
\end{definitie}

\begin{lema}\label{num6}
Soit $A\in\cal B(E)$. Alors ${\left\{e^{tA}\right\}}_{t\geq 0}$ est un
semi-groupe uniform\'ement continu d'op\'erateurs lin\'eaires born\'es sur
$\cal E$ dont le g\'en\'erateur infinit\'esimal est $A$.
\end{lema}
\dem
Soit $A\in\cal B(E)$ et $[0,\infty)\ni t\longmapsto T(t)\in\cal B(E)$ une
application d\'efinie par:
$$
T(t)=e^{tA}=\sum\limits_{k=0}^{\infty}\frac{t^kA^k}{k!}\quad.
$$
La s\'erie du membre de droite de l'\'egalit\'e est convergente pour la
topologie de la norme de ${\cal B(E)}$. De plus, il est \'evident que
$T(0)=I$ et $T(t+s)=T(t)T(s)$ quels que soient $t,s\geq 0.$ \\
Compte tenu de l'in\'egalit\'e:
$$
\left\|T(t)-I\right\|\leq e^{t\|A\|}-1\quad,\quad(\forall) t\geq 0,
$$
il r\'esulte:
$$
\lim_{t\searrow 0}\left\|T(t)-I\right\|=0\quad.
$$
Donc la famille $\left\{T(t)\right\}_{t\geq 0}\subset\cal B(E)$ est un
semi-groupe uniform\'ement continu.\\
D'autre part, puisque:
\begin{eqnarray*}
&
&\left\|\frac{T(t)-I}{t}-A\right\|=\left\|\frac{1}{t}\left(e^{tA}-I-tA\right)\right\|=\left\|\frac{1}{t}\left(\sum\limits_{k=0}^{\infty}\frac{t^kA^k}{k!}-I-tA\right)\right\|=\\
&=&\left\|\frac{1}{t}\left(I+tA+\sum\limits_{k=2}^{\infty}\frac{t^kA^k}{k!}-I-tA\right)\right\|\leq\frac{1}{t}\sum\limits_{k=2}^{\infty}\frac{t^k\|A\|^k}{k!}=\\
&=&\frac{1}{t}\left(1+t\|A\|+\sum\limits_{k=2}^{\infty}\frac{t^k\|A\|^k}{k!}-1-t\|A\|\right)=\frac{1}{t}\left(e^{t\|A\|}-1-t\|A\|\right)=\\
&=&\frac{e^{t\|A\|}-1}{t\|A\|}\|A\|-\|A\|\longrightarrow
0\quad\mbox{si}\quad t\searrow 0,
\end{eqnarray*}
nous obtenons:
$$
\lim_{t\searrow 0}\frac{T(t)-I}{t}=A\quad.
$$
Le semi-groupe $\left\{T(t)\right\}_{t\geq 0}$ admet donc pour g\'en\'erateur
infinit\'esimal l'op\'erateur $A.$\fin

\begin{lema}\label{num7}
Etant donn\'e un op\'erateur $A\in\cal B(E)$, il existe un unique semi-groupe uniform\'ement continu $\left\{T(t)\right\}_{t\geq 0}$ tel que:
$$
T(t)=e^{tA}\quad,\quad(\forall) t\geq 0.
$$
\end{lema}
\dem
Si $\left\{S(t)\right\}_{t\geq 0}$ est un autre semi-groupe uniform\'ement
continu engendr\'e par $A$, nous avons:
$$
\lim_{t\searrow 0}\frac{T(t)-I}{t}=A
$$
et:
$$
\lim_{t\searrow 0}\frac{S(t)-I}{t}=A\quad.
$$
Par cons\'equent:
$$
\lim_{t\searrow 0}\left\|\frac{T(t)-S(t)}{t}\right\|=0\quad.
$$
Pour $a\in ]0,\infty)$, nous consid\'erons l'intervalle
$I_a=[0,a[$. Comme $\left\{T(t)\right\}_{t\geq 0}$ et
$\left\{S(t)\right\}_{t\geq 0}$ sont des semi-groupes uniform\'ement
continus, nous voyons que les applications:
$$
t\longmapsto \|T(t)\|
$$
et:
$$
t\longmapsto \|S(t)\|
$$
sont continues. Il existe $c_a\in [1,\infty)$ tel que:
$$
\sup_{t\in I_a}\left\{\|T(t)\|,\|S(t)\|\right\}\leq c_a\quad.
$$
Si $\varepsilon>0$, il existe $t_0\in I_a$, $t_0>0$,
tel que:
$$
\left\|\frac{T(t)-S(t)}{t}\right\|\leq\frac{\varepsilon}{ac_a^2}\quad,\quad(\forall)
t\in ]0,t_0[.
$$
Soit $t\in I_a$ arbitrairement fix\'e et $n\in {\bf N}$ tel que
$\frac{t}{n}\in ]0,t_0[$. Alors:
\begin{eqnarray*}
& &T(t)-S(t)={\left[T\left(n\frac{t}{n}\right)\right]}-{\left[S\left(n\frac{t}{n}\right)\right]}=\\
&=&T\left(n\frac{t}{n}\right)S\left(0\frac{t}{n}\right)-T\left((n-1)\frac{t}{n}\right)S\left(1\frac{t}{n}\right)+\\
&+&T\left((n-1)\frac{t}{n}\right)S\left(1\frac{t}{n}\right)-T\left((n-2)\frac{t}{n}\right)S\left(2\frac{t}{n}\right)+\\
&+&T\left((n-2)\frac{t}{n}\right)S\left(2\frac{t}{n}\right)-\cdots-T\left(0\frac{t}{n}\right)S\left(n\frac{t}{n}\right)=\\
&=&\sum\limits_{k=0}^{n-1}\left[T\left((n-k)\frac{t}{n}\right)S\left(k\frac{t}{n}\right)-T\left((n-k-1)\frac{t}{n}\right)S\left((k+1)\frac{t}{n}\right)\right]=\\
&=&\sum\limits_{k=0}^{n-1}T\left((n-k-1)\frac{t}{n}\right)\left[T\left(\frac{t}{n}\right)-S\left(\frac{t}{n}\right)\right]S\left(k\frac{t}{n}\right)
\end{eqnarray*}
quel que soit $t\in I_a.$\\
De l'in\'egalit\'e:
$$
\left\|\frac{T\left(\frac{t}{n}\right)-S\left(\frac{t}{n}\right)}{\frac{t}{n}}\right\|\leq\frac{\varepsilon}{ac_a^2}\quad,
$$
nous obtenons:
$$
\left\|T\left(\frac{t}{n}\right)-S\left(\frac{t}{n}\right)\right\|\leq\frac{\varepsilon}{ac_a^2}\frac{t}{n}
$$
et par suite:
$$
\left\|T(t)-S(t)\right\|\leq\sum\limits_{k=0}^{n-1}c_a\frac{\varepsilon}{ac_a^2}\frac{t}{n}c_a<\varepsilon\quad,\quad(\forall) t\in I_a.
$$
Puisque $\varepsilon >0$ est arbitraire, il en r\'esulte que
$T(t)=S(t)$, pour tout $t\in I_a$. Mais, comme $a\in ]0,\infty)$ est
aussi arbitraire, il s'ensuit que $T(t)=S(t)$, $(\forall) t\in
[0,\infty)$.\fin

Pr\'esentons maintenant la condition n\'ecessaire et suffisante pour
qu'un op\'erateur soit le g\'en\'erateur infinit\'esimal d'un
semi-groupe uniform\'ement continu.

\begin{teo}\label{num42}
Un op\'erateur $A:{\cal E}\longrightarrow{\cal E}$ est le
g\'en\'erateur infinit\'esimal d'un semi-groupe uniform\'ement continu si et
seulement si $A$ est un op\'erateur lin\'eaire born\'e.
\end{teo}
\dem
$\Longrightarrow$  Soit $A:{\cal E}\longrightarrow{\cal E}$ le
g\'en\'erateur infinit\'esimal d'un semi-groupe uniform\'ement continu $\left\{T(t)\right\}_{t\geq
0}\subset {\cal B(E)}$. Alors:
$$
\lim_{t\searrow 0}\left\|T(t)-I\right\|=0\quad.
$$
L'application $[0,\infty)\ni
t\mapsto T(t)\in {\cal B(E)}$ est continue et par suite
$\int\limits_{0}^{t}\!{T(s)}\:ds\in {\cal B(E)}$. Avec le lemme
\ref{num4}, on voit que:
$$
\lim_{t\searrow 0}\frac{1}{t}\int\limits_{0}^{t}\!T(s)\:ds=T(0)=I\quad.
$$
Il existe donc $\tau >0$ tel que:
$$
\left\|\frac{1}{\tau}\int\limits_{0}^{\tau}\!T(t)\:dt-I\right\|<1\quad.
$$
Compte tenu de la remarque \ref{num5}, l'\'el\'ement
$\frac{1}{\tau}\int\limits_{0}^{\tau}\!T(t)\:dt$ est
inversible, d'o\`u il s'ensuit
que $\int\limits_{0}^{\tau}\!T(t)\:dt$ est
inversible. Nous avons:
\begin{eqnarray*}
& &\frac{T(h)-I}{h}\int\limits_{0}^{\tau}\!T(t)\:dt=\frac{1}{h}\left[\int\limits_{0}^{\tau}\!T(t+h)\:dt-\int\limits_{0}^{\tau}\!T(t)\:dt\right]=\\
&=&\frac{1}{h}\int\limits_{\tau}^{\tau
+h}\!\!T(u)\:du-\frac{1}{h}\int\limits_{0}^{h}\!\!T(u)\:du\quad.
\end{eqnarray*}
Avec le lemme \ref{num4}, nous obtenons:
\begin{eqnarray*}
& &\lim_{h\searrow
0}\frac{T(h)-I}{h}\int\limits_{0}^{\tau}\!T(t)\:dt=\\
&=&\lim_{h\searrow
0}\left[\frac{1}{h}\int\limits_{\tau}^{\tau
+h}\!\!T(u)\:du-\frac{1}{h}\int\limits_{0}^{0+h}\!\!T(u)\:du\right]=\\
&=&T(\tau)-T(0)=T(\tau)-I\quad,
\end{eqnarray*}
d'o\`u:
$$
\lim_{h\searrow
0}\frac{T(h)-I}{h}=\left[T(\tau)-I\right]\left[\int\limits_{0}^{\tau}\!T(t)\:dt\right]^{-1}.
$$
Par cons\'equent, le g\'en\'erateur infinit\'esimal du semi-groupe
uniform\'ement continue $\left\{T(t)\right\}_{t\geq 0}$ est l'op\'erateur:
$$A=\left[T(\tau)-I\right]\left[\int\limits_{0}^{\tau}\!T(t)\:dt\right]^{-1}\in
{\cal B(E)}\quad.
$$
$\Longleftarrow$  Cette implication est \'evidente compte tenu du
lemme \ref{num6} et du lemme \ref{num7}.\fin

\newtheorem{cor}[lema]{Corollaire}
\begin{cor}\label{num8}
Soient $\left\{T(t)\right\}_{t\geq 0}$ un semi-groupe uniform\'ement continu et A son g\'en\'erateur
infinit\'esimal. Alors:\\
i) il existe $\omega\geq 0$ tel que $\left\|T(t)\right\|\leq
e^{\omega t}\quad,\quad(\forall) t\geq 0;$\\
ii) l'application $[0,\infty)\ni t \longmapsto T(t) \in {\cal B(E)}$ est
diff\'erentiable pour la topologie de la norme et:
$$
\frac{dT(t)}{dt}=AT(t)=T(t)A\quad,\quad(\forall) t\geq 0.
$$
\end{cor}
\dem
i) Nous avons:
$$
\left\|T(t)\right\|=\left\|e^{tA}\right\|\leq
e^{t\|A\|}\quad,\quad(\forall) t\geq 0.
$$
Pour $\omega=\|A\|$, nous obtenons l'in\'egalit\'e:
$$
\left\|T(t)\right\|\leq e^{\omega t}\quad,\quad(\forall) t\geq 0.
$$
L'assertion (ii) provient des \'egalit\'es suivantes:
$$
A=\lim_{t\searrow 0}\frac{T(t)-I}{t}=\lim_{t\searrow
0}\frac{T(t)-T(0)}{t-0}\quad,
$$
nous en d\'eduisons que l'application consid\'er\'ee est d\'erivable
au point $t=0$.\\
Soient $t>0$ et $h>0$. Alors:
\begin{eqnarray*}
& &\left\|\frac{T(t+h)-T(t)}{h}-AT(t)\right\|\leq\\
&\leq&\left\|\frac{T(h)-I}{h}-A\right\|\left\|T(t)\right\|\leq\left\|\frac{T(h)-I}{h}-A\right\|e^{t\|A\|}\quad,
\end{eqnarray*}
d'o\`u:
$$
\lim_{h\searrow 0}\left\|\frac{T(t+h)-T(t)}{h}-AT(t)\right\|=0\quad.
$$
Par cons\'equent, l'application consid\'er\'ee dans l'\'enonc\'e est
d\'erivable \`a droite et on a:
$$
\frac{d^{+}T(t)}{dt}=AT(t)\quad,\quad(\forall)t>0.
$$
Soient $t>0$ et $h<0$ tel que $t+h>0$. Alors:
\begin{eqnarray*}
& &\left\|\frac{T(t+h)-T(t)}{h}-AT(t)\right\|\leq\\
&\leq&\left\|\frac{I-T(-h)}{h}-AT(-h)\right\|\left\|T(t+h)\right\|\leq\\
&\leq&\left\|\frac{T(-h)-I}{-h}-AT(-h)\right\|e^{(t+h)\|A\|}\quad,
\end{eqnarray*}
d'o\`u il vient:
$$
\lim_{h\nearrow 0}\frac{T(t+h)-T(t)}{h}=AT(t)\quad.
$$
Par cons\'equent l'application consid\'er\'ee dans l'\'enonc\'e est d\'erivable \`a gauche et nous avons:
$$
\frac{d^{-}T(t)}{dt}=AT(t)\quad,\quad(\forall)t>0.
$$
Finalement on voit que l'application consid\'er\'ee dans l'\'enonc\'e est d\'erivable sur $[0,\infty)$ et
nous avons:
$$
\frac{dT(t)}{dt}=AT(t)\quad,\quad(\forall)t\geq 0.
$$
On v\'erifie que $AT(t)=T(t)A$ , $(\forall)t\geq
0$.\fin

Maintenant abordons quelques probl\`emes de th\'eorie spectrale pour
un semi-groupe uniform\'ement continu $\left\{T(t)\right\}_{t\geq 0}$ ayant
pour le g\'en\'erateur infinit\'esimal l'op\'erateur $A\in{\cal B(E)}$.

\begin{teo}
Soient $\left\{T(t)\right\}_{t\geq 0}$ un semi-groupe uniform\'ement continu et
$A$ son g\'en\'erateur infinit\'esimal. Si $\lambda\in{\bf C}$
tel que $\mbox{Re}\lambda>\|A\|$, alors l'application:
$$
R_{\lambda}:{\cal E}\longrightarrow{\cal E}\quad,
$$
$$
R_{\lambda}x=\int\limits_{0}^{\infty}\!e^{-\lambda t}T(t)x\:dt
$$
d\'efinit un op\'erateur lin\'eaire born\'e, $\lambda\in\rho(A)$ et
$R_{\lambda}x=R(\lambda;A)x$ , pour tout $x\in{\cal E}$.
\end{teo}
\dem
Soit $\lambda\in{\bf C}$ avec $\mbox{\em Re}\lambda>\|A\|$. Avec le
corollaire \ref{num8} (i), on voit que:
$$
\left\|T(t)\right\|\leq e^{\|A\|t}\quad,\quad(\forall)t\geq 0.
$$
De m\^eme, nous avons:
$$
\left\|e^{-\lambda t}T(t)x\right\|\leq
e^{-({\scriptstyle{Re}}\lambda-\|A\|)t}\|x\|\quad,\quad(\forall)x\in{\cal E},
$$
et:
$$
\int\limits_{0}^{\infty}\!e^{-({\scriptstyle{Re}}\lambda-\|A\|)t}\:dt=\frac{1}{\mbox{\em
Re}
\lambda-\|A\|}\quad.
$$
L'application $R_{\lambda}$ est donc born\'ee et il est clair que $R_{\lambda}$ est lin\'eaire.\\
Pour $x\in{\cal E}$, nous avons:
\begin{eqnarray*}
R_{\lambda}Ax&=&\int\limits_{0}^{\infty}\!e^{-\lambda t}T(t)Ax\:dt
=\int\limits_{0}^{\infty}\!e^{-\lambda t}\frac{d}{dt}T(t)x\:dt=\\
&=&-x+\lambda\int\limits_{0}^{\infty}\!e^{-\lambda
t}T(t)x\:dt=-x+\lambda R_{\lambda}x\quad,
\end{eqnarray*}
d'o\`u $x=R_{\lambda}(\lambda I-A)x$, pour tout $x\in{\cal E}$. Par
cons\'equent $R_{\lambda}(\lambda I-A)=I$.\\
De m\^eme, nous avons:
\begin{eqnarray*}
AR_{\lambda}x&=&A\int\limits_{0}^{\infty}\!e^{-\lambda
t}T(t)x\:dt=\int\limits_{0}^{\infty}\!e^{-\lambda t}AT(t)x\:dt=\\
&=&\int\limits_{0}^{\infty}\!e^{-\lambda
t}T(t)Ax\:dt=R_{\lambda}Ax\quad,\quad(\forall)x\in{\cal E}.
\end{eqnarray*}
Par suite, on a $AR_{\lambda}x=R_{\lambda}Ax=-x +\lambda R_{\lambda}x$, pour tout $x\in{\cal E}$. Il en  r\'esulte que $(\lambda
I-A)R_{\lambda}=I$.\\
Par cons\'equent $\lambda\in\rho(A)$ et $R_{\lambda}=R(\lambda;A)$.\fin

\begin{definitie}
L'op\'erateur $R_{\lambda}:{\cal E}\longrightarrow {\cal E}$ d\'efini
par:
$$
R_{\lambda}x=\int\limits_0^{\infty}\!e^{-\lambda
t}T(t)x\:dt\quad,\quad\lambda\in{\bf C}\mbox{ avec Re}\lambda>\|A\|,
$$
s'appelle la transform\'ee de Laplace du semi-groupe uniform\'ement continu
$\left\{T(t)\right\}_{t\geq 0}$ ayant pour g\'en\'erateur
infinit\'esimal l'op\'erateur $A$.
\end{definitie}

\begin{obs}
\em
On a:
$$
\left\{\lambda\in{\bf C}\left|\mbox{\em Re}\lambda
>\|A\|\right.\right\}\subset\rho(A)
$$
et:
$$
\sigma(A)\subset\left\{\lambda\in{\bf
C}\left|\mbox{\em Re}\lambda\leq\|A\|\right.\right\}\:.
$$
De m\^eme, nous obtenons:
$$
\left\|R(\lambda;A)\right\|\leq\frac{1}{\mbox{\em Re}\lambda-\|A\|}
$$
pour tout $\lambda\in{\bf C}$ avec $\mbox{\em Re}\lambda>\|A\|$.
\end{obs}

Pour obtenir des repr\'esentations de type Riesz-Dunford et de type Bromwich, on a besoin d'une classe sp\'eciale de contours de Jordan.

\begin{definitie}
Un contour de Jordan lisse et ferm\'e qui entoure $\sigma(A)$, s'appelle un contour de Jordan $A$-spectral s'il est homotope avec un cercle $C_r$ de centre $O$ et de rayon $r>\|A\|$.
\end{definitie}

\begin{teo}[Riesz-Dunford]\label{num9}
Soit $A$ le g\'en\'erateur infinit\'esimal d'un semi-groupe uniform\'ement continu $\left\{T(t)\right\}_{t\geq 0}$. Si $\Gamma_A$ est un contour
de Jordan $A$-spectral, alors nous avons:
$$
T(t)=\frac{1}{2\pi i}\int\limits_{\Gamma_A}\!e^{\lambda
t}R(\lambda;A)\:d\lambda\quad,\quad(\forall)t\geq 0.
$$
\end{teo}
\dem
Soit $\Gamma_A$ un contour de Jordan $A$-spectral. Alors $\Gamma_A$ est homotope avec un cercle $C_r$ de centre $O$ et de rayon $r>\|A\|$. Par cons\'equent, on a:
$$
\frac{1}{2\pi i}\int\limits_{\Gamma_A}\!e^{\lambda
t}R(\lambda;A)\:d\lambda=\frac{1}{2\pi i}\int\limits_{C_r}\!e^{\lambda
t}R(\lambda;A)\:d\lambda\quad,\quad(\forall)t\geq 0.
$$
Compt tenu de la proposition \ref{num3} (iii), on voit que:
$$
R(\lambda;A)=\sum\limits_{n=0}^{\infty}\frac{A^n}{\lambda^{n+1}}\quad,
$$
uniform\'ement par rapport \`a $\lambda$ sur les sous-ensembles compacts de $\{\lambda\in{\bf C}|\:|\lambda|>\|A\|\}$, particuli\`erement sur le cercle $C_r$. On a:
\begin{eqnarray*}
& &\frac{1}{2\pi i}\int\limits_{C_r}\!e^{\lambda
t}R(\lambda;A)\:d\lambda=\frac{1}{2\pi
i}\int\limits_{C_r}\!e^{\lambda
t}\sum\limits_{n=0}^{\infty}\frac{A^n}{\lambda^{n+1}}\:d\lambda=\\
&=&\sum\limits_{n=0}^{\infty}\frac{1}{2\pi
i}\int\limits_{C_r}\frac{e^{\lambda t}}{\lambda^{n+1}}\:d\lambda
A^n\quad.
\end{eqnarray*}
Appliquons la formule de Cauchy (\cite[pag. 228]{dunford-schwartz}) avec la fonction
$f(\lambda)=e^{\lambda t}$, nous obtenons:
$$
\frac{1}{2\pi i}\int\limits_{C_r}\frac{e^{\lambda
t}}{\lambda^{n+1}}\:d\lambda=\frac{t^n}{n!}\quad,\quad(\forall)n\in{\bf
N}\quad.
$$
Par cons\'equent:
$$
\frac{1}{2\pi i}\int\limits_{\Gamma_A}\!e^{\lambda
t}R(\lambda;A)\:d\lambda=\sum\limits_{n=0}^{\infty}\frac{t^nA^n}{n!}=e^{tA}=T(t)\quad,\quad(\forall)t\geq 0.\fin
$$

\begin{teo}[Bromwich]\label{num38}
Soient $\left\{T(t)\right\}_{t\geq 0}$ un semi-groupe uniform\'ement
continu et $A$ son g\'en\'erateur infinit\'esimal. Si $a>\|A\|$, alors nous avons:
$$
T(t)=\frac{1}{2\pi i}\int\limits_{a-i\infty}^{a+i\infty}\!e^{zt}R(z;A)\:dz
$$
et l'int\'egrale est uniform\'ement convergente par rapport \`a $t$
sur les intervalles compacts de $]0,\infty)$.
\end{teo}
\dem
Soit $a>\|A\|$, pour $R>2a$ nous consid\'erons le contour de Jordan lisse et ferm\'e
$$
\Gamma_R=\Gamma^{'}_R\cup\Gamma^{"}_R
$$
o\`u
$$
\Gamma^{'}_R=\left\{a+i\tau|\tau\in[-R,R]\right\}
$$
et
$$
\Gamma^{"}_R=\left\{a+R(\cos\varphi+i\sin\varphi)\left|\varphi\in\left[\frac{\pi}{2},\frac{3\pi}{2}\right]\right.\right\}\quad.
$$
Remarquons que pour $z\in\Gamma^{'}_R$ on a:
$$
|z|=|a+i\tau|>a>\|A\|\quad.
$$
De m\^eme, si $z\in\Gamma^{"}_R$, alors nous avons:
\begin{eqnarray*}
|z|&=&|a+(\cos\varphi+i\sin\varphi)|=|a-[-R(\cos\varphi+i\sin\varphi)]|\geq\\
&\geq&\left||a|-|-R(\cos\varphi+i\sin\varphi)|\right|=|a-R|=R-a>\|A\|\quad.
\end{eqnarray*}
Par cons\'equent, $z\in\Gamma_R$ implique $z\in\rho(A)$. De plus, on voit que $\Gamma_R$ est homotope au cercle $C$ de centre $O$ et de rayon $R-a$. Il s'ensuit donc que $\Gamma_R$ est un contour de Jordan $A$-spectral et avec le th\'eor\`eme de Riesz-Dunford nous obtenons:
$$
T(t)=\frac{1}{2\pi i}\int\limits_{\Gamma_R}\!e^{zt}R(z;A)\:dz\quad,\quad(\forall)t\geq 0,
$$
pour tout $R>2a$. Il en r\'esulte:
$$
T(t)=I^{'}_t(R)+I^{"}_t(R)\quad,\quad(\forall)t\geq 0,
$$
pour tout $R>2a$, o\`u nous avons not\'e
$$
I^{'}_t(R)=\frac{1}{2\pi i}\int\limits_{\Gamma^{'}_R}\!e^{zt}R(z;A)\:dz
$$
et
$$
I^{"}_t(R)=\frac{1}{2\pi i}\int\limits_{\Gamma^{"}_R}\!e^{zt}R(z;A)\:dz\quad.
$$
Montrons que
$$
\lim_{R\rightarrow\infty}\frac{1}{2\pi i}\int\limits_{\Gamma^{"}_R}\!e^{zt}R(z;A)\:dz=0\quad,\quad(\forall)t\geq 0.
$$
Compte tenu de la proposition \ref{num3} (iii), on voit que:
$$
R(z;A)=\sum\limits_{n=0}^{\infty}\frac{A^n}{z^{n+1}}\quad,
$$
la s\'erie de la partie droite de l'\'egalit\'e \'etant uniform\'ement convergente par rapport \`a $z$ sur les sous-ensembles compacts de $\{z\in{\bf C}|\:|z|>\|A\|\}$, particuli\`erement sur $\Gamma^{"}_R$. Il s'ensuit que:
\begin{eqnarray*}
& &I^{"}(R)=\sum\limits_{n=0}^{\infty}\left(\frac{1}{2\pi i}\int\limits_{\Gamma^{"}_R}\frac{e^{zt}}{z^{n+1}}A^n\:dz\right)=\\
&=&\left(\frac{1}{2\pi i}\int\limits_{\Gamma^{"}_R}\frac{e^{zt}}{z}\:dz\right)I+\sum\limits_{n=1}^{\infty}\left(\frac{1}{2\pi i}\int\limits_{\Gamma^{"}_R}\frac{e^{zt}}{z^{n+1}}\:dz\right)A^n\quad,\quad(\forall)t\geq 0,
\end{eqnarray*}
pour tout $R>2a$.
Notons
$$
A_t(R)=\left(\frac{1}{2\pi i}\int\limits_{\Gamma^{"}_R}\frac{e^{zt}}{z}\:dz\right)I
$$
et
$$
B_t(R)=\sum\limits_{n=1}^{\infty}\left(\frac{1}{2\pi i}\int\limits_{\Gamma^{"}_R}\frac{e^{zt}}{z^{n+1}}\:dz\right)A^n\quad.
$$
Pour l'int\'egrale $A_t(R)$, avec la param\'etrisation suivante
$$
z=a+R(\cos\varphi+i\sin\varphi)\quad,\quad\varphi\in\left[\frac{\pi}{2},\frac{3\pi}{2}\right],
$$
on obtient:
\begin{eqnarray*}
A_t(R)&=&\left[\frac{1}{2\pi i}\int\limits_{\frac{\pi}{2}}^{\frac{3\pi}{2}}\frac{e^{t(a+R\cos\varphi+i\sin\varphi)}}{z}R(-\sin\varphi+i\cos\varphi)\:d\varphi\right]I=\\
&=&\left[\frac{R}{2\pi}e^{ta}\int\limits_{\frac{\pi}{2}}^{\frac{3\pi}{2}}e^{tR\cos\varphi}e^{itR\sin\varphi}\frac{1}{z}(\cos\varphi+i\sin\varphi)\:d\varphi\right]I\quad.
\end{eqnarray*}
Il en r\'esulte que:
\begin{eqnarray*}
\left\|A_t(R)\right\|&\leq&\frac{R}{2\pi}e^{ta}\int\limits_{\frac{\pi}{2}}^{\frac{3\pi}{2}}\left|e^{tR\cos\varphi}\right|\:\left|e^{itR\sin\varphi}\right|\frac{1}{|z|}|\cos\varphi+i\sin\varphi|\:d\varphi\leq\\
&\leq&\frac{R}{2\pi}e^{ta}\int\limits_{\frac{\pi}{2}}^{\frac{3\pi}{2}}e^{tR\cos\varphi}\frac{1}{R-a}\:d\varphi=\\
&=&\frac{1}{2\pi}\frac{R}{R-a}e^{ta}\int\limits_{\frac{\pi}{2}}^{\frac{3\pi}{2}}e^{tR\cos\varphi}\:d\varphi
\end{eqnarray*}
parce que $z\in\Gamma^{"}_R$ implique
$$
|z|=|a+R(\cos\varphi+i\sin\varphi)|>R-a
$$
donc
$$
\frac{1}{|z|}<\frac{1}{R-a}\quad.
$$
De l'in\'egalit\'e $R>2a$, on obtient $2R-2a>R$, d'o\`u
$$
\frac{R}{R-a}<2\quad.
$$
Par cons\'equent:
$$
\left\|A_t(R)\right\|\leq\frac{1}{\pi}e^{ta}\int\limits_{\frac{\pi}{2}}^{\frac{3\pi}{2}}e^{tR\cos\varphi}\:d\varphi\quad,\quad(\forall)t\geq 0,
$$
pour tout $R>2a$.
Soient $0\leq t_1<t_2$ et $t\in[t_1,t_2]$. Pour tout $R>2a$ et tout $\varphi\in\left[\frac{\pi}{2},\frac{3\pi}{2}\right]$ on a 
$$
e^{tR\cos\varphi}\leq 1\quad.
$$
Comme
$$
\lim_{R\rightarrow\infty}e^{tR\cos\varphi}=0\quad,
$$
avec le th\'eor\`eme de la convergence born\'ee de Lebesgue il r\'esulte que
$$
\lim_{R\rightarrow\infty}\int\limits_{\frac{\pi}{2}}^{\frac{3\pi}{2}}e^{tR\cos\varphi}\:d\varphi=0
$$
et par cons\'equent
$$
\lim_{R\rightarrow\infty}A_t(R)=0
$$
uniform\'ement par rapport \`a $t\in[t_1,t_2]$.

Soit maintenant l'int\'egrale
$$
B_t(R)=\sum\limits_{n=1}^{\infty}\left(\frac{1}{2\pi i}\int\limits_{\Gamma^{"}_R}\frac{e^{zt}}{z^{n+1}}\:dz\right)A^n\quad.
$$
Pour tout $t\in[t_1,t_2]$ et tout $R>2a$ on a:
$$
e^{tR\cos\varphi}\leq 1\quad,\quad(\forall)\varphi\in\left[\frac{\pi}{2},\frac{3\pi}{2}\right]\quad.
$$
On voit que:
$$
\left|\int\limits_{\Gamma^{"}_R}\frac{e^{zt}}{z^{n+1}}\:dz\right|\leq\frac{Re^{ta}}{\left(R-a\right)^{n+1}}\int\limits_{\frac{\pi}{2}}^{\frac{3\pi}{2}}e^{tR\cos\varphi}\:d\varphi\leq\pi e^{ta}\frac{R}{(R-a)^{n+1}}\quad.
$$
Puisque $R>2a>a+\|A\|$, il vient:
$$
\left\|B_t(R)\right\|\leq\sum\limits_{n=1}^{\infty}\frac{\|A\|^n}{2\pi}\left|\int\limits_{\Gamma^{"}_R}\frac{e^{zt}}{z^{n+1}}\:dz\right|\leq\frac{e^{ta}}{2}\frac{R}{R-a}\sum\limits_{n=1}^{\infty}\left(\frac{\|A\|}{R-a}\right)^n
$$
et comme
$$
\frac{\|A\|}{R-a}<1\quad,
$$
il en r\'esulte que:
$$
\left\|B_t(R)\right\|\leq e^{ta}\;\frac{\|A\|}{2}\;\frac{R}{R-a}\;\frac{1}{R-a-\|A\|}\quad,
$$
quel que soit $R>2a$.
Donc
$$
\lim_{R\rightarrow\infty}B_t(R)=0\quad,
$$
uniform\'ement par rapport \`a $t\in[t_1,t_2]$. Il s'ensuit donc que
$$
\lim_{R\rightarrow\infty}I^{"}_t(R)=0\quad,
$$
uniform\'ement par rapport \`a $t\in[t_1,t_2]$.\\
Par cons\'equent:
$$
T(t)=\lim_{R\rightarrow\infty}\frac{1}{2\pi i}\int\limits_{\Gamma^{"}_R}\!e^{zt}R(z;A)\:dz=\frac{1}{2\pi i}\int\limits_{Re\:\lambda-i\infty}^{Re\:\lambda+i\infty}\!e^{zt}R(z;A)\:dz\quad,
$$
uniform\'ement par rapport \`a $t$ sur les intervalles compacts de $]0,\infty)$.\fin

Nous finissons cette section avec le th\'eor\`eme spectral pour les
semi-groupes uniform\'ement continus.

\begin{teo}[spectral mapping]\label{num49}
Soit $A$ le g\'en\'erateur infinit\'esimal du semi-groupe uniform\'ement continu $\left\{T(t)\right\}_{t\geq 0}$. Alors:
$$
e^{t\sigma(A)}=\sigma\left(T(t)\right)\quad,\quad(\forall)t\geq 0.
$$
\end{teo}
\dem
Montrons que
$e^{t\sigma(A)}\subset\sigma\left(T(t)\right)$ , $(\forall)t\geq 0$.\\
Soit $\xi\in\sigma(A)$. Pour $\lambda\in\rho(A)$, l'application:
$$
g_\xi (\lambda)=\frac{e^{\xi t}-e^{\lambda t}}{\xi
-\lambda}
$$
est analytique dans un voisinage de $\sigma(A)$.
Compte tenu du th\'eor\`eme \ref{num9}, on voit que:
$$
e^{\xi t}I-e^{At}=(\xi I-A)g_\xi(A)\quad.
$$
Si $e^{\xi t}\in\rho\left(T(t)\right)$, alors il existe
$Q=\left[e^{\xi t}I-T(t)\right]^{-1}\in{\cal B(E)}$. Par
cons\'equent:
$$
I=(\xi I-A)g_\xi (A)Q\quad,
$$
d'o\`u il r\'esulte que $\xi\in\rho(A)$, ce qui est absurde.
Donc $e^{\xi t}\in\sigma\left(T(t)\right)$ et par suite
$e^{t\sigma(A)}\subset\sigma\left(T(t)\right)$.\\
Montrons que $\sigma\left(T(t)\right)\subset e^{t\sigma(A)}$.\\
Soit $\mu\in\sigma\left(T(t)\right)$. Supposons par absurde que
$\mu\not\in e^{t\sigma(A)}$. Alors pour $\lambda\in\rho(A)$, l'application:
$$
h(\lambda)=\left(\mu- e^{\lambda t}\right)^{-1}
$$
est d\'efinie sur un voisinage du $\sigma(A)$. Donc:
$$
h(A)\left(\mu I-e^{tA}\right)=I
$$
et il en r\'esulte que $\mu\in\rho\left(T(t)\right)$ et cela est
absurde. Par suite $\mu\in e^{t\sigma(A)}$, d'o\`u
$\sigma\left(T(t)\right)\subset e^{t\sigma(A)}$. Finalement on voit
que:
$$
e^{t\sigma(A)}=\sigma\left(T(t)\right)\quad,\quad(\forall)t\geq0\quad.\fin
$$
\vspace{2cm}

\section{Notes}

{\footnotesize
Les notions pr\'es\'ent\'ees dans cet chapitre se trouvent en majorit\'e des travaux concernant les semi-groupes d'op\'erateurs
lin\'eaires. Pour les propri\'et\'es de la pseudo-r\'esolvante, on
peut consulter \cite[pag. 36]{pazy1}.

De m\^eme, on peut trouver les
op\'erateurs dissipatifs dans \cite[pag. 13]{pazy1},
\cite[pag. 52]{davies} et \cite[pag. 30]{ahmed}. Une jolie
g\'en\'eralisation pour ces op\'erateurs est donn\'ee dans
\cite[pag. 61]{clement}.

Le th\'eor\`eme \ref{num42} a \'et\'e montr\'e pour la premi\`ere fois
ind\'ependemment par Yosida dans \cite{yosida1} et par Nathan dans
\cite{nathan}. Nous avons consult\'e aussi les preuves donn\'ees par
Pazy dans \cite[pag. 2]{pazy1}, Ahmed dans \cite[pag. 4]{ahmed} et
Davies dans \cite[pag. 19]{davies}. Compte tenu du ce th\'eor\`eme, on
peut introduire la transform\'ee de Laplace pour un semi-groupe
uniform\'ement continu et on peut montrer le th\'eor\`eme \ref{num9} et
le th\'eor\`eme \ref{num49} comme des applications du calcul
fonctionnel de Dunford (\cite[pag. 568]{dunford-schwartz}). Pour le
th\'eor\`eme \ref{num38} on peut consulter \cite[pag. 25]{pazy1}.
}

  \chapter{Semi-groupes de classe $C_0$}

    \section{D\'efinitions. Propri\'et\'es \'el\'ementaires}

\hspace{1cm}Dans le cadre de ce paragraphe, nous introduisons une classe plus
    g\'en\'erale que la classe des semi-groupes uniform\'ement continus et nous
    \'etudions leurs propri\'et\'es \'el\'ementaires.

\begin{definitie}\label{num10}
On appelle $C_0$-semi-groupe (ou semi-groupe fortement continu)
d'op\'erateurs lin\'eaires born\'es sur $\cal E$ une famille
$\left\{T(t)\right\}_{t\geq 0}\subset{\cal B(E)}$ v\'erifiant les
propri\'et\'es suivantes:\\
i) $T(0)=I;$\\
ii) $T(t+s)=T(t)T(s)\quad,\quad(\forall)t,s\geq 0;$\\
iii) $\lim_{t\searrow 0}T(t)x=x\quad,\quad(\forall)x\in{\cal E}.$
\end{definitie}

\begin{definitie}
On appelle g\'en\'erateur infinit\'esimal d'un $C_0$-semi-groupe
$\left\{T(t)\right\}_{t\geq 0}$, un op\'erateur $A$ d\'efini sur
l'ensemble:
$$
{\cal D}(A)=\left\{x\in {\cal E}\left|\lim_{t\searrow 0}\frac{T(t)x-x}{t}\mbox{ existe
}\right.\right\}
$$
par:
$$
Ax=\lim_{t\searrow 0}\frac{T(t)x-x}{t}\quad,\quad(\forall)x\in{\cal D}(A).
$$
\end{definitie}

\begin{obs}
\em
Il est clair que le g\'en\'erateur infinit\'esimal d'un $C_0$-semi-groupe
est un op\'erateur lin\'eaire.
\end{obs}

\begin{obs}
\em
Puisque:
$$
\left\|T(t)x-x\right\|\leq\left\|T(t)-I\right\|\|x\|
$$
pour tout $x\in{\cal E}$ et tout $t\geq 0$, il en r\'esulte que les
semi-groupes uniform\'ement continus sont $C_0$-semi-groupes. Mais il existe des $C_0$-semi-groupes qui ne sont pas uniform\'ement continus, comme nous pouvons le voir dans les exemples suivants.
\end{obs}

\begin{exemplu}
\em
Soit:
$$
{\cal C}[0,\infty)=\left\{\left.f:[0,\infty)\rightarrow {\bf R}\right|f\;\mbox{\em
est uniform\'ement continue et born\'ee}\right\}\quad.
$$
Avec la norme
$\|f\|_{{\cal C}[0,\infty)}=\sup_{\alpha\in[0,\infty)}|f(\alpha)|$, l'espace ${\cal
C}[0,\infty)$ devient un espace de Banach. D\'efinissons:
$$
\left(T(t)f\right)(\alpha)=f(t+\alpha)\quad,\quad(\forall)t\geq
0\mbox{ et }\alpha\in[0,\infty).
$$
Evidemment $T(t)$ est un op\'erateur lin\'eaire, et, en plus, on
a:\\
i) $\left(T(0)f\right)(\alpha)=f(0+\alpha)=f(\alpha)$. Donc $T(0)=I$;\\
ii)
$\left(T(t+s)f\right)(\alpha)=f(t+s+\alpha)=\left(T(t)f\right)(s+\alpha)=\left(T(t)T(s)f\right)(\alpha)$,
$(\forall)f\in{\cal C}[0,\infty)$. Donc $T(t+s)=T(t)T(s)$,
$(\forall)t,s\geq 0$;\\
iii) $\lim_{t\searrow 0}\left\|T(t)f-f\right\|_{{\cal C}[0,\infty)}=\lim_{t\searrow
0}\left\{\sup_{\alpha\in[0,\infty)}\left|f(t+\alpha)-f(\alpha)\right|\right\}=0$,
$(\forall)f\in{\cal C}[0,\infty)$.\\
De m\^eme, nous avons:
\begin{eqnarray*}
& &\left\|T(t)f\right\|_{{\cal C}[0,\infty)}=\sup_{\alpha\in[0,\infty)}\left|\left(T(t)f\right)(\alpha)\right|=\sup_{\alpha\in[0,\infty)}\left|f(t+\alpha)\right|=\\
&=&\sup_{\beta\in[t,\infty)}\left|f(\beta)\right|\leq\sup_{\beta\in[0,\infty)}\left|f(\beta)\right|=\|f\|_{{\cal
C}[0,\infty)}
\quad,\quad(\forall)t\geq
0.
\end{eqnarray*}
Donc $\left\|T(t)\right\|=1$, $(\forall)t\geq 0$.
Par cons\'equent $\left\{T(t)\right\}_{t\geq 0}$ est un $C_0$-semi-groupe
d'op\'erateurs lin\'eaires born\'es sur ${\cal C}[0,\infty)$, nomm\'e le $C_0$-semi-groupe de
translation \`a droite.\\
Soit $A:{\cal D}(A)\subset{\cal C}[0,\infty)\longrightarrow{\cal C}[0,\infty)$ le
g\'en\'erateur infinit\'esimal du $C_0$-semi-groupe
$\left\{T(t)\right\}_{t\geq 0}$. Si $f\in{\cal D}(A)$, alors nous avons:
$$
Af(\alpha)=\lim_{t\searrow
0}\frac{T(t)f(\alpha)-f(\alpha)}{t}=\lim_{t\searrow
0}\frac{f(\alpha+t)-f(\alpha)}{t}=f'(\alpha)\quad,
$$
uniform\'ement par rapport \`a $\alpha$.
Par cons\'equent:
$$
{\cal D}(A)\subset\left\{f\in{\cal C}[0,\infty)\left|f'\in{\cal
C}[0,\infty)\right.\right\}\quad.
$$
Si $f\in{\cal C}[0,\infty)$ tel que $f'\in{\cal
C}[0,\infty)$, alors:
$$
\left\|\frac{T(t)f-f}{t}-f'\right\|_{{\cal C}[0,\infty)}=\sup_{\alpha\in[0,\infty)}\left|\frac{\left(T(t)f\right)(\alpha)-f(\alpha)}{t}-f'(\alpha)\right|\quad.
$$
Mais:
\begin{eqnarray*}
& &\left|\frac{\left(T(t)f\right)(\alpha)-f(\alpha)}{t}-f'(\alpha)\right|=\left|\frac{f(\alpha+t)-f(\alpha)}{t}-f'(\alpha)\right|=\\
&=&\left|\frac{1}{t}f(\tau)|_\alpha^{\alpha+t}-f'(\alpha)\right|=\frac{1}{t}\left|\int\limits_{\alpha}^{\alpha+t}\!\![f'(\tau)-f'(\alpha)]\:d\tau\right|\leq\\
&\leq&\frac{1}{t}\int\limits_{\alpha}^{\alpha+t}\!\left|f'(\tau)-f'(\alpha)\right|\:d\tau\longrightarrow
0
\end{eqnarray*}
uniform\'ement par rapport \`a $\alpha$ pour $t\searrow 0$. Par suite:
$$
\left\|\frac{T(t)f-f}{t}-f'\right\|_{{\cal
C}[0,\infty)}\longrightarrow 0\quad\mbox{si}\quad t\searrow 0,
$$
d'o\`u
$f\in{\cal D}(A)$ et:
$$
\left\{f\in{\cal C}[0,\infty)\left|f'\in{\cal
C}[0,\infty)\right.\right\}\subset{\cal D}(A)\quad.
$$
Par cons\'equent
${\cal D}(A)=\left\{f\in{\cal C}[0,\infty)\left|f'\in{\cal C}[0,\infty)\right.\right\}$
et $Af=f'$. Comme cet op\'erateur est non born\'e,  compte tenu du
th\'eor\`eme \ref{num42}, il ne peut pas engendrer un semi-groupe
uniform\'ement continu.
\end{exemplu}

\begin{exemplu}
\em
Consid\'erons l'espace $L_p]0,\infty)$, $1\leq p<\infty$, avec la
norme:
$$
\|f\|_p=\left\{\int\limits_{0}^{\infty}\!\left|f(\alpha)\right|^p\:d\alpha\right\}^{\frac{1}{p}}\quad.
$$
Avec cette norme, $L_p]0,\infty)$, $1\leq p<\infty$, est un espace de
Banach. D\'efinissons:
$$
\left(T(t)f\right)(\alpha)=f(t+\alpha)\quad,\quad(\forall)t\geq
0\mbox{ et }\alpha\in]0,\infty).
$$
Nous avons:
\begin{eqnarray*}
\left\|T(t)f\right\|_p&=&\left\{\int\limits_{0}^{\infty}\!\left|\left(T(t)f\right)(\alpha)\right|^p\:d\alpha\right\}^{\frac{1}{p}}=\left\{\int\limits_{0}^{\infty}\!\left|f(\alpha+t)\right|^p\:d\alpha\right\}^{\frac{1}{p}}=\\
&=&\left\{\int\limits_{t}^{\infty}\!\left|f(\beta)\right|^p\:d\beta\right\}^{\frac{1}{p}}\leq\left\{\int\limits_{0}^{\infty}\!\left|f(\beta)\right|^p\:d\beta\right\}^{\frac{1}{p}}=\|f\|_p\quad.
\end{eqnarray*}
Donc $\left\|T(t)\right\|=1$, $(\forall)t\geq 0$.\\
Il est \'evident que $T(0)=I$ et $T(t+s)=T(t)T(s)$, $(\forall)t,s\geq
0$.\\
De plus, on a:
\begin{eqnarray*}
& &\lim_{t\searrow 0}\left\|T(t)f-f\right\|_p=\lim_{t\searrow
0}\left\{\int\limits_{0}^{\infty}\!\left|\left(T(t)f\right)(\alpha)-f(\alpha)\right|^p\:d\alpha\right\}^{\frac{1}{p}}=\\
&=&\lim_{t\searrow
0}\left\{\int\limits_{0}^{\infty}\!\left|f(\alpha+t)-f(\alpha)\right|^p\:d\alpha\right\}^{\frac{1}{p}}=0\quad.
\end{eqnarray*}
Par suite $\left\{T(t)\right\}_{t\geq 0}$ est un $C_0$-semi-groupe
d'op\'erateurs lin\'eaires born\'es sur $L_p]0,\infty)$.\\
Soit $A:{\cal D}(A)\subset L_p]0,\infty)\longrightarrow L_p]0,\infty)$
le g\'en\'erateur infinit\'esimal du $C_0$-semi-groupe
$\left\{T(t)\right\}_{t\geq 0}$. Si $f\in{\cal D}(A)$, alors nous avons:
$$
Af(\alpha)=\lim_{t\searrow
0}\frac{T(t)f(\alpha)-f(\alpha)}{t}=\lim_{t\searrow
0}\frac{f(\alpha+t)-f(\alpha)}{t}=f'(\alpha)
$$
uniform\'ement par rapport \`a $\alpha$.
Par cons\'equent:
$$
{\cal D}(A)\subset\left\{f\in L_p]0,\infty)\left|f'\in L_p]0,\infty)\right.\right\}\quad.
$$
Si $f\in L_p]0,\infty)$ tel que $f'\in L_p]0,\infty)$, alors on a:
$$
\left\|\frac{T(t)f-f}{t}-f'\right\|_p=\left\{\int\limits_{0}^{\infty}\!\left|\frac{\left(T(t)f\right)(\alpha)-f(\alpha)}{t}-f'(\alpha)\right|^p\:d\alpha\right\}^{\frac{1}{p}}.
$$
Mais:
\begin{eqnarray*}
&
&\left|\frac{\left(T(t)f\right)(\alpha)-f(\alpha)}{t}-f'(\alpha)\right|=\left|\frac{f(\alpha+t)-f(\alpha)}{t}-f'(\alpha)\right|=\\
&=&\left|\left.\left[\frac{1}{t}f(\tau)\right]\right|_{\alpha}^{\alpha+t}-\left.\left[\frac{1}{t}f'(\alpha)\tau\right]\right|_{\alpha}^{\alpha+t}\right|=\left|\frac{1}{t}\int\limits_{\alpha}^{\alpha+t}\!\left[f'(\tau)-f'(\alpha)\right]\:d\tau\right|\longrightarrow
0
\end{eqnarray*}
uniform\'ement par rapport \`a $\alpha$ si $t\searrow 0$. Alors:
$$
\left\|\frac{T(t)f-f}{t}-f'\right\|_p\longrightarrow
0\quad\mbox{si}\quad t\searrow 0
$$
et on voit que:
$$
\left\{f\in L_p]0,\infty)\left|f'\in
L_p]0,\infty)\right.\right\}\subset{\cal D}(A)\quad.
$$
Par cons\'equent:
$$
{\cal D}(A)=\left\{f\in L_p]0,\infty)\left|f'\in L_p]0,\infty)\right.\right\}
$$
et $Af=f'$.
\end{exemplu}

\begin{teo}\label{num50}
Soit $\left\{T(t)\right\}_{t\geq 0}\subset{\cal B(E)}$
une famille ayant les propri\'et\'es:\\
i) $T(0)=I;$\\
ii) $T(t+s)=T(t)T(s)\quad,\quad(\forall)t,s\geq 0.$\\
Les affirmations suivantes sont \'equivalentes:\\
iii') $\lim_{t\searrow 0}T(t)=I$ dans la topologie forte;\\
iii'') $\lim_{t\searrow 0}T(t)=I$ dans la topologie faible.
\end{teo}
\dem
$iii') \Longrightarrow iii'')$ Cette implication est \'evidente.\\
$iii'') \Longrightarrow iii')$ Supposons que:
$$
\lim_{t\searrow 0}T(t)=I
$$
dans la topologie faible. Alors, pour tout $x\in{\cal E}$ et
tout $x^*\in{\cal E}^*$ on a:
$$
\lim_{t\searrow 0}\langle T(t)x,x^*\rangle=\langle x,x^*\rangle\quad.
$$
Si $t_0>0$, alors pour tout $h>0$, nous obtenons:
\begin{eqnarray*}
& &\left|\langle T(t_0+h)x,x^*\rangle-\langle T(t_0)x,x^*\rangle\right|=\\
&=&\left|\langle T(t_0)T(h)x,x^*\rangle-\langle
T(t_0)x,x^*\rangle\right|=\\
&=&\left|\langle T(t_0)[T(h)x-x],x^*\rangle\right|\longrightarrow
0\quad\mbox{si}\quad h\searrow 0,
\end{eqnarray*}
quel que soit $x\in{\cal E}$ et $x^*\in{\cal E}^*$. Par suite,
l'application:
$$
[0,\infty)\ni t\longmapsto T(t)\in{\cal B(E)}
$$
est faiblement continue \`a droite sur $[0,\infty)$ et on voit qu'elle
est faiblement continue sur $]0,\infty)$. En particulier, elle est
faiblement mesurable sur $]0,\infty)$.
Pour $x\in{\cal E}$ arbitrairement fix\'e, consid\'erons
l'application:
$$
[0,\infty)\ni t\longmapsto T(t)x\in{\cal E}
$$
et d\'esignons par:
$$
{\cal I}\mbox{\em m }T(\:.\:)x=\{T(t)x|t\in[0,\infty)\}
$$
son image. Supposons que l'ensemble:
$$
{\cal K}_x=\{T(q)x|q\in{\bf Q}_+^*\}\subset{\cal
I}\mbox{\em m }T(\:.\:)x
$$
n'est pas dense dans ${\cal I}\mbox{\em m }T(\:.\:)x$. Alors, il
existe $t_0\in[0,\infty)$ tel que $T(t_0)x\in{\cal I}\mbox{\em m
}T(\:.\:)x$ et:
$$
d\left(T(t_0)x,{\cal K}_x\right)>0\quad.
$$
En appliquant un corollaire du th\'eor\`eme de Hahn-Banach
(\cite[Corollary II.3.13, pag. 64]{dunford-schwartz}), on
voit qu'il existe $x_0^*\in{\cal E}^*$ tel que:
$$
\langle k_n,x_0^*\rangle=0\quad,\quad(\forall)k_n\in{\cal K}_x
$$
et:
$$
\langle T(t_0)x,x_0^*\rangle =1\quad.
$$
Soit $t_{n}\in{\bf Q}_+^*$ tel que
$\lim_{n\rightarrow\infty}t_{n}=t_0$. Alors, compte tenu de la
continuit\'e faible de l'application consid\'er\'ee, il vient:
$$
0=\lim_{n\rightarrow\infty}\langle
T\left(t_{n}\right)x,x_0^*\rangle=\langle T(t_0)x,x_0^*\rangle=1\quad,
$$
ce qui est absurde. Il s'ensuit que:
$$
\overline{{\cal K}_x}={\cal I}\mbox{\em m }T(\:.\:)x\quad,
$$
pour tout $x\in{\cal E}$. Par cons\'equent, l'application
consid\'er\'ee a une image s\'eparable. En appliquant le th\'eor\`eme de
Pettis (\cite[Theorem 3.2.2, pag. 36]{hille}), il vient que cette application est fortement mesurable sur
$]0,\infty)$. Alors, il r\'esulte que pour tout $x_n\in{\cal E}$ avec
$\|x\|\leq 1$, l'application:
$$
\|T(\:.\:)\|=\sup_{n\in{\bf N}}\|T(\:.\:)x_n\|<\infty
$$
est mesurable sur $]0,\infty)$. Montrons que l'application
$\|T(\:.\:)\|$ est born\'ee sur les intevalles
$[\alpha,\beta]\subset]0,\infty)$. Compte tenu du th\'eor\`eme de
Banach-Steinhaus (\cite[Theorem II.1.11, pag. 52]{dunford-schwartz}), il est suffisant de montrer que $\|T(\:.\:)x\|$ est
born\'ee sur les intervalles $[\alpha,\beta]$, pour tout $x\in{\cal E}$. Soient
$\alpha,\beta\in]0,\infty)$. Supposons qu'il existe
$x_0\in{\cal E}$ tel que pour tout $M>0$ on puisse trouver
$s\in[\alpha;\beta]$ tel que:
$$
\|T(s)x_0\|>M\quad.
$$
Donc il existe $t_n\in[\alpha,\beta]$, $n\in{\bf N}$, tel que:
$$
\lim_{n\rightarrow\infty}t_n=\tau\in[\alpha,\beta]
$$
et:
$$
\|T(t_n)x_0\|>n\quad,\quad(\forall)n\in{\bf N}.
$$
D'autre part, l'application $\|T(\:.\:)x_0\|$
est mesurable sur $]0,\infty)$. Donc il existe une constante $K>0$ et
un ensemble mesurable ${\bf F}\subset[0,\tau]$ avec $m({\bf
F})>\frac{\tau}{2}$ tel que:
$$
\sup_{t\in{\bf F}}\|T(t)x_0\|\leq K\quad.
$$
Si nous consid\'erons:
$$
{\bf E}_n=\{t_n-\eta|\eta\in{\bf F}\cap[0,t_n]\}\quad,
$$
on voit que ${\bf E}_n$ est un ensemble mesurable et pour $n$
suffisamment grand, nous obtenons:
$$
m({\bf E}_n)\geq\frac{\tau}{2}\quad.
$$
Alors, pour tout $\eta\in{\bf F}\cap[0,t_n]$, $n\in{\bf N}$, nous
avons:
$$
n\leq\|T(t_n)x_0\|\leq\|T(t_n-\eta)\|\:\|T(\eta)x_0\|\leq\|T(t_n-\eta)\|K\quad,
$$
d'o\`u:
$$
\|T(t)\|\geq\frac{n}{K}\quad,\quad(\forall)t\in{\bf E}_n.
$$
Si nous notons:
$$
{\bf E}=\lim\sup_{n\in{\bf N}}{\bf E}_n=\bigcap_{n\geq 0}\bigcup_{k\geq n}
{\bf E}_k\quad,
$$
alors on voit que:
$$
m({\bf E})\geq\frac{\tau}{2}
$$
et:
$$
\|T(t)\|=\infty\quad,\quad(\forall)t\in{\bf E}
$$
ce qui est absurde. Par cons\'equent, il existe $M>0$ tel que:
$$
\|T(t)\|\leq M\quad,\quad(\forall)t\in[\alpha,\beta].
$$
Soient $\alpha,\beta,t,t_0\in]0,\infty)$ tel que:
$$
0<\alpha<t<\beta<t_0
$$
et $\varepsilon>0$ tel que $\beta<t_0-\varepsilon$. Alors pour tout
$x\in{\cal E}$, l'application:
$$
[\alpha,\beta]\ni t\longmapsto T(t_0)x=T(t)T(t_0-t)x\in{\cal E}
$$
ne d\'epend pas de $t$, donc elle est B\^ochner int\'egrable par rapport
\`a $t\in[\alpha,\beta]$ et pour tout $x\in{\cal E}$ on a:
\begin{eqnarray*}
& &(\beta-\alpha)\left[T(t_0\pm\varepsilon)x-T(t_0)x\right]\:dt=\\
&=&\int\limits_{\alpha}^{\beta}\!T(t)\left[T(t_0\pm\varepsilon-t)x-T(t_0-t)x\right]\:dt\quad,
\end{eqnarray*}
d'o\`u:
\begin{eqnarray*}
& &|\beta-\alpha|\|T(t_0\pm\varepsilon)x-T(t_0)x\|\leq\\
&\leq&\int\limits_{\alpha}^{\beta}\!\|T(t)\|\:\|T(t_0\pm\varepsilon-t)x-T(t_0-t)x\|\:dt\leq\\
&\leq&M\int\limits_{t_0-\beta}^{t_0-\alpha}\!\|T(\tau\pm\varepsilon)x-T(\tau)x\|\:d\tau\longrightarrow
0\quad\mbox{si}\quad\varepsilon\searrow 0\quad,
\end{eqnarray*}
compte tenu de \cite[th\'eor\`eme 3.6.3, pag.46]{hille}. Il s'ensuit
que l'application:
$$
[0,\infty)\ni t\longmapsto T(t)\in{\cal B(E)}
$$
est fortement continue sur $]0,\infty)$.\\
En particulier, pour $x\in{\cal E}$ arbitrairement fix\'e,
l'ensemble:
$$
{\cal X}=\{T(t)x|t\in[0,1]\}
$$
est s\'eparable. Donc il contient une partie d\'enombrable dense:
$$
{\cal X}_0=\{T(t_n)x|t_n\in]0,1[,\:n\in{\bf N}\}\quad.
$$
Par cons\'equent, il existe une suite $\left(x_n\right)_{n\in{\bf
N}}\subset{\cal X}_0$ tel que:
$$
\lim_{n\rightarrow\infty}\|x_n-x\|=\lim_{n\rightarrow\infty}\|T(t_n)x-x\|=0\quad.
$$
Comme:
\begin{eqnarray*}
&
&\|T(t)x-x\|\leq\\
&\leq&\|T(t)x-T(t+t_n)x\|+\|T(t+t_n)x-T(t_n)x\|+\|T(t_n)x-x\|\leq\\
&\leq&\|T(t)\|\:\|x-T(t_n)x\|+\|T(t+t_n)x-T(t_n)x\|+\|T(t_n)x-x\|\leq\\
&\leq&\left(\sup_{t\in[0,1]}\|T(t)\|+1\right)+\|T(t+t_n)x-T(t_n)x\|\quad,
\end{eqnarray*}
il vient:
$$
\lim_{t\searrow 0}T(t)x=x\quad,\quad(\forall)x\in{\cal E}
$$
et par cons\'equent:
$$
\lim_{t\searrow 0}T(t)=I
$$
dans la topologie forte.\fin

Dans la suite, nous consid\'erons la topologie forte pour \'etudier
les propri\'et\'es des $C_0$-semi-groupes.

\begin{teo}\label{num25}
Soit $\left\{T(t)\right\}_{t\geq 0}$ un $C_0$-semi-groupe
d'op\'erateurs lin\'eaires born\'es. Alors:\\
i) il existe $\tau>0$ et $M\geq 1$ tel que:
$$
\left\|T(t)\right\|\leq M\quad,\quad(\forall)t\in[0,\tau];
$$
ii) il existe $\omega\in {\bf R}$ et $M\geq 1$ tel que:
$$
\left\|T(t)\right\|\leq Me^{\omega t}\quad,\quad(\forall)t\geq 0.
$$
\end{teo}
\dem
i) Supposons  que pour tout $\tau>0$ et tout $M\geq 1$, il existe $t\in [0,\tau]$ tel que $\left\|T(t)\right\|>M$. Pour
$\tau=\frac{1}{n}$ et $M=n\in{\bf N}^{*}$, il existe
$t_n\in\left[0,\frac{1}{n}\right]$ tel que
$\left\|T(t_n)\right\|>n$.
Donc la suite $\left(\left\|T(t_n)\right\|\right)_{n\in{\bf N}^{*}}$ est
non born\'ee. Si la suite $\left(\left\|T(t_n)x\right\|\right)_{n\in{\bf
N}^{*}}$ \'etait born\'ee pour tout $x\in{\cal E}$, alors compte tenu
du th\'eor\`eme de Banach-Steinhaus (\cite[Theorem II.1.11, pag. 52]{dunford-schwartz}), il en r\'esulterait que
$\left(\left\|T(t_n)\right\|\right)_{n\in{\bf N}^{*}}$ serait born\'ee,
mais cela contredit l'affirmation pr\'ec\'edente. Donc il existe
$x_0\in{\cal E}$ tel que
$\left(\left\|T(t_n)x_0\right\|\right)_{n\in{\bf N}^{*}}$ soit non born\'ee.
D'autre part, compte tenu de la d\'efinition \ref{num10} (iii), il
r\'esulte que $\lim_{n\rightarrow\infty}\left\|T(t_n)x_0\right\|=x_0$ et
cela est contradictoire.\\
ii) Pour $h>0$ et $t>h$, nous noterons $m=\left[\frac{t}{h}\right]\in{{\bf
N}^{*}}$.
Compte tenu du th\'eor\`eme de division avec reste, il existe $r\in[0,h)$ tel que $t=mh+r$. Alors:
\begin{eqnarray*}
\left\|T(t)\right\|&=&\left\|T(mh)T(r)\right\| \leq
\left\|T(h)\right\|^m\left\|T(r)\right\| \leq\\
&\leq& M^mM\leq Me^{\frac{t}{h}\ln M}\quad.
\end{eqnarray*}
L'in\'egalit\'e de l'\'enonc\'e en r\'esulte en prenant
$\omega=\frac{1}{h}\ln M$.\fin

\begin{cor}
Si $\left\{T(t)\right\}_{t\geq 0}$ est un $C_0$-semi-groupe, alors
l'application:
$$
[0,\infty)\ni t\longmapsto T(t)x\in {\cal E}
$$
est continue sur $[0,\infty)$, quel que soit $x\in {\cal E}$.
\end{cor}
\dem
Soient $t_0,h\in [0,\infty)$ et $x\in {\cal E}$.\\
Si $t_0<h$, nous avons:
\begin{eqnarray*}
& &\left\|T(t_0+h)x-T(t_0)x\right\|\leq
\left\|T(t_0)\right\|\left\|T(h)x-x\right\|\leq\\
&\leq& Me^{\omega t_0}\left\|T(h)x-x\right\|\quad.
\end{eqnarray*}
Si $t_0>h$, nous obtenons:
\begin{eqnarray*}
& &\left\|T(t_0-h)x-T(t_0)x\right\|\leq
\left\|T(t_0-h)\right\|\left\|T(h)x-x\right\|\leq\\
&\leq& Me^{\omega (t_0-h)}\left\|T(h)x-x\right\|\quad.
\end{eqnarray*}
La continuit\'e forte en $t_0$ de l'application consid\'er\'ee dans l'\'enonc\'e est
\'evidente.\fin

\begin{definitie}
On dit que le $C_0$-semi-groupe $\left\{T(t)\right\}_{t\geq 0}$ est
uniform\'ement born\'e s'il existe $M\geq 1$ tel que:
$$
\left\|T(t)\right\|\leq M\quad,\quad(\forall)t\geq 0.
$$
\end{definitie}

\begin{teo}\label{num44}
Soit $\left\{T(t)\right\}_{t\geq 0}$ un $C_0$-semi-groupe pour lequel
il existe $\omega\in{\bf R}$ et $M\geq 1$ tel que:
$$
\left\|T(t)\right\|\leq Me^{\omega t}\quad,\quad(\forall)t\geq 0.
$$
Alors la famille $\left\{S(t)\right\}_{t\geq 0}\subset{\cal B(E)}$,
o\`u:
$$
S(t)=e^{-\omega t}T(t)\quad,\quad(\forall)t\geq 0,
$$
est un $C_0$-semi-groupe ayant la propri\'et\'e:
$$
\left\|S(t)\right\|\leq M\quad,\quad(\forall)t\geq 0.
$$
De plus, si $A$ est le g\'en\'erateur infinit\'esimal du
$C_0$-semi-groupe $\left\{T(t)\right\}_{t\geq 0}$, alors le
$C_0$-semi-groupe $\left\{S(t)\right\}_{t\geq 0}$ a pour g\'en\'erateur
infinit\'esimal l'op\'erateur $B=A-\omega I$.
\end{teo}
\dem
Dans les conditions du th\'eor\`eme, il est \'evident que
$\left\{S(t)\right\}_{t\geq 0}$ est un $C_0$-semi-groupe et:
$$
\left\|S(t)\right\|=\left\|e^{\omega t}T(t)\right\|\leq e^{-\omega
t}Me^{\omega t}=M\quad,\quad(\forall)t\geq 0.
$$
Donc $\left\{S(t)\right\}_{t\geq 0}$ est un $C_0$-semi-groupe
uniform\'ement born\'e. Soit $A$ le g\'en\'erateur infinit\'esimal du
$C_0$-semi-groupe $\left\{T(t)\right\}_{t\geq 0}$. Si $B$ est le g\'en\'erateur infinit\'esimal du $C_0$-semi-groupe
$\left\{S(t)\right\}_{t\geq 0}$, alors pour tout $x\in{\cal D}(A)$, nous
avons:
\begin{eqnarray*}
& &\lim_{h\searrow 0}\frac{S(h)x-x}{h}=\lim_{h\searrow
0}\frac{e^{-\omega h}T(h)x-x}{h}=\\
&=&\lim_{h\searrow 0}\frac{\left(e^{-\omega
h}-1\right)T(h)x}{h}+\lim_{h\searrow 0}\frac{T(h)x-x}{h}=\\
&=&-\omega x+Ax=(A-\omega I)x\quad,
\end{eqnarray*}
d'o\`u il r\'esulte que $x\in{\cal D}(B)$ et $Bx=(A-\omega I)x$. Soit
$x\in{\cal D}(A)$. Alors, nous obtenons:
\begin{eqnarray*}
& &\lim_{h\searrow 0}\frac{T(h)x-x}{h}=\lim_{h\searrow
0}\frac{e^{\omega h}S(h)x-x}{h}=\\
&=&\lim_{h\searrow 0}\frac{\left(e^{\omega
h}-1\right)S(h)}{h}+\lim_{h\searrow 0}\frac{S(h)x-x}{h}=\\
&=&(\omega I+B)x\quad,
\end{eqnarray*}
d'o\`u il vient que $x\in{\cal D}(A)$ et $Ax=(\omega I+B)x$. Par
cons\'equent ${\cal D}(A)={\cal D}(B)$ et $B=A-\omega I$.\fin

\begin{obs}
\em
Soit $\left\{T(t)\right\}_{t\geq 0}$ un $C_0$-semi-groupe pour lequel
il existe $\omega\in{\bf R}$ et $M\geq 1$ tel que:
$$
\left\|T(t)\right\|\leq Me^{\omega t}\quad,\quad(\forall)t\geq 0.
$$
Si $\omega<0$, alors nous obtenons:
$$
\left\|T(t)\right\|\leq Me^{\omega t}\leq M\quad,\quad(\forall)t\geq
0.
$$
Par cons\'equent on peut consid\'erer que $\omega\geq 0$.
\end{obs}
Nous noterons par ${\cal SG}(M,\omega)$ l'ensemble des
$C_0$-semi-groupes $\left\{T(t)\right\}_{t\geq 0}\subset {\cal B(E)}$
pour lesquels il existe $\omega\geq 0$ et $M\geq 1$ tel que:
$$
\left\|T(t)\right\|\leq Me^{\omega t}\quad,\quad(\forall)t\geq 0\:.
$$
Avec le th\'eor\`eme \ref{num44} nous voyons que le passage entre la
classe ${\cal SG}(M,\omega)$ avec $\omega>0$ et la classe ${\cal
SG}(M,0)$ est tr\`es simple.
\vspace{2cm}

  \section{Propri\'et\'es g\'en\'erales des $C_0$-semi-groupes}

\begin{prop}
Soient \semi et $A$
son g\'en\'erateur infinit\'esi-\\mal. Si $x\in {\cal D}(A)$, alors
$T(t)x\in {\cal D}(A)$ et on a l'\'egalit\'e:
$$
T(t)Ax=AT(t)x\quad,\quad(\forall)t\geq 0.
$$
\end{prop}
\dem
Soit $x\in {\cal D}(A)$. Alors pour tout $t\geq 0$, nous avons:
\begin{eqnarray*}
T(t)Ax&=&T(t)\lim_{h\searrow 0}\frac{T(h)x-x}{h}=\\
&=&\lim_{h\searrow 0}\frac{T(h)T(t)x-T(t)x}{h}\quad.
\end{eqnarray*}
Donc $T(t)x\in {\cal D}(A)$ et on a $T(t)Ax=AT(t)x$ , $(\forall)t\geq 0$.\fin

\begin{obs}
\em
On voit que:
$$
T(t){\cal D}(A)\subseteq{\cal D}(A)\quad,\quad(\forall)t\geq 0.
$$
\end{obs}

\begin{teo}\label{num22}
Soient \semi et $A$
son g\'en\'erateur infinit\'esimal. Alors l'application:
$$
[0,\infty)\ni t\longmapsto T(t)x\in {\cal E}
$$
est d\'erivable sur $[0,\infty)$, pour tout $x\in {\cal D}(A)$ et
nous avons:\\
i) ${\textstyle\frac{d}{dt}}T(t)x=T(t)Ax=AT(t)x\quad,\quad(\forall)t\geq 0;$\\
ii) $T(t)x-x=\int\limits_{0}^{t}\!T(s)Ax\:ds\quad,\quad(\forall)t\geq 0.$
\end{teo}
\dem
i) Soient $x\in {\cal D}(A)$ , $t\geq 0$ et $h>0$. Alors:
\begin{eqnarray*}
& &\left\|\frac{T(t+h)x-T(t)x}{h}-T(t)Ax\right\|\leq
\left\|T(t)\right\|\left\|\frac{T(h)x-x}{h}-Ax\right\|\leq\\
& &\leq Me^{\omega t}\left\|\frac{T(h)x-x}{h}-Ax\right\|\quad.
\end{eqnarray*}
Par cons\'equent:
$$
\lim_{h\searrow 0}\frac{T(t+h)x-T(t)x}{h}=T(t)Ax\quad,
$$
d'o\`u:
$$
\frac{d^+}{dt}T(t)x=T(t)Ax\quad,\quad(\forall)t\geq 0.
$$
Si $t-h>0$, alors nous avons:
\begin{eqnarray*}
& &\left\|\frac{T(t-h)x-T(t)x}{-h}-T(t)Ax\right\|\leq\\
& &\leq\left\|T(t-h)\right\|\left\|\frac{T(h)x-x}{h}-Ax+Ax-T(h)Ax\right\|\leq\\
& &\leq Me^{\omega
(t-h)}\left(\left\|\frac{T(h)x-x}{h}-Ax\right\|+\left\|T(h)Ax-Ax\right\|\right)\quad.
\end{eqnarray*}
Par suite:
$$
\lim_{h\searrow 0}\frac{T(t-h)x-T(t)x}{-h}=T(t)Ax
$$
et:
$$
\frac{d^-}{dt}T(t)x=T(t)Ax\quad,\quad(\forall)t\geq 0.
$$
Il s'ensuit que l'application consid\'er\'ee dans l'\'enonc\'e est d\'erivable sur
$[0,\infty)$, quel que soit $x\in {\cal D}(A)$. De plus, on a
l'\'egalit\'e:
$$
\frac{d}{dt}T(t)x=T(t)Ax=AT(t)x\quad,\quad(\forall)t\geq 0.
$$
ii) Si $x\in {\cal D}(A)$, alors nous avons:
$$
\frac{d}{ds}T(s)x=T(s)Ax\quad,\quad(\forall)s\in [0,t]\:,\:t\geq 0,
$$
d'o\`u:
$$
\int\limits_{0}^{t}\!T(s)Ax\:ds=\int\limits_{0}^{t}\frac{d}{ds}T(s)\:ds=T(t)x-x\quad,\quad(\forall)t\geq
0.\fin
$$

On peut obtenir une formule de repr\'esentation de type Taylor pour les $C_0$-semi-groupes
avec la g\'en\'eralisation du th\'eor\`eme \ref{num22} (ii).

\begin{teo}[Taylor]
Soient $\left\{T(t)\right\}_{t\geq 0}\in{\cal SG}(M,\omega)$ et $A$ son
g\'en\'erateur infinit\'esimal. Alors:
$$
T(t)x=\sum\limits_{i=0}^{n-1}\frac{t^i}{i!}A^ix+\frac{1}{(n-1)!}\int\limits_{0}^{t}\!(t-u)^{n-1}T(u)A^nx\:du
$$
quels que soient $x\in{\cal D}(A^n)$, $t\geq 0$ et $n\in{\bf N}^*$.
\end{teo}
\dem
Compte tenu du th\'eor\`eme \ref{num22} (ii), pour $x\in{\cal D}(A)$
et $t\geq 0$ on a:
$$
T(t)x=x+\int\limits_{0}^{t}\!T(u)Ax\:du\quad.
$$
Supposons que pour $t\geq 0$ et $x\in{\cal D}(A^k)$ nous ayons:
$$
T(t)x=\sum\limits_{i=0}^{k-1}\frac{t^i}{i!}A^ix+\frac{1}{(k-1)!}\int\limits_{0}^{t}\!(t-u)^{k-1}T(u)A^kx\:du\quad.
$$
Si $x\in{\cal D}(A^{k+1})$, alors $x\in{\cal D}(A^k)$ et
$A^kx\in{\cal D}(A)$. Il en r\'esulte que:
$$
T(t)x=\sum\limits_{i=0}^{k-1}\frac{t^i}{i!}A^ix+\frac{1}{(k-1)!}\int\limits_{0}^{t}\!(t-s)^{k-1}T(s)A^nx\:ds\quad.
$$
Mais:
$$
T(s)x=x+\int\limits_{0}^{s}\!T(u)Ax\:du\quad.
$$
Il vient:
$$
(t-s)^{k-1}T(s)A^kx=(t-s)^{k-1}A^kx+(t-s)^{k-1}\int\limits_{0}^{s}\!T(u)A^{k+1}x\:du
$$
et par cons\'equent:
\begin{eqnarray*}
& &\int\limits_{0}^{t}\!(t-s)^{k-1}T(s)A^kx\:ds=\\
&=&\int\limits_{0}^{t}\!(t-s)^{k-1}A^kx\:ds+\int\limits_{0}^{t}\!(t-s)^{k-1}\:\int\limits_{0}^{s}\!T(u)A^{k+1}x\:du\;ds=\\
&=&\frac{t^k}{k}A^kx+\int\limits_{0}^{t}\!\int\limits_{u}^{t}\!(t-s)^{k-1}T(u)A^{k+1}x\:ds\:du=\\
&=&\frac{t^k}{k}A^kx+\int\limits_{0}^{t}\frac{(t-u)^k}{k}T(u)A^{k+1}x\:du\quad.
\end{eqnarray*}
Nous en d\'eduisons que:
\begin{eqnarray*}
T(t)x&=&\sum\limits_{i=0}^{k-1}\frac{t^i}{i!}A^ix+\frac{1}{(k-1)!}\left[\frac{t^k}{k}A^kx+\frac{1}{k}\int\limits_{0}^{t}\!(t-u)^kT(u)A^{k+1}x\:du\right]=\\
&=&\sum\limits_{i=0}^{k}\frac{t^i}{k!}A^ix+\frac{1}{k!}\int\limits_{0}^{t}\!(t-u)^kT(u)A^{k+1}x\:du\quad,
\end{eqnarray*}
d'o\`u il
r\'esulte l'\'egalit\'e consid\'er\'ee dans l'\'enonc\'e.\fin

\begin{lema}\label{num11}
Soit $\left\{T(t)\right\}_{t\geq 0}$ un $C_0$-semi-groupe. Alors:
$$
\lim_{h\searrow 0}\frac{1}{h}\int\limits_{t}^{t+h}\!T(s)x\:ds=T(t)x
$$
quels que soient $x\in{\cal E}$ et $t\geq 0$.
\end{lema}
\dem
L'\'egalit\'e de l'\'enonc\'e r\'esulte de
l'\'evaluation:
\begin{eqnarray*}
& &\left\|\frac{1}{h}\int\limits_{t}^{t+h}\!T(s)x\:ds-T(t)x\right\|=\left\|\frac{1}{h}\int\limits_{t}^{t+h}\left(T(s)-T(t)\right)x\:ds\right\|\leq\\
&\leq& \sup_{s\in [t,t+h]}\left\|T(s)x-T(t)x\right\|
\end{eqnarray*}
et de la continuit\'e de l'application $[0,\infty)\ni t\longmapsto
T(t)x\in {\cal E}$.\fin

\begin{prop}
Soient $\left\{T(t)\right\}_{t\geq 0}\in{\cal SG}(M,\omega)$ et $A$ son
g\'en\'erateur infinit\'esimal. Si $x\in{\cal E}$, alors
$\int\limits_{0}^{t}\!T(s)x\:ds\in{\cal D}(A)$ et on a l'\'egalit\'e:
$$
A\int\limits_{0}^{t}\!T(s)x\:ds=T(t)x-x\quad,\quad(\forall)t\geq 0.
$$
\end{prop}
\dem
Soient $x\in{\cal E}$ et $h>0$. Alors:
\begin{eqnarray*}
&
&\frac{T(h)-I}{h}\int\limits_{0}^{t}\!T(s)x\:ds=\frac{1}{h}\int\limits_{0}^{t}\!T(s+h)x\:ds-\frac{1}{h}\int\limits_{0}^{t}\!T(s)x\:ds=\\
&=&\frac{1}{h}\int\limits_{h}^{t+h}\!T(u)x\:du-\frac{1}{h}\int\limits_{0}^{t}\!T(s)x\:ds=\\
&=&\frac{1}{h}\int\limits_{0}^{t+h}\!T(u)x\:du-\frac{1}{h}\int\limits_{0}^{h}\!T(u)x\:du-\frac{1}{h}\int\limits_{0}^{t}\!T(u)x\:du=\\
&=&\frac{1}{h}\int\limits_{t}^{t+h}\!T(u)x\:du-\frac{1}{h}\int\limits_{0}^{h}\!T(u)x\:du\quad.
\end{eqnarray*}
Par pasage \`a limite pour $h\searrow 0$ et compte tenu du lemme
\ref{num11}, nous obtenons:
$$
A\int\limits_{0}^{t}\!T(s)x\:ds=T(t)x-x\quad,\quad(\forall)t\geq 0
$$
et:
$$
\int\limits_{0}^{t}\!T(s)x\:ds\in{\cal D}(A).\fin
$$

\begin{teo}\label{num12}
Soient $\left\{T(t)\right\}_{t\geq 0}\in{\cal SG}(M,\omega)$ et $A$ son
g\'en\'erateur infinit\'esimal. Alors:\\
i) $\overline{{\cal D}(A)}={\cal E};$\\
ii) $A$ est un op\'erateur ferm\'e.
\end{teo}
\dem
i) Soient $x\in{\cal E}$ et $t_n>0$ , $n\in{\bf N}$, tel que
$\lim_{n\rightarrow\infty}t_n=0$. Alors:
$$
x_n=\frac{1}{t_n}\int\limits_{0}^{t_n}\!T(s)x\:ds\in{\cal
D}(A)\quad,\quad(\forall)n\in{\bf N},
$$
d'o\`u:
$$\lim_{n\rightarrow\infty}x_n=\lim_{n\rightarrow\infty}\frac{1}{t_n}\int\limits_{0}^{t_n}\!T(s)x\:ds=T(0)x=x\quad.
$$
Par cons\'equent $\overline{{\cal D}(A)}={\cal E}$.\\
ii) Soit $\left(x_n\right)_{n\in{\bf N}}\subset{\cal D}(A)$ tel que
$\lim_{n\rightarrow\infty}x_n=x$ et
$\lim_{n\rightarrow\infty}Ax_n=y$. Alors:
$$
\left\|T(s)Ax_n-T(s)y\right\|\leq\left\|T(s)\right\|\left\|Ax_n-y\right\|\leq
Me^{\omega t}\left\|Ax_n-y\right\|
$$
quel que soit $s\in [0,t]$. Par suite
$T(s)Ax_n\longrightarrow T(s)y$, pour $n\rightarrow\infty$,
uniform\'ement par rapport \`a $s\in [0,t]$.\\
D'autre part, puisque $x_n\in{\cal D}(A)$, nous avons:
$$
T(t)x_n-x_n=\int\limits_{0}^{t}\!T(s)Ax_n\:ds\quad,
$$
d'o\`u:
$$
\lim_{n\rightarrow\infty}\left[T(t)x_n-x_n\right]=\lim_{n\rightarrow\infty}\int\limits_{0}^{t}\!T(s)Ax_n\:ds\quad,
$$
ou bien:
$$
T(t)x-x=\int\limits_{0}^{t}\!T(s)y\:ds\quad.
$$
Finalement, on voit que:
$$
\lim_{t\searrow 0}\frac{T(t)x-x}{t}=\lim_{t\searrow
0}\frac{1}{t}\int\limits_{0}^{t}\!T(s)y\:ds=y\quad.
$$
Par suite $x\in{\cal D}(A)$ et $Ax=y$, d'o\`u il r\'esulte que $A$ est un
op\'erateur ferm\'e.\fin

Nous montrons maintenant un r\'esultat qui concerne l'unicit\'e
de l'engendrement pour les $C_0$-semi-groupes.

\begin{teo}[l'unicit\'e de l'engendrement]
Soient deux $C_0$-semi-groupes $\left\{T(t)\right\}_{t\geq 0}$ et $\left\{S(t)\right\}_{t\geq
0}$ ayant pour g\'en\'erateur
infinit\'esimal le m\^eme op\'erateur $A$. Alors:
$$
T(t)=S(t)\quad,\quad(\forall)t\geq 0.
$$
\end{teo}
\dem
Soient $t>0$ et $x\in{\cal D}(A)$. D\'efinissons l'application:
$$
[0,t]\ni s\longmapsto U(s)x=T(t-s)S(s)x\in{\cal D}(A).
$$
Alors:
\begin{eqnarray*}
\frac{d}{ds}U(s)x&=&\frac{d}{ds}T(t-s)S(s)x+T(t-s)\frac{d}{ds}S(s)x=\\
&=&-AT(t-s)S(s)x+T(t-s)AS(s)x=0
\end{eqnarray*}
quel que soit $x\in{\cal D}(A)$. Par suite $U(0)x=U(t)x$, pour tout $x\in{\cal D}(A)$, d'o\`u:
$$
T(t)x=S(t)x\quad,\quad(\forall)x\in{\cal D}(A)\mbox{ et }t\geq 0.
$$
Puisque $\overline{{\cal D}(A)}={\cal E}$ et $T(t),S(t)\in{\cal
B(E)}$, pour tout $t\geq 0$, il r\'esulte que:
$$
T(t)x=S(t)x\quad,\quad(\forall)t\geq 0\mbox{ et }x\in{\cal E}\quad,
$$
ou bien:
$$
T(t)=S(t)\quad,\quad(\forall)t\geq 0.\fin
$$

\begin{teo}\label{num57}
Soient $\left\{T(t)\right\}_{t\geq 0}$ un $C_0$-semi-groupe, $A$ son
g\'en\'erateur infinit\'e-\\simal et $F\in{\cal B(E)}$. Alors $T(t)F=FT(t)$
pour tout $t\geq 0$ si et seulement si:
$$
F{\cal D}(A)\subseteq{\cal D}(A)
$$
et:
$$
FAx=AFx\quad,\quad(\forall)x\in{\cal D}(A).
$$
\end{teo}
\dem
$\Longrightarrow$ Soit $F\in{\cal B(E)}$ tel que:
$$
T(t)F=FT(t)\quad,\quad(\forall)t\geq 0
$$
et $x\in{\cal D}(A)$. Alors, nous avons:
\begin{eqnarray*}
& &\lim_{t\searrow 0}\frac{T(t)Fx-Fx}{t}=\lim_{t\searrow
0}\frac{FT(t)x-Fx}{t}=\\
&=&\lim_{t\searrow 0}F\frac{T(t)x-x}{t}\quad.
\end{eqnarray*}
Par cons\'equent $Fx\in{\cal D}(A)$ et on a $AFx=FAx$, pour tout $x\in{\cal
D}(A).$\\
$\Longleftarrow$ Soit $F\in{\cal B(E)}$ tel que:
$$
F{\cal D}(A)\subseteq{\cal D}(A)
$$
et
$$
AFx=FAx\quad,\quad(\forall)x\in{\cal D}(A).
$$
Pour tout $t\geq 0$ et tout $x\in{\cal D}(A)$, d\'efinissons l'application:
$$
[0,t]\ni s\longmapsto U(s)x=T(t-s)FT(s)x\in{\cal D}(A)\quad.
$$
Alors nous avons:
\begin{eqnarray*}
&
&\frac{d}{ds}U(s)x=\frac{d}{ds}T(t-s)FT(s)x+T(t-s)\frac{d}{ds}FT(s)x=\\
&=&-AT(t-s)FT(s)x+T(t-s)FAT(s)x=0\quad,
\end{eqnarray*}
compte tenu de la commutativit\'e. Par cons\'equent:
$$
U(0)x=U(t)x\quad,\quad(\forall)x\in{\cal D}(A),
$$
d'o\`u on obtient:
$$
T(t)Fx=FT(t)x\quad,
$$
pour tout $t\geq 0$ et tout $x\in{\cal D}(A)$. Comme $\overline{{\cal D}(A)}={\cal E}$ et
$T(t)F,FT(t)\in{\cal B(E)}$ pour tout $t\geq 0$, nous obtenons:
$$
T(t)Fx=FT(t)x\quad,
$$
pour tout $t\geq 0$ et tout $x\in{\cal E}$.\fin

Nous finissons cette section avec une g\'en\'eralisation du th\'eor\`eme
\ref{num12}.

\begin{teo}\label{num37}
Soient \semi et $A$ son
g\'en\'erateur infinit\'esimal. Alors:\\
i) $\overline{{\cal D}(A^p)}={\cal E}$, quel que soit $p\in{\bf N}^*$;\\
ii) $A^p$ est un op\'erateur ferm\'e, quel que soit $p\in{\bf N}^*$;\\
iii) l'application:
$$
\|\:.\:\|_{{\cal D}(A^p)}:{\cal D}(A^p)\longrightarrow {\bf R}_+\quad,
$$
$$
\|x\|_{{\cal D}(A^p)}=\sum\limits_{i=0}^{p}\left\|A^ix\right\|
$$
est une norme avec laquelle ${\cal D}(A^p)$ devient un espace de
Banach, pour tout $p\in{\bf N}^*$.
\end{teo}
\dem
i) Pour $p=1$, compte tenu du th\'eor\`eme \ref{num12}(i), il r\'esulte
que $\overline{{\cal D}(A)}={\cal E}$.\\
Soit:
$$
{\cal C}_0^\infty=\left\{\varphi:]0,\infty)\rightarrow {\cal E}\left|
\varphi\;\mbox{\em ind\'efiniment d\'erivable avec un support compact}\right.\right\}.
$$
Notons:
$$
{\cal F}=\left\{\left.\int\limits_{0}^{\infty}\!\varphi
(t)T(t)x\:dt\right|x\in{\cal E}\:,\:\varphi \in {\cal
C}_0^\infty\right\}.
$$
Nous montrons que ${\cal F}\subset {\cal D}(A^p)$ , $(\forall)p\in{\bf
N}$.\\
Pour $y\in{\cal F}$ et $h>0$, nous obtenons:
\begin{eqnarray*}
\frac{T(h)-I}{h}\:y&=&\frac{1}{h}\left(\int\limits_{0}^{\infty}\!\varphi
(t)T(t+h)x\:dt-\int\limits_{0}^{\infty}\!\varphi(t)T(t)x\:dt\right)=\\
&=&\int\limits_{0}^{\infty}\frac{\varphi (u-h)-\varphi
(u)}{h}T(u)x\:du\quad.
\end{eqnarray*}
Puisque:
$$
\frac{\varphi (u-h)-\varphi (u)}{h}T(u)x\longrightarrow
-{\varphi(u)}^{(1)}T(u)x\quad\mbox{si}\quad h\searrow 0,
$$
uniform\'ement par rapport \`a $u\in\mbox{\em supp }\varphi$,
en passant \`a limite pour $h\searrow 0$, nous obtenons:
$$
Ay=-\int\limits_{0}^{\infty}\!{\varphi (u)}^{(1)}T(u)x\:du\quad.
$$
Donc $y\in {\cal D}(A)$. Il en r\'esulte que ${\cal F}\subset{\cal
D}(A)$ et par r\'ecurrence on peut montrer que ${\cal
F}\subset{\cal D}(A^p)$ et:
$$
A^py=(-1)^p\int\limits_{0}^{\infty}\!{\varphi (t)}^{(p)}T(t)x\:dt
$$
quel que soit $p\in{\bf N}^*$.\\
Nous montrons maintenant que ${\cal F}$ est dense dans ${\cal E}$.\\
Supposons que ${\cal F}$ n'est pas dense dans ${\cal
E}$. Alors il existe $x_0\in{\cal E}$ tel que $d(x_0,{\cal
F})>0$. En appliquant un corollaire du th\'eor\`eme de Hahn-Banach
(\cite[Corollary II.3.13, pag. 64]{dunford-schwartz}), on voit qu'il
existe $x_0^*\in{\cal E}^*$ tel que $\langle x_0,x_0^*\rangle=1$ et
$\langle y,x_0^*\rangle=0$, pour tout $y\in{\cal F}$. Alors:
$$
\int\limits_{0}^{\infty}\!\varphi
(t)\langle T(t)x,x_0^*\rangle\:dt
=\left\langle\int\limits_{0}^{\infty}\!\varphi
(t)T(t)x\:dt,x_0^*\right\rangle=0\quad,\quad(\forall)\varphi\in{\cal
C}_0^\infty\mbox{ et }x\in{\cal E}.
$$
Par cons\'equent, pour tout $x\in{\cal E}$, nous avons:
$$
\langle T(t)x,x_0^*\rangle=0\quad,\quad(\forall)t\in[0,\infty),
$$
parce que dans le cas contraire, on peut trouver $\varphi\in{\cal
C}_0^\infty$ tel que:
$$
\int\limits_{0}^{\infty}\!\varphi
(t)\langle T(t)x,x_0^*\rangle\:dt\neq 0
$$
ce qui est contradictoire. Il s'ensuit que pour tout $x\in{\cal E}$,
on a:
$$
\langle T(t)x,x_0^*\rangle=0\quad,\quad(\forall)t\in[0,\infty),
$$
d'o\`u:
$$
\langle x,x_0^*\rangle=\langle
T(0)x,x_0^*\rangle=0\quad,\quad(\forall)x\in{\cal E},
$$
ce qui est absurde. Finalement, on voit que ${\cal F}$ est dense dans ${\cal E}$ et donc $\overline{{\cal
D}(A^n)}={\cal E}$.\\
ii) Compte tenu du th\'eor\`eme \ref{num12}(ii), on voit que:
$$
A:{\cal
D}(A)\subset{\cal E}\longrightarrow{\cal E}
$$
est un op\'erateur
ferm\'e. Supposons que:
$$
A^k:{\cal D}(A^k)\subset{\cal E}\longrightarrow{\cal E}
$$
est un op\'erateur ferm\'e et montrons que:
$$
A^{k+1}:{\cal D}(A^{k+1})\subset{\cal E}\longrightarrow{\cal E}
$$
est un op\'erateur ferm\'e.\\
Soit $(x_n)_{n\in{\bf N}}\subset{\cal D}(A^{k+1})$ tel que:
$$
\lim_{n\rightarrow\infty}x_n=x
$$
et:
$$
\lim_{n\rightarrow\infty}A^{k+1}x_n=y\quad.
$$
Mais $x_n\in{\cal D}(A^{k+1})$ est \'equivalent avec $x_n\in{\cal
D}(A^{k})$ et $A^kx_n\in{\cal D}(A)$. Alors $x_n\in{\cal D}(A^k)$,
$\lim_{n\rightarrow\infty}x_n=x$, comme $A^k$ est un op\'erateur ferm\'e, ceci implique
$x\in{\cal D}(A^k)$ et $\lim_{n\rightarrow\infty}A^kx_n=A^kx$. Comme
$A^kx_n\in{\cal D}(A)$, $\lim_{n\rightarrow\infty}A^kx_n=A^kx$ et $A$
est un op\'erateur ferm\'e, il s'ensuit que $A^kx\in{\cal D}(A)$ et
$\lim_{n\rightarrow\infty}A\left(A^kx_n\right)=A\left(A^kx\right)$.
Nous avons obtenu donc que $x\in{\cal D}(A^{k+1})$, $A^kx\in{\cal
D}(A)$ et $\lim_{n\rightarrow\infty}A^{k+1}x_n=A^{k+1}x$, d'o\`u il
r\'esulte que $x\in{\cal D}(A^{k+1})$ et $A^{k+1}x=y$. Par
cons\'equent $A^{k+1}$ est un op\'erateur ferm\'e, d'o\`u on obtient (ii).\\
iii) Pour $p=1$ on peut v\'erifier facilement les propri\'et\'es de
norme de l'application:
$$
\|\:.\:\|_{{\cal D}(A)}:{\cal D}(A)\longrightarrow {\bf R}_+\quad,
$$
$$
\|x\|_{{\cal D}(A)}=\|x\|+\|Ax\|\quad.
$$
Donc ${\cal D}(A)$ est un espace norm\'e.\\
Soit $(x_n)_{n\in{\bf N}^*}\subset{\cal D}(A)$ tel que
$\|x_m-x_n\|_{{\cal D}(A)}\longrightarrow 0$ pour
$m,n\rightarrow\infty$.
Alors:
$$
\|x_m-x_n\|+\|Ax_m-Ax_n\|\longrightarrow 0\mbox{ pour
}m,n\rightarrow\infty.
$$
Donc:
$$
\|x_m-x_n\|\longrightarrow 0\mbox{ et }\|Ax_m-Ax_n\|\longrightarrow
0\mbox{ pour }m,n\rightarrow\infty.
$$
Puis que ${\cal E}$ est un espace de Banach, il r\'esulte que les
suites $(x_n)_{n\in{\bf N}}$ et $(Ax_n)_{n\in{\bf N}}$ sont
convergentes. Donc $x_n\longrightarrow x$ et $Ax_n\longrightarrow y$
pour $n\rightarrow\infty$. Comme $A$ est un op\'erateur ferm\'e, il
r\'esulte que $x\in{\cal D}(A)$ et $y=Ax$. Par cons\'equent:
$$
\|x_n-x\|_{{\cal D}(A)}=\|x_n-x\|+\|Ax_n-Ax\|\longrightarrow 0\mbox{
pour }n\rightarrow\infty.
$$
Donc la suite $(x_n)_{n\in{\bf N}}$ est convergente par rapport \`a la
norme $\|\:.\:\|_{{\cal D}(A)}$. Il s'ensuit que ${\cal D}(A)$ est un
espace de Banach avec la norme $\|\:.\:\|_{{\cal D}(A)}$.\\
Supposons que l'application:
$$
\|\:.\:\|_{{\cal D}(A^k)}:{\cal D}(A^k)\longrightarrow {\bf R}_+\quad,
$$
$$
\|x\|_{{\cal D}(A^k)}=\sum\limits_{i=0}^{k}\left\|A^ix\right\|
$$
est une norme avec laquelle ${\cal D}(A^k)$ est un espace de
Banach. Montrons que:
$$
\|\:.\:\|_{{\cal D}(A^{k+1})}:{\cal D}(A^{k+1})\longrightarrow {\bf R}_+\quad,
$$
$$
\|x\|_{{\cal D}(A^{k+1})}=\sum\limits_{i=0}^{k+1}\left\|A^ix\right\|
$$
est une norme avec laquelle ${\cal D}(A^{k+1})$ devient un espace de
Banach. On peut v\'erifier facilement les propri\'et\'es de norme de
l'application $\|\:.\:\|_{{\cal D}(A^{k+1})}$. Donc ${\cal
D}(A^{k+1})$ est un espace norm\'e. Soit $(x_n)_{n\in{\bf N}}\subset{\cal D}(A^{k+1})$
tel que:
$$
\|x_m-x_n\|_{{\cal D}(A^{k+1})}\longrightarrow 0\quad\mbox{si}\quad m,n\rightarrow\infty.
$$
Alors nous avons:
$$
\sum\limits_{i=0}^{k+1}\left\|A^ix_m-A^ix_n\right\|\longrightarrow
0\quad\mbox{si}\quad m,n\rightarrow\infty,
$$
d'o\`u il s'ensuit que:
$$
\left\|A^ix_m-A^ix_n\right\|\longrightarrow 0\quad\mbox{si}\quad m,n\rightarrow\infty,
$$
pour tout $i\in\{0,1,\ldots,k+1\}$. Mais ${\cal E}$ est un espace de
Banach. Donc pour tout $i\in\{0,1,\ldots,k+1\}$, les suites
$\left(A^ix_n\right)_{n\in{\bf N}}$ sont convergentes et comme les
op\'erateurs $A^i$ sont ferm\'es pour tout $i\in\{1,2,\ldots,k+1\}$,
on voit que:
$$
\left\|A^ix_n-A^ix\right\|\longrightarrow 0\quad\mbox{si}\quad n\rightarrow\infty,
$$
pour tout $i\in\{0,1,\ldots,k+1\}$. Par cons\'equent:
$$
\sum\limits_{i=0}^{k+1}\left\|A^ix_n-A^ix\right\|\longrightarrow
0\quad\mbox{si}\quad n\rightarrow\infty,
$$
d'o\`u:
$$
\|x_m-x\|_{{\cal D}(A^{k+1})}\longrightarrow 0\quad\mbox{si}\quad n\rightarrow\infty.
$$
Finalement, on voit que ${\cal D}(A^{k+1})$ est un espace de Banach et l'affirmation de l'\'enonc\'e en r\'esulte.\fin
\vspace{2cm}

  \section{Le th\'eor\`eme de Hille - Yosida}

\hspace{1cm}Dans ce paragraphe nous pr\'esentons un r\'esultat tr\`es important
concernant les semi-groupes de classe $C_0$. Il s'agit du c\'el\`ebre
th\'eor\`eme de Hille-Yosida qui donne une caract\'erisation pour les
op\'erateurs qui sont g\'en\'erateurs de $C_0$-semi-groupes. Nous
avons besoin de quelques r\'esultats interm\'ediaires.
Dans la suite, pour $\omega\geq 0$ nous d\'esignerons par $\Lambda_\omega$ l'ensemble $\{\lambda\in{\bf
C}\left|\mbox{\em Re}\lambda>\omega\right.\}$.

\begin{teo}\label{num13}
Soient $\left\{T(t)\right\}_{t\geq 0}\in{\cal SG}(M,\omega)$ et $A$ son
g\'en\'erateur infinit\'esimal. Si $\lambda\in\Lambda_\omega$, alors l'application:
$$
R_\lambda:{\cal E}\longrightarrow{\cal E},
$$
$$
R_\lambda x=\int\limits_{0}^{\infty}\!\!e^{-\lambda t}T(t)x\:dt
$$
d\'efinit un op\'erateur lin\'eaire born\'e sur ${\cal E}$,
$\lambda\in\rho(A)$ et $R_\lambda x=R(\lambda;A)x$ ,
pour tout $x\in{\cal E}$.
\end{teo}
\dem
Soit $\lambda\in\Lambda_\omega$. Puisque
$\left\{T(t)\right\}_{t\geq 0}\in{\cal SG}(M,\omega)$, nous avons:
$$
\left\|T(t)\right\|\leq Me^{\omega t}\quad,\quad(\forall)t\geq 0
$$
et on voit que:
$$
\left\|e^{-\lambda t}T(t)x\right\|\leq e^{-{\scriptstyle{Re}}\lambda
t}\left\|T(t)\right\|\|x\|\leq Me^{-({\scriptstyle{Re}}\lambda
-\omega)t}\|x\|\quad,\quad(\forall)x\in{\cal E}.
$$
D\'efinissons l'application:
$$
R_\lambda:{\cal E}\longrightarrow{\cal E}\quad,
$$
par:
$$
R_\lambda x=\int\limits_{0}^{\infty}\!\!e^{-\lambda t}T(t)x\:dt\quad.
$$
Il est clair que $R_\lambda$ est un op\'erateur lin\'eaire. De plus, on a:
$$
\left\|R_\lambda
x\right\|\leq\int\limits_{0}^{\infty}\left\|e^{-\lambda
t}T(t)x\right\|dt\leq\frac{M}{\mbox{\em Re}\lambda
-\omega}\|x\|\quad,\quad(\forall)x\in{\cal E},
$$
d'o\`u il r\'esulte que $R_\lambda$ est un op\'erateur lin\'eaire born\'e.\\
Si $x\in{\cal E}$, alors nous avons:
\begin{eqnarray*}
& &\frac{T(h)R_\lambda x-R_\lambda
x}{h}=\frac{1}{h}\int\limits_{0}^{\infty}\!\!e^{-\lambda
t}T(t+h)x\:dt-\frac{1}{h}\int\limits_{0}^{\infty}\!\!e^{-\lambda
t}T(t)x\:dt=\\
&=&\frac{1}{h}\int\limits_{h}^{\infty}\!\!e^{-\lambda
(s-h)}T(s)x\:ds-\frac{1}{h}\int\limits_{0}^{\infty}\!\!e^{-\lambda
t}T(t)x\:dt=\\
&=&\frac{e^{\lambda h}}{h}\int\limits_{h}^{\infty}\!\!e^{-\lambda
s}T(s)x\:ds-\frac{1}{h}\int\limits_{0}^{\infty}\!\!e^{-\lambda
t}T(t)x\:dt=\\
&=&\frac{e^{\lambda
h}}{h}\left(\int\limits_{0}^{\infty}\!\!e^{-\lambda
s}T(s)x\:ds-\int\limits_{0}^{h}\!\!e^{-\lambda
s}T(s)x\:ds\right)-\frac{1}{h}\int\limits_{0}^{\infty}\!\!e^{-\lambda
t}T(t)x\:dt=\\
&=&\frac{e^{\lambda h}-1}{h}\int\limits_{0}^{\infty}\!\!e^{-\lambda
s}T(s)x\:ds-\frac{e^{\lambda
h}}{h}\int\limits_{0}^{\infty}\!\!e^{-\lambda s}T(s)x\:ds\quad.
\end{eqnarray*}
Par passage \`a limite, on obtient:
$$
\lim_{h\searrow 0}\frac{T(h)R_\lambda x-R_\lambda x}{h}=\lambda
R_\lambda x-x\quad.
$$
Il en r\'esulte que $R_\lambda x\in{\cal D}(A)$ et
$$
AR_\lambda x=\lambda R_\lambda x-x\quad,\quad(\forall)x\in{\cal E},
$$
ou bien
$$
(\lambda I-A)R_\lambda x=x\quad,\quad(\forall)x\in{\cal E}.
$$
Si $x\in{\cal D}(A)$, alors nous obtenons:
\begin{eqnarray*}
R_\lambda Ax&=&\int\limits_{0}^{\infty}\!\!e^{-\lambda
t}T(t)Ax\:dt=\int\limits_{0}^{\infty}\!\!e^{-\lambda
t}\frac{d}{dt}T(t)x\:dt=\\
&=&\left.\left[e^{-\lambda t}T(t)x\right]\right|_0^\infty
+\lambda\int\limits_{0}^{\infty}\!\!e^{-\lambda t}T(t)x\:dt=x+\lambda
R_\lambda x\quad,
\end{eqnarray*}
d'o\`u:
$$
R_\lambda (\lambda I-A)x=x\quad,\quad(\forall)x\in{\cal D}(A).
$$
Finalement, on voit que $\lambda\in\rho(A)$ et $R_\lambda
x=R(\lambda;A)x$ , pour tout $x\in{\cal E}$.\fin

\begin{obs}
\em
On voit que pour tout $\lambda\in\Lambda_\omega$ on a:
$$
{\cal I}\mbox{\em m }R(\lambda;A)={\cal I}\mbox{\em m
}R_\lambda\subseteq{\cal D}(A)
$$
et:
$$
R(\lambda;A){\cal D}(A)=R_\lambda{\cal D}(A)\subseteq{\cal D}(A)\quad.
$$
\end{obs}

\begin{definitie}
L'op\'erateur:
$$
R_\lambda :{\cal E}\longrightarrow{\cal E}
$$
$$
R_\lambda x=\int\limits_{0}^{\infty}\!\!e^{-\lambda
t}T(t)x\:dt\quad,\quad\lambda\in\Lambda_\omega,
$$
s'appelle la transform\'ee de Laplace du semi-groupe
$\left\{T(t)\right\}_{t\geq 0}\in{\cal SG}(M,\omega)$.
\end{definitie}

\begin{obs}
\em
Soient \semi et $A$ son g\'en\'
erateur infinit\'esimal. Alors nous avons:
$$
\left\{\lambda\in{\bf
C}\left|\mbox{\em Re}\lambda>\omega\right.\right\}\subset\rho(A).
$$
et:
$$
\sigma(A)\subset\left\{\lambda\in{\bf
C}\left|\mbox{\em Re}\lambda\leq\omega\right.\right\}.
$$
\end{obs}

\begin{teo}\label{num15}
Soient \semi et $A$ son g\'en\'erateur infinit\'esimal. Pour tout $\lambda\in\Lambda_\omega$ on a:
$$
\left\|{R(\lambda;A)}^n\right\|\leq\frac{M}{(\mbox{Re}\lambda-\omega)^n}\quad,\quad(\forall)n\in{\bf
N}^*.
$$
\end{teo}
\dem
Soit \semi. Alors:
$$
\left\|T(t)\right\|\leq Me^{\omega t}\quad,\quad(\forall)t\geq 0.
$$
Compte tenu du th\'eor\`eme \ref{num13}, si $\lambda\in\Lambda_\omega$, nous avons $\lambda\in\rho(A)$ et:
$$
R(\lambda;A)x=R_\lambda x=\int\limits_{0}^{\infty}\!\!e^{-\lambda
t}T(t)x\:dt\quad,\quad(\forall)x\in{\cal E}.
$$
De plus:
$$
\left\|R(\lambda;A)\right\|\leq\frac{M}{\mbox{\em Re}\lambda-\omega}\quad.
$$
Il est clair que:
$$
\frac{d}{d\lambda}R(\lambda;A)x=-\int\limits_{0}^{\infty}\!\!te^{-\lambda
t}T(t)x\:dt\quad,\quad(\forall)x\in{\cal E}
$$
et par r\'ecurrence on peut montrer que:
$$\frac{d^n}{d\lambda^n}R(\lambda;A)x=(-1)^n\int\limits_{0}^{\infty}\!\!t^ne^{-\lambda
t}T(t)x\:dt\quad,\quad(\forall)x\in{\cal E}\mbox{ et }n\in{\bf N}^*.
$$
D'autre part, avec la proposition \ref{num14} (iii) nous obtenons:
$$
\frac{d^n}{d\lambda^n}R(\lambda;A)x=(-1)^nn!{R(\lambda;A)}^{n+1}x\quad,\quad(\forall)x\in{\cal
E}\mbox{ et }n\in{\bf N}^*.
$$
Par suite, on a:
$$
(-1)^nn!{R(\lambda;A)}^{n+1}x=(-1)^n\int\limits_{0}^{\infty}\!\!t^ne^{-\lambda
t}T(t)x\:dt\quad,\quad(\forall)x\in{\cal
E}\mbox{ et }n\in{\bf N}^*,
$$
d'o\`u il r\'esulte que:
$$
{R(\lambda;A)}^nx=\frac{1}{(n-1)!}\int\limits_{0}^{\infty}\!\!t^{n-1}e^{-\lambda
t}T(t)x\:dt\quad,\quad(\forall)x\in{\cal E}\mbox{ et }n\in{\bf N}^*.
$$
De plus:
\begin{eqnarray*}
& &\left\|{R(\lambda;A)}^nx\right\|\leq\frac{M\|x\|}{(n-1)!}\int\limits_{0}^{\infty}\!\!t^{n-1}e^{-({\scriptstyle{Re}}\lambda-\omega)t}\:dt=\\
&=&\frac{M\|x\|}{(n-1)!}\frac{n-1}{\mbox{\em
Re}\lambda-\omega}\int\limits_{0}^{\infty}\!\!t^{n-2}e^{-({\scriptstyle{Re}}\lambda-\omega)t}\:dt=\cdots=\frac{M\|x\|}{(\mbox{\em
Re}\lambda-\omega)^n}
\end{eqnarray*}
quels que soient $x\in{\cal E}$ et $n\in{\bf N}^*$. Par cons\'equent:
$$
\left\|{R(\lambda;A)}^n\right\|\leq\frac{M}{(\mbox{\em Re}\lambda-\omega)^n}\quad,\quad(\forall)n\in{\bf
N}^*.\fin
$$

\begin{lema}\label{num16}
Soit $A:{\cal D}(A)\subset{\cal E}\longrightarrow{\cal E}$ un
op\'erateur lin\'eaire v\'erifiant les propri\'et\'es suivantes:\\
i) $A$ est un op\'erateur ferm\'e et $\overline{{\cal D}(A)}={\cal
E}$;\\
ii) il existe $\omega\geq 0$ et $M\geq 1$ tel que
$\Lambda_\omega\subset\rho(A)$ et pour
$\lambda\in\Lambda_\omega$, on a:
$$
\left\|{R(\lambda;A)}^n\right\|\leq\frac{M}{(\mbox{Re}\lambda-\omega)^n}\quad,\quad(\forall)n\in{\bf
N}^*.
$$
Alors pour tout $\lambda\in\Lambda_\omega$, nous avons:
$$
\lim_{\scriptstyle{Re}\lambda\rightarrow\infty}\lambda R(\lambda;A)x=x\quad,\quad(\forall)x\in{\cal
E}.
$$
De plus $\lambda AR(\lambda;A)\in{\cal B(E)}$ et:
$$
\lim_{\scriptstyle{Re}\lambda\rightarrow\infty}\lambda
AR(\lambda;A)x=Ax\quad,\quad(\forall)x\in{\cal D}(A).
$$
\end{lema}
\dem
Soient $x\in{\cal D}(A)$ et $\lambda\in{\bf C}$ tel que
$\mbox{\em Re}\lambda>\omega$. Alors $R(\lambda;A)(\lambda I-A)x=x$. Si $\mbox{\em Re}\lambda\rightarrow\infty$, nous avons:
\begin{eqnarray*}
\left\|\lambda
R(\lambda;A)x-x\right\|&=&\left\|R(\lambda;A)Ax\right\|\leq\left\|R(\lambda;A)\right\|\|Ax\|\leq\\
&\leq&\frac{M}{\mbox{\em Re}\lambda-\omega}\|Ax\|\longrightarrow 0\quad,
\end{eqnarray*}
d'o\`u il r\'esulte que:
$$
\lim_{\scriptstyle{Re}\lambda\rightarrow\infty}\lambda
R(\lambda;A)x=x\quad,\quad(\forall)x\in{\cal D}(A).
$$
Soit $x\in{\cal E}$, puisque $\overline{{\cal D}(A)}={\cal E}$,
il existe une suite $(x_n)_{n\in{\bf N}}\subset{\cal D}(A)$ telle que
$x_n\longrightarrow x$ si $n\rightarrow\infty$. Nous avons:
\begin{eqnarray*}
& &\left\|\lambda R(\lambda;A)x-x\right\|\leq\\
&\leq&\left\|\lambda
R(\lambda;A)x-\lambda R(\lambda;A)x_n\right\|+\left\|\lambda
R(\lambda;A)x_n-x_n\right\|+\|x_n-x\|\leq\\
&\leq&\left\|\lambda R(\lambda;A)\right\|\|x-x_n\|+\left\|\lambda
R(\lambda;A)x_n-x_n\right\|+\|x_n-x\|\leq\\
&\leq&\frac{|\lambda|M}{\mbox{\em
Re}\lambda-\omega}\|x-x_n\|+\frac{M}{\mbox{\em
Re}\lambda-\omega}\|Ax_n\|+\|x_n-x\|=\\
&=&\frac{|\lambda|M+\mbox{\em Re}\lambda-\omega}{\mbox{\em
Re}\lambda-\omega}\:\|x_n-x\|+\frac{M}{\mbox{\em
Re}\lambda-\omega}\:\|Ax_n\|\quad.
\end{eqnarray*}
Mais $x_n\longrightarrow x$ si $n\rightarrow\infty$. Donc pour tout
$\varepsilon>0$ , il existe $n_\varepsilon\in{\bf N}$ tel
que:
$$
\|x_{n_\varepsilon}-x\|<\varepsilon\:\frac{\mbox{\em
Re}\lambda-\omega}{|\lambda|M+\mbox{\em Re}\lambda-\omega}\quad.
$$
Par cons\'equent:
$$
\left\|\lambda R(\lambda;A)x-x\right\|<\varepsilon+\frac{M}{\mbox{\em
Re}\lambda-\omega}\left\|Ax_{n_\varepsilon}\right\|\quad,
$$
d'o\`u:
$$
\limsup_{\scriptstyle{Re}\lambda\rightarrow\infty}\left\|\lambda
R(\lambda;A)x-x\right\|<\varepsilon\quad,\quad(\forall)x\in{\cal E},
$$
ou bien:
$$
\lim_{\scriptstyle{Re}\lambda\rightarrow\infty}\lambda
R(\lambda;A)x=x\quad,\quad(\forall)x\in{\cal E}.
$$
De plus:
$$
\lambda AR(\lambda;A)=\lambda\left[\lambda I-(\lambda I-A)\right]R(\lambda;A)=
\lambda\left[\lambda
R(\lambda;A)-I\right]=\lambda^2R(\lambda;A)-\lambda I.
$$
Par suite, on a:
\begin{eqnarray*}
& &\left\|\lambda AR(\lambda;A)x\right\|=\left\|\lambda\left[\lambda
R(\lambda;A)-I\right]x\right\|\leq\\
&\leq&|\lambda|\left\|\lambda
R(\lambda;A)x-x\right\|\leq|\lambda|\left(\left\|\lambda
R(\lambda;A)x\right\|+\|x\|\right)\leq\\
&\leq&|\lambda|\left(\frac{|\lambda|M}{\mbox{\em
Re}\lambda-\omega}+1\right)\|x\|\quad,\quad(\forall)x\in{\cal E}
\end{eqnarray*}
et on voit que $\lambda AR(\lambda;A)\in{\cal B(E)}$.\\
Si $x\in{\cal D}(A)$, alors nous avons:
$$
\lambda R(\lambda;A)Ax=\lambda^2R(\lambda;A)-\lambda I=\lambda
AR(\lambda;A)\quad,
$$
d'o\`u il r\'esulte que:
$$
\lim_{\scriptstyle{Re}\lambda\rightarrow\infty}\lambda
AR(\lambda;A)x=\lim_{\scriptstyle{Re}\rightarrow\infty}\lambda
R(\lambda;A)Ax=Ax\quad,\quad(\forall)x\in{\cal D}(A).\fin
$$

\begin{obs}
\em
On peut dire que les op\'erateurs born\'es $\lambda AR(\lambda;A)$
sont des\\approximations pour l'op\'erateur non born\'e $A$. C'est le motif
pour lequel on introduit la d\'efinition suivante.
\end{obs}

\begin{definitie}\label{num55}
La famille
$\left\{A_\lambda\right\}_{\lambda\in\Lambda_\omega}\subset{\cal B(E)}$, o\`u $A_\lambda=\lambda AR(\lambda;A)$, pour tout $\lambda\in\Lambda_\omega$,
s'appelle l'approximation g\'en\'eralis\'ee de Yosida de l'op\'erateur $A$.
\end{definitie}

\begin{obs}
\em
Evidemment, pour $\lambda\in\Lambda_\omega$, on voit que $A_\lambda$ est le g\'en\'erateur
infinit\'esimal d'un semi-groupe uniform\'ement continu $\left\{e^{tA_\lambda}\right\}_{t\geq 0}$. Nous utiliserons cette
famille pour montrer l'existence d'un $C_0$-semi-groupe engendr\'e par $A$.
\end{obs}

\begin{lema}\label{num17}
Soit $A:{\cal D}(A)\subset{\cal E}\longrightarrow{\cal E}$ un
op\'erateur lin\'eaire v\'erifiant les propri\'et\'es suivantes:\\
i) $A$ est un op\'erateur ferm\'e et $\overline{{\cal D}(A)}={\cal
E}$;\\
ii) il existe $\omega\geq 0$ et $M\geq 1$ tel que
$\Lambda_\omega\subset\rho(A)$ et pour
$\lambda\in\Lambda_\omega$, on a:
$$
\left\|{R(\lambda;A)}^n\right\|\leq\frac{M}{(\mbox{Re}\lambda-\omega)^n}\quad,\quad(\forall)n\in{\bf
N}^*.
$$
Si $\left\{A_\lambda\right\}_{\lambda\in{\Lambda_\omega}}$ est
l'approximation g\'en\'eralis\'ee de Yosida de l'op\'erateur $A$,
alors pour tous $\alpha,\beta\in\Lambda_\omega$ nous avons:
$$
\left\|e^{tA_\alpha}x-e^{tA_\beta}x\right\|\leq M^2te^{\omega
t}\left\|A_\alpha x-A_\beta x\right\|\quad,\quad(\forall)x\in{\cal
E}\mbox{ et }t\geq 0.
$$
\end{lema}
\dem
Soient $\alpha,\beta\in\Lambda_\omega$, $v\in[0,1]$ et $x\in{\cal
E}$. Alors:
$$
\frac{d}{dv}\left(e^{vtA_\alpha}e^{(1-v)tA_\beta}x\right)=tA_\alpha
e^{vtA_\alpha}e^{(1-v)tA_\beta}x-te^{vtA_\alpha}A_\beta
e^{(1-v)tA_\beta}x\quad.
$$
On peut facilement v\'erifier que $A_\alpha$, $A_\beta$,
$e^{vtA_\alpha}$ et $e^{(1-v)tA_\beta}$ commutent quels que soient
$\alpha,\beta\in\Lambda_\omega$ et $t\geq 0$. Nous obtenons:
\begin{eqnarray*}
&
&\int\limits_{0}^{1}\frac{d}{dv}\left(e^{vtA_\alpha}e^{(1-v)tA_\beta}x\right)\:dv=\\
&=&\int\limits_{0}^{1}\left(te^{vtA_\alpha}A_\alpha
e^{(1-v)tA_\beta}x-te^{vtA_\alpha}A_\beta
e^{(1-v)tA_\beta}x\right)\:dv\quad,
\end{eqnarray*}
d'o\`u:
$$
\left.e^{vtA_\alpha}e^{(1-v)tA_\beta}x\right|_0^1=\int\limits_{0}^{1}\left(te^{vtA_\alpha}e^{(1-v)tA_\beta}A_\alpha
x-te^{vtA_\alpha}e^{(1-v)tA_\beta}A_\beta x\right)\:dv\quad,
$$
ou bien:
$$
e^{tA_\alpha}x-e^{tA_\beta}x=t\int\limits_{0}^{1}\!\!e^{vtA_\alpha}e^{(1-v)tA_\beta}\left(A_\alpha
x-A_\beta x\right)\:dv
$$
quels que soient $t\geq 0$ et $x\in{\cal E}$. Nous en d\'eduisons que:
$$
\left\|e^{tA_\alpha}x-e^{tA_\beta}x\right\|\leq
t\int\limits_{0}^{1}\left\|e^{vtA_\alpha}\right\|\left\|e^{(1-v)tA_\beta}\right\|\left\|A_\alpha
x-A_\beta x\right\|\:dv\quad.
$$
D'autre part, nous avons:
\begin{eqnarray*}
& &\left\|e^{tA_\alpha}\right\|=\left\|e^{t\left(\alpha^2R(\alpha;A)-\alpha
I\right)}\right\|=\left\|e^{-\alpha
tI}e^{\alpha^2tR(\alpha;A)}\right\|\leq\\
&\leq& e^{-{\scriptstyle{Re}}\alpha
t}\left\|\sum\limits_{k=0}^{\infty}\frac{t^k\alpha^{2k}{R(\alpha;A)}^k}{k!}\right\|\leq
e^{-{\scriptstyle{Re}}\alpha
t}\sum\limits_{k=o}^{\infty}\frac{t^k|\alpha|^{2k}\left\|{R(\alpha;A)}^k\right\|}{k!}\leq\\
&\leq& e^{-{\scriptstyle{Re}}\alpha
t}\sum\limits_{k=0}^{\infty}\frac{t^k|\alpha|^{2k}M}{k!(\mbox{\em Re}\alpha-\omega)^k}=Me^{-{\scriptstyle{Re}}\alpha
t}\sum\limits_{k=0}^{\infty}\frac{\left(\frac{t|\alpha|^2}{{\scriptstyle{Re}}\alpha-\omega}\right)^k}{k!}=\\
&=&Me^{-{\scriptstyle{Re}}\alpha
t}e^{\frac{t|\alpha|^2}{{\scriptstyle{Re}}\alpha-\omega}}=Me^{\frac{\omega{\scriptstyle{Re}}\alpha+{\scriptstyle{Im}}^2\alpha}{{\scriptstyle{Re}}\alpha-\omega}t}\quad,
\end{eqnarray*}
quel que soient $\alpha\in\Lambda_\omega$ et $t\geq 0$. Soit $r>1$ tel
que:
$$
\frac{\omega\mbox{\em Re}\alpha+\mbox{\em Im}^2\alpha}{\mbox{\em
Re}\alpha-\omega}<\omega r\quad.
$$
Alors, nous avons:
$$
\omega\mbox{\em Re}\alpha+\mbox{\em Im}^2\alpha<\omega r\mbox{\em
Re}\alpha-\omega^2 r\quad,
$$
d'o\`u:
$$
\omega\mbox{\em Re}\alpha<\omega r\mbox{\em Re}\alpha-\omega^2 r\quad,
$$
ou bien:
$$
\omega^2 r<\omega(r-1)\mbox{\em Re}\alpha\quad.
$$
Il en d\'ecoule:
$$
\mbox{\em Re}\alpha>\frac{r}{r-1}\omega\quad.
$$
Par cons\'equent, pour tout $r>1$ et tout $\alpha\in\Lambda_\omega$ tel que
$\mbox{\em Re}\alpha>\frac{r}{r-1}\omega$, on obtient:
$$
\left\|e^{tA_\alpha}\right\|\leq Me^{r\omega
t}\quad,\quad(\forall)t\geq 0
$$
et par passage \`a limite pour $r\searrow 1$, nous obtenons:
$$
\left\|e^{tA_\alpha}\right\|\leq Me^{\omega
t}\quad,\quad(\forall)t\geq 0,
$$
pour tout $\alpha\in\Lambda_\omega$. Il vient:
\begin{eqnarray*}
\left\|e^{tA_\alpha}x-e^{tA_\beta}x\right\|&\leq&t\int\limits_{0}^{1}\!Me^{\omega
vt}Me^{\omega(1-v)t}\left\|A_\alpha x-A_\beta x\right\|\:dv=\\
&=&M^2te^{\omega t}\left\|A_\alpha x-A_\beta x\right\|
\end{eqnarray*}
quels que soient $x\in{\cal E}$ et $t\geq 0$.\fin

Maintenant nous pr\'esentons une variante du c\'el\`ebre th\'eor\`eme
de Hille - Yosida pour les semi-groupes de classe ${\cal SG}(M,\omega)$.

\begin{teo}[Hille - Yosida]
Un op\'erateur lin\'eaire:
$$
A:{\cal D}(A)\subset{\cal E}\longrightarrow{\cal E}
$$
est le g\'en\'erateur infinit\'esimal d'un
semi-groupe \semi si et seulement si:\\
i) $A$ est un op\'erateur ferm\'e et $\overline{{\cal D}(A)}={\cal
E}$;\\
ii) il existe $\omega\geq 0$ et $M\geq 1$ tel que
$\Lambda_\omega\subset\rho(A)$ et pour
$\lambda\in\Lambda_\omega$, on a:
$$
\left\|{R(\lambda;A)}^n\right\|\leq\frac{M}{(\mbox{Re}\lambda-\omega)^n}\quad,\quad(\forall)n\in{\bf
N}^*.
$$
\end{teo}
\dem
$\Longrightarrow$ On obtient cette implication en tenant
compte du th\'eor\`eme \ref{num12} et du th\'eor\`eme \ref{num15}.\\
$\Longleftarrow$  Supposons que l'op\'erateur $A:{\cal D}(A)\subset{\cal
E}\longrightarrow{\cal E}$ poss\'ede les propri\'et\'es (i) et (ii). Soit
$\left\{A_\lambda\right\}_{\lambda\in\Lambda_\omega}$,
l'approximation g\'en\'eralis\'ee de Yosida de l'op\'erateur $A$. Compte tenu du lemme
\ref{num16}, il r\'esulte que $A_\lambda\in{\cal B(E)}$ et:
$$
\lim_{\scriptstyle{Re}\lambda\rightarrow\infty}A_\lambda
x=Ax\quad,\quad(\forall)x\in{\cal D}(A).
$$
Pour $\lambda\in\Lambda_\omega$, soit
$\left\{T_\lambda(t)\right\}_{t\geq
0}=\left\{e^{tA_\lambda}\right\}_{t\geq 0}$ le semi-groupe uniform\'ement
continu engendr\'e par $A_\lambda$. Avec le lemme
\ref{num17}, on a:
$$
\left\|T_\alpha(t)x-T_\beta(t)x\right\|\leq M^2te^{\omega
t}\left\|A_\alpha x-A_\beta x\right\|\quad,\quad(\forall)\alpha,\beta\in\Lambda_\omega,\:x\in{\cal D}(A)\mbox{ et
}t\geq 0.
$$
Soient $[{\cal D}(A)]$ l'espace de Banach ${\cal D}(A)$ avec la norme
$\|\:.\:\|_{{\cal D}(A)}$, et ${\cal B}([{\cal D}(A)],{\cal E})$
l'espace des op\'erateurs lin\'eaires born\'es d\'efinis sur $[{\cal
D}(A)]$ avec valeur dans ${\cal E}$, dot\'e de la topologie forte.
Notons par ${\cal C}\left([0,\infty);{\cal B}([{\cal D}(A)],{\cal E})\right)$
l'espace des fonctions continues d\'efinies sur $[0,\infty)$ \`a
valeurs dans ${\cal B}([{\cal D}(A)],{\cal E})$ dot\'e de la topologie de la
convergence uniforme sur les intervalles compacts de $[0,\infty)$. Si
$[a,b]\subset[0,\infty)$, alors pour tout $x\in{\cal D}(A)$ nous avons:
$$
\sup_{t\in[a,b]}\left\|T_\alpha(t)x-T_\beta(t)x\right\|\leq M^2be^{\omega
b}\left(\|A_\alpha x-Ax\|+\|A_\beta x-Ax\|\right)\longrightarrow 0
$$
si $\mbox{\em Re}\alpha,\mbox{\em Re}\beta\rightarrow\infty$, d'o\`u il r\'esulte que
$\left(\left\{T_\lambda(t)\right\}_{t\geq 0}\right)_{\lambda\in\Lambda_\omega}$ est une
suite de Cauchy dans  ${\cal C}\left([0,\infty);{\cal B}([{\cal
D}(A)],{\cal E})\right)$. Donc, il existe un unique $T_0\in{\cal
C}\left([0,\infty);{\cal B}({\cal D}(A),{\cal E})\right)$ tel
que $T_\lambda(t)x\longrightarrow T_0(t)x$, si $\mbox{\em Re}\lambda\rightarrow\infty$, quel
que soit $x\in{\cal D}(A)$, pour la topologie de la convergence
uniforme sur les intervalles compacts de $[0,\infty)$. Puisque:
$$
\left\|T_\lambda(t)\right\|\leq Me^{\omega t}\quad,\quad(\forall)t\geq
0,
$$
on obtient:
$$
\left\|T_0(t)x\right\|\leq Me^{\omega
t}\|x\|\quad,\quad(\forall)t\geq 0\mbox{ et }x\in{\cal D}(A)
$$
Consid\'erons l'application lin\'eaire:
$$
\Theta_0:{\cal D}(A)\longrightarrow{\cal C}\left([a,b];{\cal E}\right)
$$
$$
\Theta_0x=T_0(\:.\:)x
$$
quel que soit $[a,b]\subset[0,\infty)$. Comme nous avons:
$$
\left\|\Theta_0x\right\|_{{\cal C}\left([a,b];{\cal
E}\right)}=\sup_{t\in[a,b]}\left\|T_0(t)x\right\|\leq Me^{\omega
b}\|x\|\leq Me^{\omega b}\|x\|_{{\cal
D}(A)}\quad,\quad(\forall)x\in{\cal D}(A),
$$
on voit que $\Theta_0$ est une application continue et puisque $\overline{{\cal D}(A)}={\cal E}$, elle se prolonge de fa\c con unique en une
application lin\'eaire continue:
$$
\Theta:{\cal E}\longrightarrow{\cal C}\left([a,b];{\cal E}\right)
$$
telle que:
$$
\left.\Theta\right|_{{\cal D}(A)}=\Theta_0
$$
et:
$$
\left\|\Theta x\right\|_{{\cal C}\left([a,b];{\cal E}\right)}\leq
Me^{\omega b}\|x\|
$$
quel que soit $x\in{\cal E}$. Par cons\'equent, il existe un seul
op\'erateur $T\in{\cal C}\left([a,b];{\cal B(E)}\right)$ tel que:
$$
\Theta x=T(\:.\:)x\quad,\quad(\forall)x\in{\cal E}.
$$
On peut r\'ep\'eter ce proc\'ed\'e pour tous les intervalles compacts
de $[0,\infty)$ et on voit qu'il existe un seul op\'erateur,
not\'e aussi par $T\in{\cal C}\left([0,\infty);{\cal B(E)}\right)$,
tel que pour tout $x\in{\cal E}$ on ait:
$$
T_\lambda(t)x\longrightarrow T(t)x\quad\mbox{si}\quad\mbox{\em Re}\lambda\rightarrow\infty,
$$
uniform\'ement par rapport \`a $t$ sur les
intervalles compacts de $[0,\infty)$. De plus:
$$
\left\|T(t)\right\|\leq Me^{\omega t}\quad,\quad(\forall)t\geq 0.
$$
Il est \'evident que:
$$
T(0)x=\lim_{{\scriptstyle{Re}}\lambda\rightarrow\infty}T_\lambda(0)x=x\quad,\quad(\forall)x\in{\cal E}
$$
et:
$$
\lim_{t\searrow 0}T(t)x=\lim_{t\searrow
0}\left(\lim_{{\scriptstyle{Re}}\lambda\rightarrow\infty}T_\lambda(t)x\right)=\lim_{{\scriptstyle{Re}}\lambda\rightarrow\infty}\left(\lim_{t\searrow
0}T_\lambda(t)x\right)=x\quad,\quad(\forall)x\in{\cal E}.
$$
Soient $t,s\in[0,\infty)$ et $x\in{\cal E}$. Alors, nous avons:
\begin{eqnarray*}
& &\left\|T(t+s)x-T(t)T(s)x\right\|\leq\left\|T(t+s)x-T_\lambda(t+s)x\right\|+\\
&+&\left\|T_\lambda(t+s)x-T_\lambda(t)T(s)x\right\|+\left\|T_\lambda(t)T(s)x-T(t)T(s)x\right\|\leq\\
&\leq&\left\|T(t+s)x-T_\lambda(t+s)x\right\|+\left\|T_\lambda(t)\right\|\left\|T_\lambda(s)x-T(s)x\right\|+\\
&+&\left\|T_\lambda(t)\left(T(s)x\right)-T(t)\left(T(s)x\right)\right\|\quad.
\end{eqnarray*}
Puisque $T_\lambda(t)\longrightarrow T(t)$, si $\mbox{\em Re}\lambda\rightarrow\infty$, pour
la topologie forte de ${\cal B(E)}$, il s'ensuit que
$T(t+s)x=T(t)T(s)x$, pour tout $x\in{\cal E}$.\\
Par cons\'equent $\left\{T(t)\right\}_{t\geq 0}\in{\cal SG}(M,\omega)$.\\
Montrons que $A$ est le g\'en\'erateur infinit\'esimal du semi-groupe
$\left\{T(t)\right\}_{t\geq 0}$.\\
Pour tout $x\in{\cal D}(A)$ on a:
\begin{eqnarray*}
& &\left\|T_\lambda(s)A_\lambda x-T(s)Ax\right\|\leq\\
&\leq&\left\|T_\lambda(s)\right\|\left\|A_\lambda x-Ax\right\|+\left\|T_\lambda(s)Ax-T(s)Ax\right\|\leq\\
&\leq& Me^{\omega
t}\left\|A_\lambda x-Ax\right\|+\left\|T_\lambda(s)Ax-T(s)Ax\right\|\longrightarrow
0
\end{eqnarray*}
si $\mbox{\em Re}\lambda\rightarrow\infty$, uniform\'ement par rapport \`a $s\in[0,t]$, d'o\`u:
$$
T(t)x-x=\lim_{{\scriptstyle{Re}}\lambda\rightarrow\infty}\left[T_\lambda(t)x-x\right]
=\lim_{{\scriptstyle{Re}}\lambda\rightarrow\infty}\int\limits_{0}^{t}\!T_\lambda(s)A_\lambda
x\:ds=\int\limits_{0}^{t}\!T(t)Ax\:ds
$$
quels que soient $x\in{\cal D}(A)$ et $t\geq 0$.\\
Soit $B$ le g\'en\'erateur infinit\'esimal du $C_0$-semigroupe
$\left\{T(t)\right\}_{t\geq 0}$. Si $x\in{\cal D}(A)$, alors:
$$
\lim_{t\searrow 0}\frac{T(t)x-x}{t}=\lim_{t\searrow
0}\frac{1}{t}\int\limits_{0}^{t}\!T(s)Ax\:ds=Ax
$$
et nous voyons que $x\in{\cal D}(B)$. Par cons\'equent ${\cal D}(A)\subset{\cal
D}(B)$ et $\left.B\right|_{{\cal D}(A)}=A$.\\
D'autre part, nous avons l'in\'egalit\'e:
$$
\left\|T(t)\right\|\leq Me^{\omega t}\quad,\quad(\forall)t\geq 0.
$$
Si $\lambda\in\Lambda_\omega$, alors $\lambda\in\rho(A)\cap\rho(B)$. Soit
$x\in{\cal D}(B)$, on a donc $\left(\lambda I-B\right)x\in{\cal E}$ et comme
l'op\'erateur $\lambda I-A:{\cal D}(A)\longrightarrow{\cal E}$ est
bijectif, il existe $x'\in{\cal D}(A)$ tel que
$\left(\lambda I-A\right)x'=\left(\lambda I-B\right)x$. Puisque
$\left.B\right|_{{\cal D}(A)}=A$, il vient que $\left(\lambda
I-B\right)x'=\left(\lambda I-B\right)x$ et comme
$\lambda\in\rho(B)$, il en r\'esulte que $x'=x$. Par suite $x\in{\cal D}(A)$ et
donc ${\cal D}(B)\subset{\cal D}(A)$.\\
Finalement on voit que ${\cal D}(A)={\cal D}(B)$ et $A=B$.\\
Nous avons montr\'e donc que A est le g\'en\'erateur infinit\'esimal
du $C_0$-semi-groupe $\left\{T(t)\right\}_{t\geq 0}$ et compte tenu
du th\'eor\`eme de l'unicit\'e de l'engendrement, il r\'esulte que
$\left\{T(t)\right\}_{t\geq 0}$ est l'unique $C_0$-semi-groupe
engendr\'e par $A$.\fin

\begin{cor}\label{num40}
Soient \semi, $A$ son g\'en\'erateur infinit\'esimal et
$\left\{A_\lambda\right\}_{\lambda\in{\Lambda_\omega}}$ l'approximation
g\'en\'eralis\'ee de Yosida de l'op\'erateur $A$. Alors:
$$
T(t)x=\lim_{{\scriptstyle{Re}}\lambda\rightarrow\infty}e^{tA_\lambda}x\quad,\quad(\forall)x\in{\cal
E},
$$
uniform\'ement par rapport \`a $t$ sur les intervalles compacts de
$[0,\infty)$.
\end{cor}
\dem
Elle r\'esulte du th\'eor\`eme de Hille-Yosida.\fin

Dans la suite nous noterons par ${\cal GI}({\cal E})$ l'ensemble des
op\'erateurs lin\'eaires qui sont des g\'en\'erateurs infinit\'esimaux de
$C_0$-semi-groupes sur l'espace de Banach ${\cal E}$. De m\^eme, pour
$\omega\geq 0$ et $M\geq 1$, nous noterons par ${\cal
GI}(M,\omega)$ l'ensemble des g\'en\'erateurs infinit\'esimaux $A\in{\cal GI}({\cal
E})$ pour lesquels:
$$
\left\|{R(\lambda;A)}^n\right\|\leq\frac{M}{(\mbox{\em
Re}\lambda-\omega)^n}\quad(\forall)\lambda\in\Lambda_\omega\mbox{ et }n\in{\bf N}^*.
$$
\vspace{2cm}

   \section{La repr\'esentation de Bromwich}

\hspace{1cm}Dans la section 1.3, avec le th\'eor\`eme \ref{num38} nous avons vu que pour les semi-groupes
uniform\'ement continus on peut obtenir une repr\'esentation par la
transform\'ee de Laplace inverse. Dans ce paragraphe nous montrerons qu'il
existe une repr\'esentation du m\^eme type pour les
$C_0$-semi-groupes. Nous commen\c cons avec quelques propri\'et\'es sur
l'approximation g\'en\'eralis\'ee de Yosida.
\begin{lema}\label{num36}
Soient \semi, $A$ son g\'en\'erateur infinit\'esimal et
$\left\{A_\mu\right\}_{\mu\in{\Lambda_\omega}}$ l'approximation
g\'en\'eralis\'ee de Yosida de l'op\'erateur $A$. Alors pour tout
$\mu\in\Lambda_\omega$, il existe $\Omega>\omega$ tel que
$\Lambda_\Omega\subset\rho(A_\mu)$
et pour tout $\lambda\in\Lambda_\Omega$ on a:
$$
\left\|R(\lambda;A_\mu)\right\|\leq\frac{M}{\mbox{\em
Re}\lambda-\Omega}\quad.
$$
De plus, pour $\varepsilon>0$, il existe une constante $C>0$ (qui
d\'epend de $M$ et $\varepsilon$) tel que:
$$
\left\|R(\lambda;A_\mu)x\right\|\leq\frac{C}{|\lambda|}\left(\|x\|+\|Ax\|\right)\quad,\quad(\forall)x\in{\cal
D}(A),
$$
quels que soient $\lambda,\mu\in{\bf C}$, avec
$\mbox{Re}\lambda>\Omega+\varepsilon$ et $\mbox{Re}\mu>\omega+\frac{|\mu|}{2}$.
\end{lema}
\dem
Soit $\mu\in\Lambda_\omega$ arbitrairement fix\'e. Nous avons vu que $A_\mu$ est le
g\'en\'erateur infinit\'esimal du semi-groupe uniform\'ement continu
$\left\{e^{tA_\mu}\right\}_{t\geq 0}$. En ce cas, nous avons:
$$
\left\|e^{tA_\mu}\right\|\leq Me^{\frac{\omega{\scriptstyle{Re}}\mu+{\scriptstyle{Im}}^2\mu}{{\scriptstyle{Re}}\mu-\omega}t}\quad,\quad(\forall)t\geq 0.
$$
Si nous notons:
$$
\Omega=\frac{\omega\mbox{\em Re}\mu+\mbox{\em Im}^2\mu}{\mbox{\em Re}\mu-\omega}\quad,
$$
alors il est clair que:
$$
\Omega=\omega+\frac{\omega^2+\mbox{\em Im}^2\mu}{\mbox{\em
Re}\mu-\omega}>\omega
$$
et que $\Lambda_\Omega=\left\{\lambda\in{\bf C}\left|\mbox{Re}\lambda>\Omega\right.\right\}\subset\rho(A_\mu)$. De plus, pour tout
$\lambda\in\Lambda_\Omega$, nous avons:
$$
\left\|R(\lambda;A_\mu)\right\|\leq\frac{M}{\mbox{\em
Re}\lambda-\Omega}\quad.
$$
Si nous consid\'erons $\lambda\in{\bf C}$ tel que $\mbox{\em
Re}\lambda>\Omega+\varepsilon$, o\`u $\varepsilon>0$, alors on voit
que:
$$
\left\|R(\lambda;A_\mu)\right\|\leq\frac{M}{\varepsilon}\quad.
$$
D'autre part, pour $x\in{\cal D}(A)$ et $\mu\in\Lambda_\omega$ tel que $\mbox{\em
Re}\mu>\omega+\frac{|\mu|}{2}$, nous obtenons:
\begin{eqnarray*}
& &\left\|A_\mu x\right\|=\left\|\mu
R(\mu;A)Ax\right\|\leq|\mu|\|R(\mu;A)\|\|Ax\|\leq\\
&\leq&|\mu|\frac{M}{\mbox{\em Re}\mu-\omega}\|Ax\|\leq 2M\|Ax\|\quad.
\end{eqnarray*}
De l'\'egalit\'e:
$$
(\lambda I-A_\mu)R(\lambda;A_\mu)=I\quad,
$$
il vient:
$$
R(\lambda;A_\mu)=\frac{1}{\lambda}I+\frac{1}{\lambda}R(\lambda;A_\mu)A_\mu
$$
et par cons\'equent:
\begin{eqnarray*}
\left\|R(\lambda;A_\mu)x\right\|&\leq&\frac{1}{|\lambda|}\left(\|x\|+\|R(\lambda;A_\mu)\|\|A_\mu
x\|\right)\leq\\
&\leq&\frac{1}{|\lambda|}\left(\|x\|+\frac{2M^2}{\varepsilon}\|Ax\|\right)\leq\\
&\leq&\frac{C}{|\lambda|}\left(\|x\|+\|Ax\|\right)\quad,\quad(\forall)x\in{\cal
D}(A),
\end{eqnarray*}
o\`u la constante $C$ ne d\'epend que de $M$ et de $\varepsilon$.\fin

\begin{lema}\label{num39}
Soient \semi, $A$ son g\'en\'erateur infinit\'esimal,
$\left\{A_\mu\right\}_{\mu\in{\Lambda_\omega}}$ l'approximation
g\'en\'eralis\'ee de
Yosida de l'op\'erateur $A$ et $\lambda\in{\bf C}$ tel que
$\mbox{Re}\lambda>\omega+\varepsilon$, arbitrairement fix\'e pour
$\varepsilon>0$. Alors:
$$
\lim_{{\scriptstyle{Re}}\mu\rightarrow\infty}R(\lambda;A_\mu)x=R(\lambda;A)x\quad,\quad(\forall)x\in{\cal
E},
$$
uniform\'ement par rapport \`a $\mbox{Im}\lambda\in[-k,k]$, o\`u
$k>0$.
\end{lema}
\dem
Compte tenu du lemme \ref{num36}, pour $\mu\in\Lambda_\omega$, il
existe $\Omega>\omega$ tel que
$\Lambda_\Omega\subset\rho(A_\mu)$. Nous avons:
$$
\Omega=\frac{\omega\mbox{\em Re}\mu+\mbox{\em Im}^2\mu}{\mbox{\em
Re}\mu-\omega}\quad.
$$
Donc l'in\'egalit\'e $\mbox{\em Re}\lambda>\Omega$ est \'equivalente
avec:
$$
\mbox{\em Re}\lambda>\omega+\frac{\omega^2+\mbox{\em Im}\mu}{\mbox{\em
Re}\mu-\omega}\quad.
$$
Soit $\varepsilon>0$. Si $\mu\in\Lambda_\omega$ tel que $\frac{\omega^2+{\scriptstyle{Im}}\mu}{{\scriptstyle{Re}}\mu-\omega}<\varepsilon$, alors $\mbox{\em Re}\lambda>\omega+\varepsilon$
implique $\mbox{\em Re}\lambda>\Omega$. Par suite,
$\lambda\in\rho(A_\mu)$. Donc il existe $R(\lambda;A_\mu)$ et avec le
lemme \ref{num36} on voit que:
$$
\left\|R(\lambda;A_\mu)\right\|\leq\frac{M}{\mbox{\em
Re}\lambda-\omega}\quad.
$$
D'autre part, nous avons:
\begin{eqnarray*}
\mbox{\em Re}\frac{\lambda\mu}{\lambda+\mu}&=&\mbox{\em
Re}\left(\lambda-\frac{\lambda^2}{\lambda+\mu}\right)=\mbox{\em
Re}\lambda-\mbox{\em Re}\frac{\lambda^2}{\lambda+\mu}>\\
&>&\omega+\varepsilon-\mbox{\em Re}\frac{\lambda^2}{\lambda+\mu}\quad.
\end{eqnarray*}
Etant donn\'e $k>0$ tel que $|\mbox{\em Im}\lambda|\leq k$, il
existe $\mu\in\Lambda_\omega$ tel que $\mbox{\em
Re}\frac{\lambda^2}{\lambda+\mu}<\frac{\varepsilon}{2}$. Il s'ensuit
que $\mbox{\em
Re}\frac{\lambda\mu}{\lambda+\mu}>\omega+\frac{\varepsilon}{2}$. Par
cons\'equent, $\frac{\lambda\mu}{\lambda+\mu}\in\rho(A)$ et donc $R\left(\frac{\lambda\mu}{\lambda+\mu};A\right)$ existe bien. Nous avons:
\begin{eqnarray*}
& &\frac{1}{\lambda+\mu}(\lambda I-A_\mu)(\mu
I-A)R\left(\frac{\lambda\mu}{\lambda+\mu};A\right)=\\
&=&\frac{1}{\lambda+\mu}\left[\lambda I-\mu^2R(\mu;A)+\mu I\right](\mu
I-A)R\left(\frac{\lambda\mu}{\lambda+\mu};A\right)=\\
&=&\left(\mu
I-A-\frac{\mu^2}{\lambda+\mu}I\right)R\left(\frac{\lambda\mu}{\lambda+\mu};A\right)=\\
&=&\left(\frac{\lambda\mu}{\lambda+\mu}I-A\right)R\left(\frac{\lambda\mu}{\lambda+\mu};A\right)=I\quad.
\end{eqnarray*}
Par un calcul analogue, on peut obtenir:
$$
\frac{1}{\lambda+\mu}(\mu
I-A)R\left(\frac{\lambda\mu}{\lambda+\mu};A\right)(\lambda
I-A_\mu)=I\quad.
$$
Il s'ensuit que:
$$
R(\lambda;A_\mu)=\frac{1}{\lambda+\mu}(\mu
I-A)R\left(\frac{\lambda\mu}{\lambda+\mu};A\right)\quad.
$$
Par cons\'equent:
\begin{eqnarray*}
& &R(\lambda;A_\mu)-R(\lambda;A)=\frac{1}{\lambda+\mu}(\mu
I-A)R\left(\frac{\lambda\mu}{\lambda+\mu};A\right)-R(\lambda;A)=\\
&=&\frac{1}{\lambda+\mu}\left[(\mu
I-A)R\left(\frac{\lambda\mu}{\lambda+\mu};A\right)-(\lambda+\mu)R(\lambda;A)\right]=\\
&=&\left.\frac{1}{\lambda+\mu}(\mu
I-A)R\left(\frac{\lambda\mu}{\lambda+\mu};A\right)\right[(\mu
I-A)(\lambda I-A)-\\
&-&\left.(\lambda+\mu)\left(\frac{\lambda\mu}{\lambda+\mu}I-A\right)\right]R(\mu;A)R(\lambda;A)=\\
&=&\frac{1}{\lambda+\mu}(\mu
I-A)R\left(\frac{\lambda\mu}{\lambda+\mu};A\right)A^2R(\mu;A)R(\lambda;A)=\\
&=&\frac{1}{\lambda+\mu}R\left(\frac{\lambda\mu}{\lambda+\mu};A\right)R(\lambda;A)A^2\quad.
\end{eqnarray*}
Mais:
$$
\left\|R(\lambda;A)\right\|\leq\frac{M}{\varepsilon}
$$
et:
$$
\left\|R\left(\frac{\lambda\mu}{\lambda+\mu};A\right)\right\|\leq\frac{2M}{\varepsilon}\quad.
$$
Si $x\in{\cal D}(A^2)$, alors on voit que:
\begin{eqnarray*}
& &\left\|R(\lambda;A_\mu)x-R(\lambda;A)x\right\|\leq\\
&\leq&\frac{1}{|\lambda+\mu|}\left\|R\left(\frac{\lambda\mu}{\lambda+\mu};A\right)\right\|\|R(\lambda;A)\|\|A^2x\|\leq\\
&\leq&\frac{1}{|\mu|}\frac{2M}{\varepsilon}\frac{M}{\varepsilon}\|A^2x\|\leq\frac{1}{\mbox{\em
Re}\mu}\frac{2M^2}{\varepsilon^2}\|A^2x\|\quad.
\end{eqnarray*}
Il s'ensuit que:
$$
\lim_{{\scriptstyle{Re}}\mu\rightarrow\infty}R(\lambda;A_\mu)x=R(\lambda;A)x\quad,\quad(\forall)x\in{\cal
D}(A^2),
$$
uniform\'ement par rapport \`a $\mbox{\em Im}\lambda\in[-k,k]$, o\`u
$k>0$. Avec le th\'eor\`eme \ref{num37}, on sait que $\overline{{\cal
D}(A^2)}={\cal E}$. Comme $R(\lambda;A)$ et $R(\lambda;A_\mu)$ sont
uniform\'ement born\'es, on obtient:
$$
\lim_{{\scriptstyle{Re}}\mu\rightarrow\infty}R(\lambda;A_\mu)x=R(\lambda;A)x\quad,\quad(\forall)x\in{\cal
E},
$$
uniform\'ement par rapport \`a $\mbox{\em Im}\lambda\in[-k,k]$, o\`u
$k>0$.\fin

\begin{teo}\label{num41}
Soit $A$ le g\'en\'erateur infinit\'esimal du semi-groupe \semi et
$\lambda\in\Lambda_\omega$. Alors pour tout $x\in{\cal D}(A)$ on a:
$$
\int\limits_{0}^{t}\!T(s)x\:ds=\int\limits_{{\scriptstyle{Re}}\lambda-i\infty}^{{\scriptstyle{Re}}\lambda-i\infty}\!e^{zt}R(z;A)x\:\frac{dz}{z}
$$
et l'int\'egrale de la partie droite de l'\'egalit\'e est
uniform\'ement convergente par\\rapport \`a $t$ sur les intervalles
compacts de $]0,\infty)$.
\end{teo}
\dem
Soit $\left\{A_\mu\right\}_{\mu\in{\Lambda_\omega}}$ l'approximation
g\'en\'eralis\'ee de Yosida de l'op\'erateur $A$. Soit $\mu\in\Lambda_\omega$ tel que $Re\:\mu>\omega+\frac{|\mu|}{2}$. Avec le lemme \ref{num36}, nous
d\'eduissons qu'il existe
$\Omega=\frac{\omega{\scriptstyle{Re}}\mu+{\scriptstyle{Im}}^2\mu}{{\scriptstyle{Re}}\mu-\omega}>\omega$
tel que $\Lambda_\Omega=\{\lambda\in{\bf C}|Re\:\lambda>\Omega\}\subset\rho(A_\mu)$.
Soit $\lambda\in\Lambda_\Omega$. En utilisant le th\'eor\`eme \ref{num38}, pour $R>2Re\:\lambda$ on peut consid\'erer le contour de Jordan $A_\mu$-spectral
$$
\Gamma^1_R=\Gamma^{1^{'}}_R\cup\Gamma^{1^{"}}_R\quad,
$$
o\`u
$$
\Gamma^{1^{'}}_R=\{Re\:\lambda+i\tau|\tau\in[-R,R]\}
$$
et
$$
\Gamma^{1^{"}}_R=\left\{Re\:\lambda+R(\cos\varphi+i\sin\varphi)\left|\varphi\in\left[\frac{\pi}{2},\frac{3\pi}{2}\right]\right.\right\}\quad.
$$
Pour le semi-groupe uniform\'ement continu $\left\{e^{tA_\mu}\right\}_{t\geq 0}$ engendr\'e par $A_\mu$ il en r\'esulte:
$$
e^{tA_\mu}=\lim_{R\rightarrow\infty}\frac{1}{2\pi i}\int\limits_{Re\:\lambda-iR}^{Re\:\lambda+iR}\!e^{zt}R(z;A_\mu)\:dz=\frac{1}{2\pi i}\int\limits_{Re\:\lambda-i\infty}^{Re\:\lambda+i\infty}\!e^{zt}R(z;A_\mu)\:dz\quad,
$$
uniform\'ement par rapport \`a $t$ sur les intervalles compacts de $[0,\infty)$.
Pour $R>2Re\:\lambda$ et $x\in{\cal D}(A)$ nous notons
$$
I_R(s)=\frac{1}{2\pi i}\int\limits_{Re\:\lambda-iR}^{Re\:\lambda+iR}\!e^{zs}R(z;A_\mu)x\:dz\quad.
$$
Soient $0<a<b$. Pour tout $t\in[a,b]$ nous obtenons:
\begin{eqnarray*}
& &\int\limits_{0}^{t}\!I_R(s)\:ds=\frac{1}{2\pi
i}\int\limits_{0}^{t}\!\int\limits_{Re\:\lambda-iR}^{Re\:\lambda+iR}\!e^{zs}R(z;A_\mu)x\:dz\:ds=\\
&=&\frac{1}{2\pi
i}\int\limits_{Re\:\lambda-iR}^{Re\:\lambda+iR}\!\int\limits_{0}^{t}\!e^{zs}\:dsR(z;A_\mu)x\:dz=\\
&=&\frac{1}{2\pi i}\int\limits_{Re\:\lambda-iR}^{Re\:\lambda+iR}\!e^{zt}R(z;A_\mu)x\:\frac{dz}{z}-\frac{1}{2\pi
i}\int\limits_{Re\:\lambda-iR}^{Re\:\lambda+iR}\!R(z;A_\mu)x\:\frac{dz}{z}\quad.
\end{eqnarray*}
Montrons que pour l'int\'egrale
$$
I(R)=\frac{1}{2\pi
i}\int\limits_{Re\:\lambda-iR}^{Re\:\lambda+iR}\!R(z;A_\mu)x\:\frac{dz}{z}\quad,
$$
on a
$$
\lim_{R\rightarrow\infty}I(R)=0\quad.
$$
Soit le contour de Jordan lisse et ferm\'e
$$
\Gamma^2_R=\Gamma^{2^{'}}_R\cup\Gamma^{2^{"}}_R\quad,
$$
o\`u
$$
\Gamma^{2^{'}}_R=\{Re\:\lambda+i\tau|\tau\in[-R,R]\}
$$
et
$$
\Gamma^{2^{"}}_R=\left\{Re\:\lambda+R(\cos\varphi+i\sin\varphi)\left|\varphi\in\left[-\frac{\pi}{2},\frac{\pi}{2}\right]\right.\right\}\quad.
$$
Avec le th\'eor\`eme de Cauchy (\cite[pag. 225]{dunford-schwartz}), on voit que
$$
\frac{1}{2\pi
i}\int\limits_{\Gamma^2_R}\!R(z;A_\mu)x\:\frac{dz}{z}=0\quad,
$$
ou bien
$$
\frac{1}{2\pi
i}\int\limits_{\Gamma^{2^{'}}_R}\!R(z;A_\mu)x\:\frac{dz}{z}+\frac{1}{2\pi
i}\int\limits_{\Gamma^{2^{"}}_R}\!R(z;A_\mu)x\:\frac{dz}{z}=0\quad.
$$
Soit $z\in\Gamma^{2^{"}}$. Compte tenu du lemme \ref{num36}, il existe $C>0$ tel que
$$
\left\|R(z;A_\mu)x\right\|\leq\frac{C}{|z|}\left(\|x\|+\|Ax\|\right)\quad.
$$
De plus, pour $z\in\Gamma^2{^{"}}_R$ on a:
\begin{eqnarray*}
|z|&=&|Re\:\lambda+R(\cos\varphi+isin\varphi)|=|Re\:\lambda-[-R(\cos\varphi+i\sin\varphi)]|\geq\\
&\geq&|\:|Re\:\lambda|-|-R(\cos\varphi+i\sin\varphi)|\:|=|Re\:\lambda-R|=R-Re\:\lambda\quad,
\end{eqnarray*}
d'o\`u il r\'esulte
$$
\frac{1}{|z|}\leq\frac{1}{R-Re\:\lambda}\quad.
$$
Par cons\'equent, on a:
\begin{eqnarray*}
& &\left\|\frac{1}{2\pi
i}\int\limits_{\Gamma^{2^{"}}_R}\!R(z;A_\mu)x\:\frac{dz}{z}\right\|\leq\frac{1}{2\pi}\int\limits_{\Gamma^{2^{"}}_R}\!\|R(z;A_\mu)x\|\:\frac{|dz|}{|z|}\leq\\
&\leq&\frac{1}{2\pi}\int\limits_{\Gamma^{2^{"}}_R}\!\frac{C}{|z|}(\|x\|+\|Ax\|)\frac{1}{|z|}|dz|\leq\frac{C}{2\pi}\frac{\|x\|+\|Ax\|}{(R-Re\:\lambda)^2}\int\limits_{\Gamma^{2^{"}}_R}\!|dz|=\\
&=&\frac{C}{2}\frac{R}{(R-Re\:\lambda)^2}(\|x\|+\|Ax\|)\quad.
\end{eqnarray*}
Il s'ensuit donc que:
$$
\lim_{R\rightarrow\infty}\frac{1}{2\pi
i}\int\limits_{\Gamma^{2^{"}}_R}\!R(z;A_\mu)x\:\frac{dz}{z}=0\quad.
$$
Par suite, on a:
$$
\lim_{R\rightarrow\infty}\frac{1}{2\pi
i}\int\limits_{\Gamma^{2^{'}}_R}\!R(z;A_\mu)x\:\frac{dz}{z}=0
$$
ou bien
$$
\lim_{R\rightarrow\infty}I(R)=0\quad.
$$
Alors nous avons:
$$
\lim_{R\rightarrow\infty}\int\limits_{0}^{t}\!I_R(s)\:ds=\lim_{R\rightarrow\infty}\frac{1}{2\pi i}\int\limits_{Re\:\lambda-iR}^{Re\:\lambda+iR}\!e^{zt}R(z;A_\mu)x\:\frac{dz}{z}
$$
d'o\`u
$$
\int\limits_{0}^{t}\!I_R(s)\:ds=\frac{1}{2\pi i}\int\limits_{Re\:\lambda-i\infty}^{Re\:\lambda+i\infty}\!e^{zt}R(z;A_\mu)x\:\frac{dz}{z}\quad.
$$
Avec le corollaire \ref{num40}
et le lemme \ref{num39}, on obtient:
\begin{eqnarray*}
& &\int\limits_{0}^{t}\!T(s)x\:ds=\lim_{{\scriptstyle{Re}}\mu\rightarrow\infty}\int\limits_{0}^{t}\!e^{sA_\mu}x\:ds=\\
&=&\lim_{{\scriptstyle{Re}}\mu\rightarrow\infty}\frac{1}{2\pi
i}\int\limits_{Re\:\lambda-i\infty}^{Re\:\lambda+i\infty}\!e^{zt}R(z;A_\mu)x\:\frac{dz}{z}=\\
&=&\frac{1}{2\pi
i}\int\limits_{Re\:\lambda-i\infty}^{Re\:\lambda+i\infty}\!e^{zt}R(z;A)x\:\frac{dz}{z}
\end{eqnarray*}
et comme:
$$
\lim_{{\scriptstyle{Re}}\mu\rightarrow\infty}\Omega=\lim_{{\scriptstyle{Re}}\mu\rightarrow\infty}\frac{\omega\mbox{\em
Re}\mu+\mbox{\em
Im}^2\mu}{\mbox{\em Re}\mu-\omega}=
\omega\quad,
$$
nous obtenons le r\'esultat d\'esir\'e.\fin

\begin{teo}[Bromwich]\label{num47}
Soit $A$ le g\'en\'erateur infinit\'esimal d'un semi-groupe \semi et
$\lambda\in\Lambda_\omega$. Alors:
$$
T(t)x=\frac{1}{2\pi
i}\int\limits_{{\scriptstyle{Re}}\lambda-i\infty}^{{\scriptstyle{Re}}\lambda+i\infty}\!e^{zt}R(z;A)x\:dz\quad,\quad(\forall)x\in{\cal D}(A^2)
$$
et pour tout $\delta>0$, l'int\'egrale est uniform\'ement convergente
par rapport \`a $t\in\left[\delta,\frac{1}{\delta}\right]$.
\end{teo}
\dem
Si $x\in{\cal D}(A^2)$, alors $Ax\in{\cal D}(A)$. Compte tenu du
th\'eor\`eme \ref{num41}, on voit que:
$$
\int\limits_{0}^{t}\!T(s)Ax\:ds=\frac{1}{2\pi i}\int\limits_{{\scriptstyle{Re}}\lambda-i\infty}^{{\scriptstyle{Re}}\lambda+i\infty}\!e^{zt}R(z;A)Ax\:\frac{dz}{z}\quad.
$$
D'o\`u il r\'esulte que:
$$
T(s)x-x=\frac{1}{2\pi i}\int\limits_{{\scriptstyle{Re}}\lambda-i\infty}^{{\scriptstyle{Re}}\lambda+i\infty}\!e^{zt}R(z;A)Ax\:\frac{dz}{z}\quad.
$$
De l'\'egalit\'e:
$$
R(z;A)(zI-A)=I\quad,
$$
nous d\'eduisons:
$$
\frac{1}{z}R(z;A)A=R(z;A)-\frac{1}{z}I
$$
et par suite:
$$
T(s)x-x=\frac{1}{2\pi i}\int\limits_{{\scriptstyle{Re}}\lambda-i\infty}^{{\scriptstyle{Re}}\lambda+i\infty}\!e^{zt}R(z;A)x\:dz-\frac{1}{2\pi i}\int\limits_{{\scriptstyle{Re}}\lambda-i\infty}^{{\scriptstyle{Re}}\lambda+i\infty}\!e^{zt}x\:\frac{dz}{z}\quad.
$$
Compte tenu que:
$$
\frac{1}{2\pi i}\int\limits_{{\scriptstyle{Re}}\lambda-i\infty}^{{\scriptstyle{Re}}\lambda+i\infty}\!e^{zt}x\:\frac{dz}{z}=x
$$
et que pour tout $\delta>0$, l'int\'egrale est uniform\'ement convergente
par rapport \`a $t\in\left[\delta,\frac{1}{\delta}\right]$, nous
obtenons l'\'egalit\'e de l'\'enonc\'e.\fin
\vspace{2cm}

    \section{Conditions suffisantes d'appartenances \`a ${\cal GI}(M,0)$}

\hspace{1cm}Nous pr\'esentons dans la suite deux conditions suffisantes pour qu'un
op\'erateur soit le g\'en\'erateur infinit\'esimal d'un
$C_0$-semi-groupe uniform\'ement born\'e.

\begin{teo}\label{num45}
Soit $A$ un op\'erateur lin\'eaire ferm\'e d\'efini sur un sous
espace dense de ${\cal E}$ et v\'erifiant les propri\'et\'es
suivantes:\\
i) il existe $\delta\in\left]0,\frac{\pi}{2}\right[$ tel que:
$$
\rho(A)\supset\Sigma_\delta=\left\{\lambda\in{\bf C}\left||\arg
z|<\frac{\pi}{2}+\delta\right.\right\}\cup\{0\};
$$
ii) il existe une constante $K>1$ tel que:
$$
\|R(\lambda;A)\|\leq
\frac{K}{|\lambda|}\quad,\quad(\forall)\lambda\in\Sigma_\delta-\{0\}.
$$
Alors $A$ est le g\'en\'erateur infinit\'esimal d'un semi-groupe
$\left\{T(t)\right\}_{t\geq 0}$ pour lequel il existe $M>1$ tel que $\|T(t)\|\leq M$, pour tout $t\geq 0$.\\ 
De plus, pour tout
$\nu\in\left]\frac{\pi}{2},\frac{\pi}{2}+\delta\right[$ et
$\Gamma_\nu=\Gamma_\nu^{(1)}\cup\Gamma_\nu^{(2)}$, o\`u:
$$
\Gamma_\nu^{(1)}=\left\{\left.r(\cos\nu-i\sin\nu)\right|\:r\in[0,\infty)\right\}
$$
et:
$$
\Gamma_\nu^{(2)}=\left\{\left.r(\cos\nu+i\sin\nu)\right|\:r\in[0,\infty)\right\}\quad,
$$
on a:
$$
T(t)=\frac{1}{2\pi i}\int\limits_{\Gamma_\nu}\!e^{zt}R(z;A)\:dz
$$
et l'int\'egrale est uniform\'ement convergente par rapport \`a $t>0$.
\end{teo}
\dem
Soit $\delta\in\left]0,\frac{\pi}{2}\right[$. Pour $\nu\in\left]\frac{\pi}{2},\frac{\pi}{2}+\delta\right[$ consid\'erons le chemin d'int\'egration $\Gamma_\nu=\Gamma_\nu^{(1)}\cup\Gamma_\nu^{(2)}$, o\`u:
$$
\Gamma_\nu^{(1)}=\left\{\left.r(\cos\nu-i\sin\nu)\right|\:r\in[0,\infty)\right\}
$$
et:
$$
\Gamma_\nu^{(2)}=\left\{\left.r(\cos\nu+i\sin\nu)\right|\:r\in[0,\infty)\right\}\quad.
$$
Soit:
$$
U(t)=\frac{1}{2\pi i}\int\limits_{\Gamma_\nu}\!e^{zt}R(z;A)\:dz\quad.
$$
Compte tenu du (ii), on voit que l'int\'egrale est uniform\'ement
convergente par rapport \`a $t>0$. Pour $R>0$, nous d\'efinissons le contour de Jordan lisse et ferm\'e $\Gamma_\nu=\Gamma_{R,\nu}^{(1)}\cup\Gamma_{R,\nu}^{(2)}\cup\Gamma_{R,\nu}^{(3)}$ o\`u
$$
\Gamma_{R,\nu}^{(1)}=\left\{\left.r(\cos\nu-i\sin\nu)\right|\:r\in[0,R]\right\}\quad,
$$
$$
\Gamma_{R,\nu}^{(2)}=\left\{\left.R(\cos\nu+i\sin\nu)\right|\:\theta\in[-\nu,\nu]\right\}\quad,
$$
et
$$
\Gamma_{R,\nu}^{(3)}=\left\{\left.r(\cos\nu+i\sin\nu)\right|\:r\in[0,R]\right\}\quad.
$$
D'apr\'es le th\'eor\`eme de Cauchy (\cite[pag. 225]{dunford-schwartz}), on a
$$
\int\limits_{\Gamma_{R,\nu}}\!e^{zt}R(z;A)\:dz=0\quad.
$$
Par cons\'equent, on peut changer le chemin d'int\'egration $\Gamma_\nu$ par $\Gamma_t=\Gamma_t^{(1)}\cup\Gamma_t^{(2)}\cup\Gamma_t^{(3)}$ o\`u
$$
\Gamma_t^{(1)}=\left\{r(\cos\nu-i\sin\nu)\left|\:r\in\left[\frac{1}{t},\infty\right)\right.\right\}\quad,
$$
$$
\Gamma_t^{(2)}=\left\{\frac{1}{t}(\cos\theta+i\sin\theta)\left|\:\theta\in[-\nu,\nu]\right.\right\}\quad
$$
et
$$
\Gamma_t^{(3)}=\left\{r(\cos\nu+i\sin\nu)\left|\:r\in\left[\frac{1}{t},\infty\right)\right.\right\}\quad.
$$
Alors:
$$
U(t)=\frac{1}{2\pi i}\int\limits_{\Gamma_t^{(1)}}\!e^{zt}R(z;A)\:dz+\frac{1}{2\pi i}\int\limits_{\Gamma_t^{(2)}}\!e^{zt}R(z;A)\:dz+\frac{1}{2\pi i}\int\limits_{\Gamma_t^{(3)}}\!e^{zt}R(z;A)\:dz
$$
et si nous notons
$$
U_1(t)=\frac{1}{2\pi i}\int\limits_{\Gamma_t^{(1)}}\!e^{zt}R(z;A)\:dz\quad,
$$
$$
U_2(t)=\frac{1}{2\pi i}\int\limits_{\Gamma_t^{(2)}}\!e^{zt}R(z;A)\:dz
$$
et
$$
U_3(t)=\frac{1}{2\pi i}\int\limits_{\Gamma_t^{(3)}}\!e^{zt}R(z;A)\:dz\quad,
$$
il vient:
$$
\left\|U(t)\right\|\leq\left\|U_1(t)\right\|+\left\|U_2(t)\right\|+\left\|U_3(t)\right\|\quad.
$$
Comme $\nu\in\left]\frac{\pi}{2},\frac{\pi}{2}+\delta\right[$, on en d\'eduit que $\cos\nu<0$. 
Avec le changement de variable
$$
z=r(\cos\nu-i\sin\nu)\quad,\quad r\in\left[\frac{1}{t},\infty\right),
$$
nous avons:
\begin{eqnarray*}
& &\left\|U_1(t)\right\|=\left\|\frac{1}{2\pi i}\int\limits_{\Gamma_t^{(1)}}\!e^{zt}R(z;A)\:dz\right\|=\\
&=&\left\|\frac{1}{2\pi i}\int\limits_{\frac{1}{t}}^{\infty}\!e^{rt(\cos\nu-i\sin\nu)}R(r(\cos\nu-i\sin\nu);A)(\cos\nu-isin\nu)\:dr\right\|\leq\\
&\leq&\frac{K}{2\pi}\int\limits_{\frac{1}{t}}^{\infty}\!\!e^{rt\cos\nu}\frac{K}{r}\:dr=\frac{K}{2\pi}\int\limits_{\frac{1}{t}}^{\infty}\!\!e^{-rt(-\cos\nu)}\frac{1}{r}\:dr\quad.
\end{eqnarray*}
Pour 
$$
s=rt(-\cos\nu)\quad,\quad s\in[-\cos\nu,\infty),
$$
il vient $ds=-t\cos\nu dr$. Donc:
\begin{eqnarray*}
\left\|U_1(t)\right\|&\leq&\frac{K}{2\pi}\int\limits_{-\cos\nu}^{\infty}\!\!e^{-s}\frac{-t\cos\nu}{s}\frac{1}{-t\cos\nu}\:ds=\frac{K}{2\pi}\int\limits_{-\cos\nu}^{\infty}\!\!\frac{e^{-s}}{s}\:ds\leq\\
&\leq&\frac{K}{2\pi}\int\limits_{-\cos\nu}^{\infty}\frac{1}{s^2}\:ds=\frac{K}{2\pi}\left(\frac{1}{-\cos\nu}\right)=M^{'}\quad.
\end{eqnarray*}
De fa\c con analogue, nous obtenons:
\begin{eqnarray*}
& &\left\|U_3(t)\right\|=\left\|\frac{1}{2\pi i}\int\limits_{\Gamma_t^{(3)}}\!e^{zt}R(z;A)\:dz\right\|=\\
&=&\left\|\frac{1}{2\pi i}\int\limits_{\frac{1}{t}}^{\infty}\!e^{rt(\cos\nu+i\sin\nu)}R(r(\cos\nu+i\sin\nu);A)(\cos\nu+isin\nu)\:dr\right\|\leq\\
&\leq&\frac{K}{2\pi}\int\limits_{\frac{1}{t}}^{\infty}\!\!e^{rt\cos\nu}\frac{K}{r}\:dr\leq M^{'}\quad.
\end{eqnarray*}
De m\^eme, pour l'int\'egrale $U_2(t)$, avec le changement de variable
$$
z=\frac{1}{t}(\cos\theta+i\sin\theta)\quad,\quad\theta\in[-\nu,\nu],
$$
on a
$$
dz=\frac{1}{t}(-\sin\theta+i\cos\theta)\:d\theta
$$
et
\begin{eqnarray*}
& &\left\|U_2(t)\right\|=\left\|\frac{1}{2\pi i}\int\limits_{\Gamma_t^{(2)}}\!e^{zt}R(z;A)\:dz\right\|=\\
&=&\left\|\frac{1}{2\pi i}\int\limits_{-\nu}^{\nu}\!e^{t\frac{1}{t}(\cos\theta+i\sin\theta)}R\left(\frac{1}{t}(\cos\theta+i\sin\theta);A\right)\frac{1}{t}(-sin\theta+i\cos\theta)\:d\theta\right\|\leq\\
&\leq&\frac{1}{2\pi}\int\limits_{-\nu}^{\nu}\!\!e^{\cos\theta}\frac{K}{\frac{1}{t}}\frac{1}{t}\:d\theta=\frac{K}{2\pi}\int\limits_{-\nu}^{\nu}\!\!e^{\cos\theta}\:d\theta\leq\frac{Ke}{2\pi}\int\limits_{-\nu}^{\nu}\:d\theta=M^{"}\quad.
\end{eqnarray*}
Par cons\'equent, il existe $M\geq 1$ tel que:
$$
\|U(t)\|\leq M\quad,\quad(\forall)t\geq 0.
$$
Nous allons maintenant montrer que pour tout $\lambda\in\Lambda_0=\left\{\lambda\in{\bf C}\left|\mbox{Re}\lambda>0\right.\right\}$, on
a:
$$
R(\lambda;A)=\int\limits_{0}^{\infty}\!e^{-\lambda t}U(t)\:dt\quad.
$$
Nous avons successivement:
\begin{eqnarray*}
& &\int\limits_{0}^{\tau}\!e^{-\lambda t}U(t)\:dt=\frac{1}{2\pi
i}\int\limits_{0}^{\tau}\int\limits_{\Gamma_\nu}\!\!e^{-(\lambda-z)t}R(z;A)\:dz
dt=\\
&=&\frac{1}{2\pi
i}\int\limits_{\Gamma_\nu}\int\limits_{0}^{\tau}\!e^{(z-\lambda)t}\:dtR(z;A)\:dz=\frac{1}{2\pi
i}\int\limits_{\Gamma_\nu}\frac{e^{(z-\lambda)\tau}-1}{z-\lambda}R(z;A)\:dz=\\
&=&\frac{1}{2\pi
i}\int\limits_{\Gamma_\nu}\frac{e^{(z-\lambda)\tau}}{z-\lambda}R(z;A)\:dz+\frac{1}{2\pi
i}\int\limits_{\Gamma_\nu}\frac{R(z;A)}{\lambda-z}\:dz=\\
&=&\frac{1}{2\pi
i}\int\limits_{\Gamma_\nu}\frac{e^{(z-\lambda)\tau}}{z-\lambda}R(z;A)\:dz+(\lambda I-A)^{-1}=\\
&=&\frac{1}{2\pi
i}\int\limits_{\Gamma_\nu}\frac{e^{(z-\lambda)\tau}}{z-\lambda}R(z;A)\:dz+R(\lambda;A)\quad.
\end{eqnarray*}
Par cons\'equent:
\begin{eqnarray*}
& &\left\|\int\limits_{0}^{\tau}\!e^{-\lambda t}U(t)\:dt-R(\lambda;A)\right\|=\left\|\frac{1}{2\pi
i}\int\limits_{\Gamma_\nu}\frac{e^{(z-\lambda)\tau}}{z-\lambda}R(z;A)\:dz\right\|\leq\\
&\leq&\frac{1}{2\pi}\int\limits_{\Gamma_\nu}\frac{|e^{(z-\lambda)\tau}|}{|z-\lambda|}\|R(z;A)\|\:|dz|\leq
\frac{1}{2\pi}\int\limits_{\Gamma_\nu}\!\frac{e^{({\scriptstyle{Re}}z-{\scriptstyle{Re}}\lambda)\tau}}{|z-\lambda|}\frac{K}{|z|}|dz|=\\
&=&
\frac{K}{2\pi}e^{-\tau{\scriptstyle{Re}}\lambda}\int\limits_{\Gamma_\nu}\!e^{\tau{\scriptstyle{Re}}}\frac{1}{|z-\lambda|}\frac{1}{|z|}|dz|\quad.
\end{eqnarray*}
Soit
$$
C_\tau=\frac{K}{2\pi}\int\limits_{\Gamma_\nu}\!\!e^{\tau{\scriptstyle{Re}}z}\frac{|dz|}{|z||z-\lambda|}\quad.
$$
Pour $z\in\Gamma_\nu$ on a $|z|=r$ et
$$
|z-\lambda|\geq|\:|z|-|\lambda|\:|=|r-|\lambda|\:|\quad.
$$
Donc:
$$
C_\tau\leq\frac{K}{2\pi}\int\limits_{0}^{\infty}\!e^{\tau{\scriptstyle{Re}}z}\frac{1}{|r-|\lambda||}\frac{1}{r}\:dr=\frac{K}{2\pi}e^{\tau{\scriptstyle{Re}}z}\int\limits_{0}^{\infty}\frac{1}{r|r-|\lambda||}\:dr<\infty\quad.
$$
Si nous notons
$$
C=\sup_{\tau\geq 0}C_\tau\quad,
$$
alors nous obtenons:
$$
\left\|\int\limits_{0}^{\tau}\!e^{-\lambda t}U(t)\:dt-R(\lambda;A)\right\|\leq Ce^{-\tau{\scriptstyle{Re}}\lambda}\quad.
$$
En passant \`a limite pour $\tau\rightarrow\infty$, on obtient
$$
R(\lambda;A)=\int\limits_{0}^{\infty}\!e^{-\lambda t}U(t)\:dt
$$
pour tout $\lambda\in\Lambda_0$. Comme $\|U(t)\|\leq M$, pour tout
$t\geq 0$, par r\'ecurrence on peut obtenir:
$$
\frac{d^{n-1}}{d\lambda^{n-1}}R(\lambda;A)=(-1)^{n-1}\int\limits_{0}^{\infty}\!t^{n-1}e^{-\lambda
t}U(t)\:dt\quad,\quad(\forall)n\in{\bf N}^*.
$$
Mais avec la proposition \ref{num14} (iii), on voit que:
$$
\frac{d^{n-1}}{d\lambda^{n-1}}R(\lambda;A)=(-1)^{n-1}(n-1)!{R(\lambda;A)}^n\quad,\quad(\forall)n\in{\bf
N}^*.
$$
Par cons\'equent, nous obtenons:
\begin{eqnarray*}
&
&\left\|{R(\lambda;A)}^n\right\|=\left\|\frac{1}{(n-1)!}\int\limits_{0}^{\infty}\!t^{n-1}e^{-\lambda
t}U(t)\:dt\right\|\leq\\
&\leq&\frac{M}{(n-1)!}\int\limits_{0}^{\infty}\!t^{n-1}e^{-{\scriptstyle{Re}}\lambda
t}\:dt=\\
&=&\frac{M}{(n-1)!}\frac{n-1}{\mbox{\em
Re}\lambda}\int\limits_{0}^{\infty}\!t^{n-2}e^{-{\scriptstyle{Re}}\lambda
t}\:dt=\cdots=\frac{M}{(\mbox{\em Re}\lambda)^n}
\end{eqnarray*}
pour tout $n\in{\bf N}^*$. Avec le th\'eor\`eme de Hille-Yosida, on
voit que l'op\'erateur $A$ est le g\'en\'erateur infinit\'esimal d'un
semi-groupe $\left\{T(t)\right\}_{t\geq 0}\in{\cal SG}(M,0)$.\\
Soit $x\in{\cal D}(A^2)$ et $\lambda\in\Lambda_0$. Compte tenu du
th\'eor\`eme \ref{num47}, nous avons:
$$
T(t)x=\frac{1}{2\pi
i}\int\limits_{{\scriptstyle{Re}}\lambda-i\infty}^{{\scriptstyle{Re}}\lambda+i\infty}\!\!e^{\lambda
t}R(\lambda;A)x\:d\lambda
$$
et compte tenu du (ii) et du th\'eor\`eme de Cauchy, on peut remplacer le contour d'int\'egration par $\Gamma_\nu$. Donc:
$$
T(t)x=\frac{1}{2\pi
i}\int\limits_{\Gamma_\nu}\!\!e^{\lambda
t}R(\lambda;A)x\:d\lambda\quad,\quad(\forall)x\in{\cal D}(A^2).
$$
Comme $\overline{{\cal D}(A^2)}={\cal E}$ et l'int\'egrale
$\int\limits_{\Gamma_\nu}e^{\lambda t}R(\lambda;A)d\lambda$ est
uniform\'ement convergente, nous obtenons:
$$
T(t)x=\frac{1}{2\pi
i}\int\limits_{\Gamma_\nu}\!\!e^{\lambda
t}R(\lambda;A)x\:d\lambda\quad,\quad(\forall)x\in{\cal E},
$$
d'o\`u il r\'esulte l'affirmation de l'\'enonc\'e.\fin

\begin{teo}\label{num51}
Soit $A$ un op\'erateur lin\'eaire ferm\'e d\'efini sur un sous
espace dense de ${\cal E}$ et v\'erifiant les propri\'et\'es
suivantes:\\
i) $\sigma(A)\subset\{z\in{\bf C}|\mbox{Re}z\leq 0\}$;\\
ii) pour tout $x\in{\cal E}$ et tout $x^*\in{\cal E}^*$ on a:
$$
\sup_{\gamma>0}\:\gamma\!\!\!\int\limits_{\gamma-i\infty}^{\gamma+i\infty}\left|\left\langle
{R(\lambda;A)}^2x,x^*\right\rangle\right||d\lambda|<\infty\quad.
$$
Alors $A$ est le g\'en\'erateur infinit\'esimal d'un semi-groupe
$\left\{T(t)\right\}_{t\geq 0}\in{\cal SG}(M,0)$.
\end{teo}
\dem
Pour $\gamma>0$ arbitrairement fix\'e, l'application:
$$
\{z\in{\bf C}|\mbox{\em Re}z>0\}\ni z\longmapsto\left\langle
R(z+\gamma;A)^2x,x^*\right\rangle\in{\bf C}
$$
se trouve dans l'espace de Hardy $H^1$. Par cons\'equent elle admet un
repr\'esentation par une int\'egrale de Cauchy et en particulier
pour $\alpha>0$ et $\gamma\in]0,\alpha[$, on a:
$$
\left\langle R(\alpha;A)^2x,x^*\right\rangle=\frac{1}{2\pi
i}\int\limits_{\gamma-i\infty}^{\gamma+i\infty}\!\!\frac{\left\langle
R(\lambda;A)^2x,x^*\right\rangle}{\alpha-\lambda}\:d\lambda\quad.
$$
De plus, par r\'ecurrence on voit que:
$$
\left\langle \frac{d^n}{d\alpha^n}R(\alpha;A)^2x,x^*\right\rangle=(-1)^n\frac{n!}{2\pi
i}\int\limits_{\gamma-i\infty}^{\gamma+i\infty}\!\!\frac{\left\langle
R(\lambda;A)^2x,x^*\right\rangle}{(\alpha-\lambda)^{n+1}}\:d\lambda\quad,
$$
pour tout $n\in{\bf N}^*$. D'autre part, avec la proposition
\ref{num14} (iii) il vient:
$$
\frac{d^n}{d\alpha^n}R(\alpha;A)=(-1)^nn!{R(\alpha;A)}^{n+1}\quad,\quad(\forall)n\in{\bf
N}^*
$$
Par it\'eration nous obtenons:
$$
\left\langle R(\alpha;A)^{n+1}x,x^*\right\rangle=\frac{1}{2n\pi i}\int\limits_{\gamma-i\infty}^{\gamma+i\infty}\!\!\frac{\left\langle
R(\lambda;A)^2x,x^*\right\rangle}{(\alpha-\lambda)^n}\:d\lambda\quad,
$$
pour tout $n\in{\bf N}^*$ et tout $\gamma\in]0,\alpha[$. Il s'ensuit que
pour tout $n\in{\bf N}^*$ on a:
\begin{eqnarray*}
& &\left|\left\langle R(\alpha;A)^{n+1}x,x^*\right\rangle\right|\leq\\
&\leq&\frac{1}{2n\pi}\frac{1}{(\alpha-\gamma)^n}\int\limits_{\gamma-i\infty}^{\gamma+i\infty}\!\!\left|\left\langle
R(\lambda;A)^2x,x^*\right\rangle\right|\:|d\lambda|\leq\\
&\leq&\frac{1}{2n\pi}\frac{C}{\gamma}\frac{1}{(\alpha-\gamma)^n}\quad,\quad(\forall)\gamma\in]0,\alpha[,
\end{eqnarray*}
o\`u la constante $C$ ne d\'epend que de $x\in{\cal E}$ et de
$x^*\in{\cal E}^*$. Si nous prenons:
$$
\gamma=\frac{\alpha}{n+1}\quad,
$$
alors on voit que:
\begin{eqnarray*}
& &\left|\left\langle
R(\alpha;A)^{n+1}x,x^*\right\rangle\right|\leq\frac{1}{2n\pi}\frac{C}{\frac{\alpha}{n+1}}\frac{1}{\left(\alpha-\frac{\alpha}{n+1}\right)^n}=\\
&=&\frac{C}{2\pi}\frac{1}{\alpha^{n+1}}\frac{1}{\left(\frac{n}{n+1}\right)^{n+1}}=\frac{C}{\pi}\frac{1}{\alpha^{n+1}}\frac{1}{2}\left(1-\frac{1}{n+1}\right)^{-(n+1)}\leq\\
&\leq&\frac{Ce}{\pi}\frac{1}{\alpha^{n+1}}\quad.
\end{eqnarray*}
En appliquant le th\'eor\`eme de Banach-Steinhaus (\cite[Theorem
II.1.11, pag. 52]{dunford-schwartz}), on obtient pour tout $\alpha>0$:
$$
\left\|R(\alpha;A)^m\right\|\leq\frac{M}{\alpha^m}\quad,\quad(\forall)m\in\{2,3,\ldots\}.
$$
Prouvons que cette in\'egalit\'e reste valable pour $m=1$. On a:
\begin{eqnarray*}
& &\int\limits_\alpha^\beta\!R(\tau;A)^2\:d\tau=-\int\limits_\beta^\alpha\!R(\tau;A)^2\:d\tau=\\
&=&\int\limits_\beta^\alpha\!\frac{d}{d\tau}R(\tau;A)\:d\tau=R(\alpha;A)-R(\beta;A)\quad,\quad(\forall)\alpha,\beta>0.
\end{eqnarray*}
Comme pour $m=2$ nous avons:
$$
\left\|R(\beta;A)^2\right\|\leq\frac{M}{\beta^2}\quad,
$$
on en d\'eduit que la limite existe et on pose
$$
R_0x:=\lim_{\beta\rightarrow\infty}R(\beta;A)x\quad,\quad(\forall)x\in{\cal
E}.
$$
D'autre part, si $x\in{\cal D}(A)$, alors nous avons:
$$
R(\beta;A)x=R(\beta;A)^2(\beta I-A)x\longrightarrow 0\quad\mbox{si}\quad\beta\rightarrow\infty.
$$
Comme $\overline{{\cal D}(A)}={\cal E}$, il s'ensuit que $R_0=0$ et on
voit que:
$$
R(\alpha;A)=\int\limits_\alpha^\infty\!R(\tau;A)^2\:d\tau\quad.
$$
Par cons\'equent:
$$
\left\|R(\alpha;A)^m\right\|\leq\frac{M}{\alpha^m}\quad,\quad(\forall)m\in{\bf
N}^*.
$$
En appliquant le th\'eor\`eme de Hille-Yosida, on obtient le r\'esultat d\'esir\'e.\fin
\vspace{2cm}

    \section{Propri\'et\'es spectrales des $C_0$-semi-groupes}

\hspace{1cm}Nous terminons ce chapitre avec quelques propri\'et\'es spectrales pour
les $C_0$-semi-groupes.

\begin{lema}\label{num23}
Soient \semi et $A$ son g\'en\'erateur infinit\'esimal. Alors pour tout
$\lambda\in\Lambda_\omega$ et $t>0$, l'application:
$$
B_\lambda(t):{\cal E}\longrightarrow{\cal E}
$$
$$
B_\lambda(t)x=\int\limits_{0}^{t}\!\!e^{\lambda(t-s)}T(s)x\:ds
$$
d\'efinit un op\'erateur lin\'eaire born\'e sur ${\cal E}$ v\'erifiant les
propri\'et\'es suivantes:
$$
(\lambda I-A)B_\lambda(t)x=e^{\lambda
t}x-T(t)x\quad,\quad(\forall)x\in{\cal E}
$$
et:
$$
B_\lambda(t)(\lambda I-A)x=e^{\lambda
t}x-T(t)x\quad,\quad(\forall)x\in{\cal D}(A).
$$
De plus $B_\lambda(t)T(t)=T(t)B_\lambda(t)$.
\end{lema}
\dem
Pour tout $x\in{\cal E}$ nous avons successivement:
\begin{eqnarray*}
\left\|B_\lambda(t)x\right\|&=&\left\|\int\limits_{0}^{t}\!\!e^{\lambda(t-s)}T(s)x\:ds\right\|\leq\int\limits_{0}^{t}\!\!e^{{\scriptstyle
Re}\lambda(t-s)}\|T(s)\|\|x\|\:ds\leq\\
&\leq&Me^{{\scriptstyle Re}\lambda
t}\|x\|\int\limits_{0}^{t}\!\!e^{-({\scriptstyle
Re}\lambda-\omega)s}\:ds<\infty\quad.
\end{eqnarray*}
Comme la lin\'earit\'e est \'evidente, il en r\'esulte que
$B_\lambda(t)\in{\cal B(E)}$, quels que soient $\lambda\in\Lambda_\omega$ et
$t>0$.\\
Si $x\in{\cal E}$ et $h>0$, alors nous obtenons:
\begin{eqnarray*}
&
&\frac{T(h)-I}{h}B_\lambda(t)x=\frac{T(h)-I}{h}\int\limits_{0}^{t}\!e^{\lambda(t-s)}T(s)x\:ds=\\
&=&\frac{1}{h}\int\limits_{0}^{t}\!e^{\lambda(t-s)}T(h+s)x\:ds-\frac{1}{h}\int\limits_{0}^{t}\!e^{\lambda(t-s)}T(s)x\:ds=\\
&=&\frac{1}{h}\int\limits_{h}^{t+h}\!e^{\lambda(t-\tau+h)}T(\tau)x\:d\tau-\frac{1}{h}\int\limits_{0}^{t}\!e^{\lambda(t-s)}T(s)x\:ds=\\
&=&\frac{e^{\lambda
h}}{h}\int\limits_{h}^{t+h}e^{\lambda(t-\tau)}T(\tau)x\:d\tau-\frac{1}{h}\int\limits_{0}^{t}\!e^{\lambda(t-s)}T(s)x\:ds=\\
&=&\frac{e^{\lambda
h}}{h}\left[\int\limits_{0}^{t+h}\!e^{\lambda(t-\tau)}T(\tau)x\:d\tau-\int\limits_{0}^{h}\!e^{\lambda(t-\tau)}T(\tau)x\:d\tau\right]-\\
&-&\frac{1}{h}\int\limits_{0}^{t}\!e^{\lambda(t-s)}T(s)x\:ds=\\
&=&\frac{e^{\lambda
h}}{h}\int\limits_{0}^{t+h}\!e^{\lambda(t-\tau)}T(\tau)x\:d\tau-\frac{e^{\lambda
h}
}{h}\int\limits_{0}^{t}\!e^{\lambda(t-\tau)}T(\tau)x\:d\tau+\\
&+&\frac{e^{\lambda
h}}
{h}\int\limits_{0}^{t}\!e^{\lambda(t-\tau)}T(\tau)x\:d\tau-\frac{1}{h}\int\limits_{0}^{t}\!e^{\lambda(t-s)}T(s)x\:ds-\\
&-&
\frac{e^{\lambda h}}{h}\int\limits_{0}^{h}\!e^{\lambda(t-\tau)}T(\tau)x\:d\tau=\\
&=&\frac{e^{\lambda
h}}{h}\int\limits_{t}^{t+h}\!e^{\lambda(t-\tau)}T(\tau)x\:d\tau+\frac{e^{\lambda
h}-1}{h}\int\limits_0^t\!e^{\lambda(t-s)}T(s)x\:ds-\\
&-&\frac{e^{\lambda
h}}{h}\int\limits_{0}^{h}\!e^{\lambda(t-\tau)}T(\tau)x\:d\tau=\\
&=&\frac{e^{\lambda
h}}{h}\int\limits_{t}^{t+h}\!e^{\lambda(t-\tau)}T(\tau)x\:d\tau+\frac{e^{\lambda
h}-1}{h}B_\lambda(t)x-\frac{e^{\lambda
h}}{h}\int\limits_{0}^{h}\!e^{\lambda(t-\tau)}T(\tau)x\:d\tau\quad.
\end{eqnarray*}
En passant \`a limite, on a:
$$
\lim_{h\searrow
0}\frac{T(h)B_\lambda(t)x-B_\lambda(t)x}{h}=T(t)x+\lambda
B_\lambda(t)x-e^{\lambda t}x\quad,
$$
d'o\`u $B_\lambda(t)x\in{\cal D}(A)$ et:
$$
(\lambda I-A)B_\lambda(t)x=e^{\lambda
t}x-T(t)x\quad,\quad(\forall)x\in{\cal E}.
$$
Si $x\in{\cal D}(A)$, alors nous avons:
\begin{eqnarray*}
& &B_\lambda(t)Ax=\int\limits_{0}^{t}\!\!e^{\lambda(t-s)}T(s)Ax\:ds=\\
&=&\int\limits_{0}^{t}\!\!e^{\lambda(t-s)}\frac{d}{ds}T(s)x\:ds=T(t)x-e^{\lambda t}x+\lambda B_\lambda(t)x\quad,
\end{eqnarray*}
d'o\`u l'on tire:
$$
B_\lambda(t)(\lambda I-A)x=e^{\lambda
t}x-T(t)x\quad,\quad(\forall)x\in{\cal D}(A).
$$
De plus, nous obtenons que:
$$
B_\lambda(t){\cal D}(A)\subseteq{\cal D}(A)
$$
et:
$$
(\lambda I-A)B_\lambda(t)x=B_\lambda(t)(\lambda
I-A)x\quad,\quad(\forall)x\in{\cal D}(A)\quad,
$$
d'o\`u:
$$
AB_\lambda(t)x=B_\lambda(t)Ax\quad,\quad(\forall)x\in{\cal D}(A).
$$
Compte tenu du th\'eor\`eme \ref{num57}, on voit que:
$$
B_\lambda(t)T(t)=T(t)B_\lambda(t)\quad,\quad(\forall)t\geq 0.\fin
$$

\begin{teo}[spectral mapping]\label{num24}
Soient \semi et A son g\'en\'erateur infinit\'esimal. Alors:
$$
e^{t\sigma(A)}=\left\{\left.e^{\lambda t}\right|\lambda\in\sigma(A)\right\}\subseteq\sigma(T(t))\quad,\quad(\forall)t\geq 0.
$$
\end{teo}
\dem
Soit $\lambda\in{\bf C}$ tel que $e^{\lambda t}\in\rho(T(t))$.\\
Alors on peut consid\'erer l'op\'erateur $Q=\left(e^{\lambda t}I-T(t)\right)^{-1}\in{\cal
B(E)}$. Compte tenu du lemme \ref{num23} (i), on a:
$$
(\lambda I-A)B_\lambda(t)x=e^{\lambda
t}x-T(t)x\quad,\quad(\forall)x\in{\cal E}
$$
et:
$$
B_\lambda(t)(\lambda I-A)x=e^{\lambda
t}x-T(t)x\quad,\quad(\forall)x\in{\cal D}(A).
$$
Par multiplication avec $Q$ \`a droite dans la premi\`ere \'egalit\'e et
\`a gauche dans la seconde, nous obtenons:
$$
(\lambda I-A)B_\lambda(t)Qx=x\quad,\quad(\forall)x\in{\cal E}
$$
et:
$$
QB_\lambda(t)(\lambda I-A)x=x\quad,\quad(\forall)x\in{\cal D}(A).
$$
Mais, avec le lemme \ref{num23}, il en r\'esulte que:
$$
\left(e^{\lambda t}I-T(t)\right)B_\lambda(t)=B_\lambda(t)\left(e^{\lambda
t}I-T(t)\right)\quad,
$$
et nous voyons que $QB_\lambda(t)=B_\lambda(t)Q$. Par cons\'equent:
$$
(\lambda I-A)B_\lambda(t)Qx=x\quad,\quad(\forall)x\in{\cal E}
$$
et:
$$
B_\lambda(t)Q(\lambda I-A)x=x\quad,\quad(\forall)x\in{\cal D}(A).
$$
Il s'ensuit que $\lambda\in\rho(A)$ et finalement on voit que:
$$
\rho(T(t))\subset e^{t\rho(A)}\quad,\quad(\forall)t\geq 0,
$$
ou bien:
$$
e^{t\sigma(A)}\subseteq\sigma(T(t))\quad,\quad(\forall)t\geq 0.\fin
$$

\begin{obs}
\em
Nous avons vu que pour les semi-groupes uniform\'ement continus on a
l'\'egalit\'e:
$$
e^{t\sigma(A)}=\sigma(T(t))\quad,\quad(\forall)t\geq 0.
$$
Mais il existe des $C_0$-semi-groupes pour lesquels l'inclusion du
th\'eor\`eme \ref{num24} est stricte.
\end{obs}

\begin{definitie}
On dit que le $C_0$-semi-groupe $\left\{T(t)\right\}_{t\geq 0}$ est
nilpotent s'il existe $t_0>0$ tel que $T(t)=0$, pour tout $t>t_0$.
\end{definitie}

\begin{prop}
Soient \semi un semi-groupe nilpotent et
$A$ son g\'en\'erateur infinit\'esimal. Alors $\sigma(A)=\emptyset$.
\end{prop}
\dem
Comme le $C_0$-semi-groupe $\left\{T(t)\right\}_{t\geq 0}$ est nilpotent, il existe $t_0>0$ tel que $T(t)=0$, $(\forall)t>t_0$. 
Pour tout $\lambda\in{\bf C}$ et tout $x\in{\cal E}$, on a:
$$
\left\|e^{-\lambda t}T(t)x\right\|\leq e^{-{\scriptstyle Re}\lambda
t}Me^{\omega t}\|x\|\quad,\quad(\forall)t\in[0,t_0]
$$
et comme:
$$
\int\limits_{t_0}^{\infty}\!e^{-\lambda t}T(t)x\:dt=0\quad,
$$
on peut d\'efinir la transform\'ee de Laplace:
$$
R_\lambda:{\cal E}\longrightarrow{\cal E}
$$
$$
R_\lambda x=\int\limits_{0}^{\infty}\!\!e^{-\lambda t}T(t)x\:dt=\int\limits_{0}^{t_0}\!e^{-\lambda t}T(t)x\:dt
$$
pour tout $\lambda\in{\bf C}$. Avec le th\'eor\`eme \ref{num13},
il vient $\lambda\in\rho(A)$ et $R_\lambda x=R(\lambda;A)x$, pour tout $x\in{\cal E}$. Donc $\rho(A)={\bf C}$, c'est-\`a-dire 
$\sigma(A)=\emptyset$.\fin

\begin{obs}
\em
Pour un semi-groupe nilpotent
\semi ayant pour g\'en\'erateur
infinit\'esimal l'op\'erateur $A$, l'inclusion du th\'eor\`eme \ref{num24} est
stricte.
\end{obs}

\begin{teo}\label{num56}
Soient \semi et $A$ son g\'en\'erateur infinit\'esimal. Alors:
$$
\sigma(R(\lambda;A))=\left\{\left.\frac{1}{\lambda-\zeta}\right|\zeta\in\sigma(A)\right\}\cup\{0\}
$$
quel que soit $\lambda\in\Lambda_\omega$.
\end{teo}
\dem
Soient $\lambda\in\Lambda_\omega$ et
$\mu\in\rho(A)$, $\mu\neq\lambda$. D\'efinissons:
$$
S:{\cal E}\longrightarrow{\cal E}
$$
par:
$$
S=(\lambda-\mu)(\lambda I-A)R(\mu;A).
$$
Comme $S$ est un op\'erateur ferm\'e, avec le th\'eor\`eme du graphe
ferm\'e, on voit que $S\in{\cal B(E)}$. De plus, pour tout $x\in{\cal E}$
nous avons:
\begin{eqnarray*}
& &SR(\lambda;A)x=(\lambda-\mu)(\lambda I-A)R(\mu;A)R(\lambda;A)x=\\
&=&(\lambda-\mu)(\lambda
I-A)R(\lambda;A)R(\mu;A)x=(\lambda-\mu)R(\mu;A)x
\end{eqnarray*}
et:
\begin{eqnarray*}
& &R(\lambda;A)Sx=R(\lambda;A)(\lambda-\mu)(\lambda I-A)R(\mu;A)x=\\
&=&(\lambda-\mu)R(\lambda;A)(\lambda
I-A)R(\mu;A)x=(\lambda-\mu)R(\mu;A)x\quad.
\end{eqnarray*}
Par cons\'equent $SR(\lambda;A)=R(\lambda;A)S$.\\
De m\^eme, pour $x\in{\cal E}$ on a:
\begin{eqnarray*}
& &S\left[\frac{1}{\lambda-\mu}I-R(\lambda;A)\right]x=\\
&=&(\lambda-\mu)(\lambda
I-A)R(\mu;A)\left[\frac{1}{\lambda-\mu}I-R(\lambda;A)\right]x=\\
&=&\left[(\lambda I-A)R(\mu;A)-(\lambda-\mu)R(\mu;A)\right]x=\\
&=&(\lambda I-A-\lambda I+\mu
I)R(\mu;A)x=\\
&=&(\mu I-A)R(\mu;A)x=x\quad.
\end{eqnarray*}
De fa\c con analogue, pour tout $x\in{\cal E}$ on peut montrer que:
$$
\left[\frac{1}{\lambda-\mu}I-R(\lambda;A)\right]Sx=x\quad.
$$
Par cons\'equent:
$$
\frac{1}{\lambda-\mu}\in\rho(R(\lambda;A))\quad,
$$
d'o\`u:
$$
\left\{\left.\frac{1}{\lambda-\mu}\right|\mu\in\rho(A)\right\}\subset\rho(R(\lambda;A))\quad.
$$
Il s'ensuit que:
$$
\sigma(R(\lambda;A))\subset\left\{\left.\frac{1}{\lambda-\zeta}\right|\zeta\in\sigma(A)\right\}\quad.
$$
R\'eciproquement, soit $\lambda\in\Lambda_\omega$ et $\mu\in{\bf C}$, $\mu\neq\lambda$, tel que
$\frac{1}{\lambda-\mu}\in\rho(R(\lambda;A))$. Alors il existe
$R\left(\frac{1}{\lambda-\mu};R(\lambda;A)\right)\in{\cal B(E)}$ et
pour tout $x\in{\cal D}(A)$ nous avons:
\begin{eqnarray*}
& &R(\lambda;A)R\left(\frac{1}{\lambda-\mu};R(\lambda;A)\right)x=R(\lambda;A)\left[\frac{1}{\lambda-\mu}I-R(\lambda;A)\right]^{-1}x=\\
&=&\left[R(\lambda;A)^{-1}\right]^{-1}\left[\frac{1}{\lambda-\mu}I-R(\lambda;A)\right]^{-1}x=\\
&=&\left\{\left[\frac{1}{\lambda-\mu}I-R(\lambda;A)\right]{R(\lambda;A)}^{-1}\right\}^{-1}x=\\
&=&\left[\frac{1}{\lambda-\mu}{R(\lambda;A)}^{-1}-I\right]^{-1}x=\left\{{R(\lambda;A)}^{-1}\left[\frac{1}{\lambda-\mu}I-R(\lambda;A)\right]\right\}^{-1}x=\\
&=&\left[\frac{1}{\lambda-\mu}I-R(\lambda;A)\right]^{-1}R(\lambda;A)x=R\left(\frac{1}{\lambda-\mu};R(\lambda;A)\right)R(\lambda;A)x\quad.
\end{eqnarray*}
Posons:
$$
Q=R(\lambda;A)R\left(\frac{1}{\lambda-\mu};R(\lambda;A)\right).
$$
Pour tout $x\in{\cal D}(A)$, nous avons:
\begin{eqnarray*}
& &(\mu I-A)Qx=(\mu I-\lambda I+\lambda
I-A)R(\lambda;A)R\left(\frac{1}{\lambda-\mu};R(\lambda;A)\right)x=\\
&=&\left[(\lambda
I-A)R(\lambda;A)-(\lambda-\mu)R(\lambda;A)\right]R\left(\frac{1}{\lambda-\mu};R(\lambda;A)\right)x=\\
&=&\left[I-(\lambda-\mu)R(\lambda;A)\right]R\left(\frac{1}{\lambda-\mu};R(\lambda;A)\right)x=\\
&=&(\lambda-\mu)\left[\frac{1}{\lambda-\mu}I-R(\lambda;A)\right]R\left(\frac{1}{\lambda-\mu};R(\lambda;A)\right)x=(\lambda-\mu)x\quad,
\end{eqnarray*}
d'o\`u il r\'esulte que:
$$
\frac{1}{\lambda-\mu}(\mu I-A)Qx=x\quad,\quad(\forall)x\in{\cal D}(A).
$$
De m\^eme, nous obtenons:
\begin{eqnarray*}
& &Q(\mu
I-A)x=R(\lambda;A)R\left(\frac{1}{\lambda-\mu};R(\lambda;A)\right)(\mu
I-A)x=\\
&=&R\left(\frac{1}{\lambda-\mu};R(\lambda;A)\right)R(\lambda;A)(\mu
I-\lambda I+\lambda I-A)x=\\
&=&R\left(\frac{1}{\lambda-\mu};R(\lambda;A)\right)\left[R(\lambda;A)(\lambda
I-A)-R(\lambda;A)(\lambda-\mu)\right]x=\\
&=&R\left(\frac{1}{\lambda-\mu};R(\lambda;A)\right)\left[I-(\lambda-\mu)R(\lambda;A)\right]x=\\
&=&(\lambda-\mu)R\left(\frac{1}{\lambda-\mu};R(\lambda;A)\right)\left[\frac{1}{\lambda-\mu}I-R(\lambda;A)\right]x=(\lambda-\mu)x\quad,
\end{eqnarray*}
d'o\`u:
$$
\frac{1}{\lambda-\mu}Q(\mu I-A)x=x\quad,\quad(\forall)x\in{\cal D}(A).
$$
Par cons\'equent $\mu\in\rho(A)$. Il s'ensuit que:
$$
\rho(R(\lambda;A))\subset\left\{\left.\frac{1}{\lambda-\mu}\right|\mu\in\rho(A)\right\}\quad,
$$
ou bien:
$$
\left\{\left.\frac{1}{\lambda-\zeta}\right|\zeta\in\sigma(A)\right\}\subset\sigma(R(\lambda;A))\quad,\quad(\forall)\lambda\in{\bf
C}\mbox{ avec }\mbox{\em Re}\lambda>\omega.
$$
Finalement, nous voyons que:
$$
\sigma(R(\lambda;A))=\left\{\left.\frac{1}{\lambda-\zeta}\right|\zeta\in\sigma(A)\right\}\quad,\quad(\forall)\lambda\in\Lambda_\omega.
$$
Si $0\in\rho (R(\lambda;A))$, alors il existe
$\left(0I-R(\lambda;A)\right)^{-1}\in{\cal B(E)}$, d'o\`u $A\in{\cal
B(E)}$ ce qui est absurde. Par
cons\'equent $0\in\sigma(R(\lambda;A))$.\fin
\vspace{2cm}

	\section{Notes}

{\footnotesize
Les notions et les r\'esultats de ce chapitre se trouvent dans les
monographies concernant les $C_0$-semi-groupes d'op\'erateurs
lin\'eaires born\'es. Le th\'eor\`eme \ref{num50} se trouve dans
\cite[pag.184]{hille}, mais une preuve \'el\'egante utilisant le
th\'eor\`eme de Krein-\v Smulian se trouve dans
\cite[pag. 15]{davies}. De m\^eme, dans \cite[pag. 43]{pazy1} on peut
trouver une caract\'erisation du g\'en\'erateur infinit\'esimal d'un
$C_0$-semi-groupe pour la topologie faible.

Pour le th\'eor\`eme de l'unicit\'e de l'engendrement nous avons
utilis\'e \cite[pag. 6]{pazy1} et le th\'eor\`eme \ref{num57} se trouve
dans \cite[pag. 11]{davies}.

Le r\'esultat le plus important de ce chapitre est le th\'eor\`eme de
Hille-Yosida. Il a \'et\'e montr\'e pour la premi\`ere fois
ind\'ependamment par Hille dans \cite{hille} et par Yosida dans
\cite{yosida2} pour les $C_0$-semi-groupes de contractions. Quelques
ann\'ees plus tard, Feller dans \cite{feller}, Miyadera dans
\cite{miyadera} et Phillips dans \cite{phillips} donnent une preuve
pour le cas g\'en\'eral d'un $C_0$-semi-groupe. Nous avons utilis\'e
les id\'ees du livre de Pazy \cite[pag. 8]{pazy1} pour obtenir une preuve
dans le cas le plus g\'en\'eral, en utilisant l'approximation
g\'en\'eralis\'ee de Yosida que nous avons introduit dans la
d\'efinition \ref{num55}.

Pour obtenir le repr\'esentation de Bromwich d'un $C_0$-semi-groupe,
nous avons utilis\'e aussi les id\'ees de Pazy de
\cite[pag. 29]{pazy1}. Une variante du lemme \ref{num36} se trouve
dans \cite{lemle2}.

Pour le th\'eor\`eme \ref{num45} on peut consulter
\cite[pag. 30]{pazy1} ou bien \cite[pag. 76]{ahmed}. Le t\'eor\`eme
\ref{num51} a \'et\'e montr\'e par Gomilko dans \cite{gomilko}.

Pour les propri\'et\'es spectrales des $C_0$-semi-groupes on peut
consulter \cite[pag. 44]{pazy1} o\`u on peut trouver aussi des autres
r\'esultats sur cette probl\`eme. Finalement, pour le th\'eor\`eme
\ref{num56} on pourra consulter \cite[pag. 39]{davies}.
}

  \chapter{$C_0$-semigroupes avec propri\'et\'es sp\'eciales}

    \section{$C_0$-semi-groupes diff\'erentiables}

\hspace{1cm}Par la suite, nous \' etudierons les propri\'et\'es des $C_0$-semi-groupes
pour lesquels l'application $]0,\infty)\ni t\longmapsto T(t)x\in{\cal
E}$ est diff\'erentiable, quel que soit $x\in{\cal E}$.

\begin{definitie}
On dit que $\left\{T(t)\right\}_{t\geq 0}$ est un $C_0$-semi-groupe  diff\'erentiable (et notons
$\left\{T(t)\right\}_{t\geq 0}\in{\cal SGD}(M,\omega)$) si
l'application:
$$
]0,\infty)\ni t\longmapsto T(t)x\in{\cal E}
$$
est diff\'erentiable, quel que soit $x\in{\cal E}$.
\end{definitie}

\begin{teo}\label{num27}
Soient \semi et $A$ son g\'en\'erateur infinit\'esimal. Les affirmations
suivantes sont \'equivalentes:\\
i) \dif;\\
ii) ${\cal I}\mbox{m }T(t)\subset{{\cal D}(A)}\quad,\quad(\forall)t>0$.
\end{teo}
\dem
$i)\Longrightarrow ii)$ Soient $x\in{\cal E}$ et $t,h>0$. Puisque
l'application:
$$
]0,\infty)\ni t\longmapsto T(t)x\in{\cal
E}
$$
est diff\'erentiable, la limite du rapport
$$
\frac{T(t+h)x-T(t)x}{h}\quad,
$$
lorsque $h\searrow 0$,
existe et est \'egale par d\'efinition avec $AT(t)x$. Par cons\'equent, $T(t)x\in{\cal D}(A)$.\\
$ii)\Longrightarrow i)$ Consid\'erons  $x\in{\cal E}$ et $t,h>0$. Comme $T(t)x\in{\cal D}(A)$, nous avons:
$$
\frac{d^+T(t)x}{dt}=\lim_{h\searrow
0}\frac{T(t+h)x-T(t)x}{h}=AT(t)x\quad.
$$
D'autre part, pour $h\in]0,t[$ et $\delta\in]0,t-h[$ on a:
\begin{eqnarray*}
& &\left\|\frac{T(t-h)x-T(t)x}{-h}-AT(t)x\right\|=\\
&=&\left\|\frac{T(t-\delta)T(\delta)x-T(t-h-\delta)T(\delta)x}{h}-AT(\delta)T(t-\delta)x\right\|=\\
&=&\left\|\frac{1}{h}\left[\int\limits_{t-h-\delta}^{t-\delta}\!\frac{d}{d\tau}T(\tau)T(\delta)x\:d\tau-\int\limits_{t-h-\delta}^{t-\delta}\!AT(\delta)T(t-\delta)x\:d\tau\right]\right\|=\\
&=&\left\|\frac{1}{h}\int\limits_{t-h-\delta}^{t-\delta}\left[AT(\delta)T(\tau)-AT(\delta)T(t-\delta)\right]x\:d\tau\right\|\leq\\
&\leq&\frac{1}{h}\left\|AT(\delta)\right\|\int\limits_{t-h-\delta}^{t-\delta}\left\|T(\tau)-T(t-\delta)\right\|\:d\tau\|x\|=\\
&=&\frac{1}{h}\left\|AT(\delta)\right\|h\left\|T(c)-T(t-\delta)\right\|\|x\|=\\
&=&\left\|AT(\delta)\right\|\left\|T(c)-T(t-\delta)\right\|\|x\|\quad,
\end{eqnarray*}
o\`u $c\in[t-h-\delta,t-\delta]$. Par cons\'equent:
$$
\frac{d^-T(t)x}{dt}=\lim_{h\searrow
0}\frac{T(t-h)x-T(t)x}{-h}=AT(t)x\quad.
$$
Donc $\left\{T(t)\right\}_{t\geq 0}$ est un $C_0$-semi-groupe
diff\'erentiable.\fin

\begin{prop}\label{num26}
Soit \dif. Alors l'application:
$$
]0,\infty)\ni t\longmapsto T(t)\in{\cal B(E)}
$$
est continue pour la topologie de la convergence uniforme.
\end{prop}
\dem
Soient $x\in{\cal E}$ et $t_1,t_2\in]0,\infty)$ tel que
$t_1<t_2$. Compte tenu du th\'eor\`eme \ref{num27}, nous obtenons:
\begin{eqnarray*}
\left\|T(t_1)x-T(t_2)x\right\|&=&\left\|\int\limits_{t_1}^{t_2}\!\frac{d}{ds}T(s)x\:ds\right\|=\left\|\int\limits_{t_1}^{t_2}\!AT(t_1)T(s-t_1)x\:ds\right\|\leq\\
&\leq&\left\|AT(t_1)\right\|\int\limits_{t_1}^{t_2}\!Me^{(s-t_1)\omega}\|x\|\:ds\quad.
\end{eqnarray*}
Par suite, nous avons:
$$
\left\|T(t_1)-T(t_2)\right\|\leq\left\|AT(t_1)\right\|M\int\limits_{t_1}^{t_2}\!\!e^{(s-t_1)\omega}\:ds\quad,
$$
d'o\`u r\'esulte la continuit\'e uniforme de l'application
consid\'er\'ee dans l'\'enonc\'e.\fin

\begin{teo}\label{num29}
Soient \dif et $A$ son g\'en\'erateur infinit\'esi-\\mal. Alors:\\
i) pour tout $n\in{\bf N}^*$ et tout $x\in{\cal E}$, on a
$T(t)x\in{\cal D}(A^n)$ et:
$$
A^nT(t)x=\left[AT\left(\frac{t}{n}\right)\right]^nx\quad,\quad(\forall)t>0;
$$
ii) pour tout $n\in{\bf N}^*$ l'application:
$$
]0,\infty)\ni t\longmapsto T(t):{\cal E}\rightarrow{\cal D}(A^n)
$$
est $n$ fois diff\'erentiable pour la topologie de la convergence
uniforme et:
$$
{T(t)}^{(n)}=\frac{d^n}{dt^n}T(t)=A^nT(t)\in{\cal
B(E)}\quad,\quad(\forall)t>0;
$$
iii) pour tout $n\in{\bf N}^*$ l'application:
$$
]0,\infty)\ni t\longmapsto {T(t)}^{(n)}\in{\cal B(E)}
$$
est continue pour la topologie de la convergence uniforme.
\end{teo}
\dem
Prouvons les affirmations de l'\'enonc\'e par r\'ecurrence.\\
i) Avec le th\'eor\`eme \ref{num27}, on voit que pour tout $x\in{\cal
E}$ on a $T(t)x\in{\cal D}(A)$ et:
$$
AT(t)x=\left[AT\left(\frac{t}{1}\right)\right]^1x\quad,\quad(\forall)t>0.
$$
Supposons que pour tout $x\in{\cal E}$ on ait $T(t)x\in{\cal D}(A^k)$
et:
$$
A^kT(t)x=\left[AT\left(\frac{t}{k}\right)\right]^kx\quad,\quad(\forall)t>0.
$$
Soient $x\in{\cal E}$ et $\delta\in]0,t[$. On voit que $T(t-\delta)T(\delta)x\in{\cal D}(A)$ et:
$$
AT(t)x=AT(t-\delta)T(\delta)x=T(t-\delta)AT(\delta)x\in{\cal D}(A^k)\quad.
$$
Par cons\'equent $T(t)x\in{\cal D}(A^{k+1})$, $(\forall)t>0$. De plus:
\begin{eqnarray*}
A^{k+1}T(t)x&=&A\left[A^kT(t-\delta)T(\delta)\right]x=A\left[T(t-\delta)A^kT(\delta)\right]x=\\
&=&AT(t-\delta)\left[AT\left(\frac{\delta}{k}\right)\right]^kx\quad.
\end{eqnarray*}
Si $\delta=\frac{kt}{k+1}$, il vient:
$$
A^{k+1}T(t)x=\left[AT\left(\frac{t}{k+1}\right)\right]^{k+1}x\quad.
$$
Finalement, nous
obtenons (i).\\
ii) Pour $n=1$, compte tenu du th\'eor\`eme \ref{num27} et de la proposition \ref{num26}, il r\'esulte que l'application:
$$
]0,\infty)\ni t\longmapsto T(t):{\cal E}\rightarrow{\cal D}(A)
$$
est diff\'erentiable pour la topologie de la convergence uniforme et:
$$
{T(t)}'=AT(t)\quad,\quad(\forall)t>0.
$$
Comme $A$ est un op\'erateur ferm\'e et $T(t)\in{\cal B(E)}$, il
r\'esulte que $AT(t)$ est un op\'erateur ferm\'e d\'efini sur ${\cal
E}$. Avec le th\'eor\`eme du graphe ferm\'e (\cite[Theorem II.2.4,
pag. 57]{dunford-schwartz}), on voit que
$AT(t)\in{\cal B(E)}$, $(\forall)t>0$. Supposons que l'application:
$$
]0,\infty)\ni t\longmapsto T(t):{\cal E}\rightarrow{\cal D}(A^k)
$$
est $k$ fois diff\'erentiable pour la topologie de la convergence
uniforme et:
$$
{T(t)}^{(k)}=A^kT(t)\in{\cal B(E)}\quad,\quad(\forall)t>0.
$$
De plus, avec la preuve pr\'ec\'edente, on voit que $T(t)x\in{\cal D}\left(A^{k+1}\right)$, pour tout $t>0$.
Soient $x\in{\cal E}$, $\|x\|\leq 1$ et $t>0$. Si $h>0$ et
$\delta\in]0,t[$, on a:
\begin{eqnarray*}
& &
\left\|\frac{{T(t+h)}^{(k)}x-{T(t)}^{(k)}x}{h}-A^{k+1}T(t)x\right\|=\\
&=&\left\|\frac{A^kT(\delta)T(t+h-\delta)x-A^kT(\delta)T(t-\delta)x}{h}-A^{k+1}T(\delta)T(t-\delta)x\right\|=\\
&=&\left\|A^kT(\delta)\frac{1}{h}\left[T(t+h-\delta)-T(t-\delta)\right]x-A^{k+1}T(\delta)T(t-\delta)x\right\|=\\
&=&\left\|A^kT(\delta)\frac{1}{h}\int\limits_{t-\delta}^{t+h-\delta}\!\frac{d}{d\tau}T(\tau)x\:d\tau-A^{k+1}T(\delta)\frac{1}{h}\int\limits_{t-\delta}^{t+h-\delta}\!T(t-\delta)x\:d\tau\right\|=\\
&=&\left\|A^kT(\delta)\frac{1}{h}\int\limits_{t-\delta}^{t+h-\delta}\!AT(\tau)x\:d\tau-A^{k+1}T(\delta)\frac{1}{h}\int\limits_{t-\delta}^{t+h-\delta}\!T(t-\delta)x\:d\tau\right\|=\\
&=&\left\|\frac{1}{h}A^{k+1}T(\delta)\int\limits_{t-\delta}^{t+h-\delta}\left[T(\tau)-T(t-\delta)\right]x\:d\tau\right\|\leq\\
&\leq&\frac{\left\|A^{k+1}T(\delta)\right\|}{h}\int\limits_{t-\delta}^{t+h-\delta}\left\|T(\tau)-T(t-\delta)\right\|\|x\|\:d\tau=\\
&=&\left\|A^{k+1}T(\delta)\right\|\left\|T(c)-T(t-\delta)\right\|\|x\|\quad,
\end{eqnarray*}
o\`u $c\in[t-\delta,t+h-\delta]$. Il s'ensuit que:
$$
\left\|\frac{{T(t+h)}^{(k)}-{T(t)}^{(k)}}{h}-A^{k+1}T(t)\right\|\leq\left\|A^{k+1}T(\delta)\right\|\left\|T(c)-T(t-\delta)\right\|\quad,
$$
o\`u $c\in[t-\delta,t+h-\delta]$. Par cons\'equent:
$$
\lim_{h\searrow
0}\frac{{T(t+h)}^{(k)}-{T(t)}^{(k)}}{h}=A^{k+1}T(t)\quad,\quad(\forall)t>0.
$$
Si $h>0$ tel que $t-h>0$ et $\delta\in]0,t-h[$, alors nous avons:
\begin{eqnarray*}
&
&\left\|\frac{{T(t-h)}^{(k)}x-{T(t)}^{(k)}x}{-h}-A^{k+1}T(t)x\right\|=\\
&=&\left\|\frac{A^kT(\delta)T(t-\delta)x-A^kT(\delta)T(t-h-\delta)x}{h}-A^{k+1}T(\delta)T(t-\delta)x\right\|=\\
&=&\left\|A^kT(\delta)\frac{1}{h}\left[T(t-\delta)-T(t-h-\delta)\right]x-A^{k+1}T(\delta)T(t-\delta)x\right\|=\\
&=&\left\|A^kT(\delta)\frac{1}{h}\int\limits_{t-h-\delta}^{t-\delta}\frac{d}{d\tau}T(\tau)x\:d\tau-A^{k+1}T(\delta)\frac{1}{h}\int\limits_{t-h-\delta}^{t-\delta}T(t-\delta)x\:d\tau\right\|=\\
&=&\left\|A^kT(\delta)\frac{1}{h}\int\limits_{t-h-\delta}^{t-\delta}AT(\tau)x\:d\tau-A^{k+1}T(\delta)\frac{1}{h}\int\limits_{t-h-\delta}^{t-\delta}T(t-\delta)x\:d\tau\right\|=\\
&=&\left\|\frac{1}{h}A^{k+1}T(\delta)\int\limits_{t-h-\delta}^{t-\delta}\left[T(\tau)-T(t-\delta)\right]x\:d\tau\right\|\leq\\
&\leq&\frac{\left\|A^{k+1}T(\delta)\right\|}{h}\int\limits_{t-h-\delta}^{t-\delta}\left\|T(\tau)-T(t-\delta)\right\|\|x\|\:d\tau=\\
&=&\left\|A^{k+1}T(\delta)\right\|\:\left\|T(c)-T(t-\delta)\right\|\|x\|\quad,
\end{eqnarray*}
o\`u $c\in[t-h-\delta,t-\delta]$. Il vient:
$$
\left\|\frac{{T(t-h)}^{(k)}-{T(t)}^{(k)}}{-h}-A^{k+1}T(t)\right\|\leq\left\|A^{k+1}T(\delta)\right\|\left\|T(c)-T(t-\delta)\right\|\quad,
$$
o\`u $c\in[t-h-\delta,t-\delta]$. Par cons\'equent:
$$
\lim_{h\searrow
0}\frac{{T(t-h)}^{(k)}-{T(t)}^{(k)}}{-h}=A^{k+1}T(t)\quad,\quad(\forall)t>0.
$$
Il s'ensuit que ${T(t)}^{(k)}$ est diff\'erentiable pour la topologie de la
convergence uniforme et:
$$
\left({T(t)}^{(k)}\right)'={T(t)}^{(k+1)}=A^{k+1}T(t)\quad,\quad(\forall)t>0.
$$
Comme $A$ est un op\'erateur ferm\'e et $A^kT(t)\in{\cal B(E)}$, il
r\'esulte que $A\left(A^kT(t)\right)$ est un op\'erateur ferm\'e
d\'efini sur ${\cal E}$. Avec le th\'eor\`eme du graphe ferm\'e
(\cite[Theorem II.2.4, pag. 57]{dunford-schwartz}), on voit que ${T(t)}^{(k+1)}=A^{k+1}T(t)\in{\cal B(E)}$,
$(\forall)t>0$.\\
Finalement, on a obtenu (ii).\\
iii) Soient $x\in{\cal E}$ avec $\|x\|\leq 1$ et $t>0$. Pour $h>0$ et
$\delta\in]0,t[$ nous obtenons:
\begin{eqnarray*}
& &\left\|{T(t+h)}'x-{T(t)}'x\right\|=\left\|AT(t+h)x-AT(t)x\right\|\leq\\
&\leq&\left\|AT(\delta)\right\|\left\|T(t+h-\delta)-T(t-\delta)\right\|\|x\|\quad,
\end{eqnarray*}
d'o\`u il r\'esulte:
$$
\left\|{T(t+h)}'-{T(t)}'\right\|\leq\left\|AT(\delta)\right\|\left\|T(t+h-\delta)-T(t-\delta)\right\|\quad.
$$
De fa\c con analogue, pour $h>0$ et $\delta\in]0,t-h[$ nous obtenons:
\begin{eqnarray*}
& &\left\|{T(t-h)}'x-{T(t)}'x\right\|=\left\|AT(t-h)x-AT(t)x\right\|\leq\\
&\leq&\left\|AT(\delta)\right\|\left\|T(t-h-\delta)-T(t-\delta)\right\|\|x\|\quad,
\end{eqnarray*}
d'o\`u:
$$
\left\|{T(t-h)}'-{T(t)}'\right\|\leq\left\|AT(\delta)\right\|\left\|T(t-h-\delta)-T(t-\delta)\right\|\quad.
$$
Il est clair que l'application:
$$
]0,\infty)\ni t\longmapsto {T(t)}'\in{\cal B(E)}
$$
est continue pour la topologie de la convergence uniforme. Supposons
que l'application:
$$
]0,\infty)\ni t\longmapsto {T(t)}^{(k)}\in{\cal B(E)}
$$
est continue pour la topologie de la convergence uniforme. Si $h>0$
et $\delta\in]0,t[$, alors nous avons:
\begin{eqnarray*}
& &\left\|{T(t+h)}^{(k+1)}x-{T(t)}^{(k+1)}x\right\|=\left\|A^{k+1}T(t+h)x-A^{k+1}T(t)x\right\|\leq\\
&\leq&\left\|A^{k+1}T(\delta)\right\|\left\|T(t+h-\delta)-T(t-\delta)\right\|\|x\|\quad,
\end{eqnarray*}
d'o\`u il s'ensuit que:
$$
\left\|{T(t+h)}^{(k+1)}-{T(t)}^{(k+1)}\right\|\leq\left\|A^{k+1}T(\delta)\right\|\left\|T(t+h-\delta)-T(t-\delta)\right\|\quad.
$$
D'autre part, pour $h>0$ et $\delta\in]0,t-h[$ nous obtenons:
\begin{eqnarray*}
& &\left\|{T(t-h)}^{(k+1)}x-{T(t)}^{(k+1)}x\right\|=\left\|A^{k+1}T(t-h)x-A^{k+1}T(t)x\right\|\leq\\
&\leq&\left\|A^{k+1}T(\delta)\right\|\left\|T(t-h-\delta)-T(t-\delta)\right\|\|x\|
\end{eqnarray*}
et on voit que:
$$
\left\|{T(t-h)}^{(k+1)}-{T(t)}^{(k+1)}\right\|\leq\left\|A^{k+1}T(\delta)\right\|\left\|T(t-h-\delta)-T(t-\delta)\right\|\quad.
$$
Donc l'application:
$$
]0,\infty)\ni t\longmapsto {T(t)}^{(k+1)}\in{\cal B(E)}
$$
est continue pour la topologie de la convergence uniforme. La
propri\'et\'e (iii) en d\'ecoule imm\'ediatement.\fin

\begin{obs}
\em
Si \dif, alors l'application:
$$
]0,\infty)\ni t\longmapsto T(t)\in{\cal B(E)}
$$
est de classe ${\cal C}_{]0,\infty)}^\infty$.
\end{obs}

\begin{obs}\label{num46}
\em
Si \dif, alors pour tout $n\in{\bf N}^*$ on a:
$$
{T(t)}^{(n)}=A^nT(t)=\left[AT\left(\frac{t}{n}\right)\right]^n\quad,\quad(\forall)t>0.
$$
\end{obs}

Nous finissons cette section avec le th\'eor\`eme spectral pour les $C_0$-semi-groupes diff\'erentiables. Soit \semi. Pour tout $\lambda\in{\bf C}$ et tout $t>0$, nous avons d\'efini l'op\'erateur lin\'eaire born\'e:
$$
B_\lambda(t):{\cal E}\longrightarrow{\cal E}
$$
$$
B_\lambda(t)x=\int\limits_{0}^{t}\!\!e^{\lambda(t-s)}T(s)x\:ds
$$
et nous avons \'etudi\'e ses propri\'et\'es avec le lemme
\ref{num23}. Si le $C_0$-semi-groupe $\left\{T(t)\right\}_{t\geq 0}$
est diff\'erentiable, on peut montrer le r\'esultat suivant.

\begin{lema}\label{num30}
Soient \dif et $A$ son g\'en\'erateur infinit\'esimal. Alors:\\
i) pour tout $\lambda\in{\bf C}$ et tout $t>0$, l'op\'erateur
$B_\lambda(t)\in{\cal B(E)}$ est ind\'efiniment d\'erivable et:
$$
{B_\lambda(t)}^{(n)}=\lambda^n\left(B_\lambda(t)+\sum\limits_{i=0}^{n}\frac{{T(t)}^{(i)}}{\lambda^{i+1}}\right)\quad,\quad(\forall)n\in{\bf
N}^*;
$$
ii) pour tout $\lambda\in{\bf C}$ et tout $t>0$ on a:
$$
{B_\lambda(t)}^{(n)}{T(t)}^{(n)}={T(t)}^{(n)}{B_\lambda(t)}^{(n)}\quad,\quad(\forall)n\in{\bf
N}^*.
$$
\end{lema}
\dem
Montrons les affirmations de l'\'enonc\'e par r\'ecurrence.\\
i) Soient $x\in{\cal E}$, $\lambda\in{\bf C}$ et $t>0$. Alors:
$$
{B_\lambda(t)}'x=\lambda\left(B_\lambda(t)x+\frac{T(t)x}{\lambda}\right)\quad.
$$
Supposons que:
$$
{B_\lambda(t)}^{(k)}x=\lambda^k\left(B_\lambda(t)x+\sum\limits_{i=0}^{k}\frac{{T(t)}^{(i)}x}{\lambda^{i+1}}\right)\quad.
$$
Alors:
\begin{eqnarray*}
{B_\lambda(t)}^{(k+1)}x&=&{\left({B_\lambda(t)}^{(k)}x\right)}'=\lambda^k\left(\lambda
B_\lambda(t)x+T(t)x+\sum\limits_{i=0}^{k}\frac{{T(t)}^{(i+1)}x}{\lambda^{i+1}}\right)=\\
&=&\lambda^{k+1}\left(B_\lambda(t)x+\sum\limits_{i=0}^{k+1}\frac{{T(t)}^{(i)}x}{\lambda^{i+1}}\right)
\end{eqnarray*}
et nous obtenons (i).\\
ii) Soient  $\lambda\in{\bf C}$ et $t>0$. Compte tenu du th\'eor\`eme
\ref{num29}, pour $x\in{\cal E}$, on voit que $T(t)x\in{\cal D}(A^n)$ et:
\begin{eqnarray*}
& &A^nT(t)x=A^nT\left(\frac{nt}{n}\right)x=A^nT\left(\underbrace{\frac{t}{n}+\frac{t}{n}+\cdots+\frac{t}{n}}_{\mbox{\scriptsize{n fois}}}\right)x=\\
&=&A^n\underbrace{T\left(\frac{t}{n}\right)T\left(\frac{t}{n}\right)\cdots T\left(\frac{t}{n}\right)}_{\mbox{\scriptsize{n fois}}}x=\left[AT\left(\frac{t}{n}\right)\right]^nx=\\
&=&\left[T\left(\frac{t}{n}\right)A\right]^nx=T\left(\frac{nt}{n}\right)A^n=T(t)A^nx\quad,\quad(\forall
)n\in{\bf
N}^*,
\end{eqnarray*}
parce que le semi-groupe commute avec son g\'en\'erateur infinit\'esimal.
De m\^eme, avec le lemme \ref{num23} il r\'esulte que:
$$
B_\lambda(t)T(t)=T(t)B_\lambda(t)\quad.
$$
Alors pour $x\in{\cal E}$, nous avons:
\begin{eqnarray*}
& &{B_\lambda(t)}^{(n)}T(t)x=\lambda^n\left(B_\lambda(t)+\sum\limits_{i=0}^{n}\frac{{T(t)}^{(i)}}{\lambda^{i+1}}\right)T(t)x=\\
&=&\lambda^n\left(B_\lambda(t)T(t)+\sum\limits_{i=0}^{n}\frac{A^iT(t)T(t)}{\lambda^{i+1}}\right)x=\\
&=&\lambda^n\left(T(t)B_\lambda(t)+\sum\limits_{i=0}^{n}\frac{T(t)A^iT(t)}{\lambda_{i+1}}\right)x=\\
&=&T(t)\lambda^n\left(B_\lambda(t)+\sum\limits_{i=0}^{n}\frac{{T(t)}^{(i)}}{\lambda^{i+1}}\right)x=T(t){B_\lambda(t)}^{(n)}x\quad,\quad(\forall)n\in{\bf N}^*.
\end{eqnarray*}
D'autre part, pour $x\in{\cal D}(A)$, nous avons:
$$
B_\lambda(t)(\lambda I-A)x=(\lambda I-A)B_\lambda(t)x\quad,
$$
d'o\`u il r\'esulte:
$$
B_\lambda(t)Ax=AB_\lambda(t)x\quad.
$$
Supposons que pour $x\in{\cal D}(A^k)$ nous avons:
$$
B_\lambda(t)A^kx=A^kB_\lambda(t)x\quad.
$$
Si $x\in{\cal D}(A^{k+1})$, il vient:
$$
B_\lambda(t)A^{k+1}x=B_\lambda(t)A^k(Ax)=A^kB_\lambda(t)Ax=A^kAB_\lambda(t)x=A^{k+1}B_\lambda(t)x.
$$
Il s'ensuit donc que:
$$
B_\lambda(t)A^nx=A^nB_\lambda(t)x\quad,
$$
pour tout $x\in{\cal D}(A^n)$ et tout $n\in{\bf N}^*$.
De m\^eme, si $x\in{\cal D}(A^n)$, on a:
\begin{eqnarray*}
& &{B_\lambda(t)}^{(n)}A^nx=\lambda^n\left(B_\lambda(t)+\sum\limits_{i=0}^{n}\frac{{T(t)}^{(i)}}{\lambda^{i+1}}\right)A^nx=\\
&=&\lambda^n\left(B_\lambda(t)A^n+\sum\limits_{i=0}^{n}\frac{A^iT(t)A^n}{\lambda^{i+1}}\right)x=\\
&=&\lambda^n\left(A^nB_\lambda(t)+\sum\limits_{i=0}^{n}\frac{A^nA^iT(t)}{\lambda^{i+1}}\right)x=\\
&=&A^n\lambda^n\left(B_\lambda(t)+\sum\limits_{i=0}^{n}\frac{{T(t)}^{(i)}}{\lambda^{i+1}}\right)x=A^nB_\lambda(t)x\quad,\quad(\forall)n\in{\bf N}^*.
\end{eqnarray*}
Finalement, pour $x\in{\cal E}$ nous obtenons:
\begin{eqnarray*}
{B_\lambda(t)}^{(n)}{T(t)}^{(n)}x&=&{B_\lambda(t)}^{(n)}A^nT(t)x=A^n{B_\lambda(t)}^{(n)}T(t)x=\\
&=&A^nT(t){B_\lambda(t)}^{(n)}x={T(t)}^{(n)}{B_\lambda(t)}^{(n)}x\quad,\quad(\forall)n\in{\bf
N}^*.\fin
\end{eqnarray*}

\begin{teo}[spectral mapping]\label{num52}
Soient \dif et $A$ son g\'en\'erateur infinit\'esimal. Alors pour tout
$n\in{\bf N}^*$ on a:
$$
\left(e^{t\sigma(A)}\right)^{(n)}=\left\{\left.\lambda^ne^{\lambda
t}\right|\lambda\in\sigma(A)\right\}
\subseteq\sigma\left({T(t)}^{(n)}\right)\quad,\quad(\forall)t>0.
$$
\end{teo}
\dem
Pour $\lambda\in{\bf C}$ et $t>0$, nous consid\'erons l'op\'erateur:
$$
B_\lambda(t):{\cal E}\longrightarrow{\cal E}
$$
$$
B_\lambda(t)x=\int\limits_{0}^{t}\!\!e^{\lambda(t-s)}T(s)x\:ds\quad.
$$
Avec le lemme \ref{num30}, on d\'eduit que l'op\'erateur $B_\lambda(t)\in{\cal B(E)}$ est ind\'efiniment d\'erivable et:
$$
{B_\lambda(t)}^{(n)}=\lambda^n\left(B_\lambda(t)+\sum\limits_{i=0}^{n}\frac{{T(t)}^{(i)}}{\lambda^{i+1}}\right)\quad,\quad(\forall)n\in{\bf
N}^*.
$$
Compte tenu du lemme \ref{num23}, il r\'esulte que:
$$
(\lambda I-A)B_\lambda(t)x=e^{\lambda
t}x-T(t)x\quad,\quad(\forall)x\in{\cal E}
$$
et que:
$$
B_\lambda(t)(\lambda I-A)x=e^{\lambda
t}x-T(t)x\quad,\quad(\forall)x\in{\cal D}(A).
$$
Pour tout $n\in{\bf N}^*$ il s'ensuit que:
$$
(\lambda I-A){B_\lambda(t)}^{(n)}x=\lambda^ne^{\lambda
t}x-{T(t)}^{(n)}x\quad,\quad(\forall)x\in{\cal E}
$$
et:
$$
{B_\lambda(t)}^{(n)}(\lambda I-A)x=\lambda^ne^{\lambda
t}x-{T(t)}^{(n)}x\quad,\quad(\forall)x\in{\cal D}(A).
$$
Si $\lambda\in{\bf C}$ est tel que $\lambda^ne^{\lambda
t}\in\rho\left({T(t)}^{(n)}\right)$, alors on peut consid\'erer:
$$
Q=\left(\lambda^ne^{\lambda t}I-{T(t)}^{(n)}\right)^{-1}\in{\cal
B(E)}\quad,
$$
pour tout $n\in{\bf N}^*$.
Par cons\'equent:
$$
(\lambda I-A){B_\lambda(t)}^{(n)}Qx=x\quad,\quad(\forall)x\in{\cal E}
$$
et:
$$
Q{B_\lambda(t)}^{(n)}(\lambda I-A)x=x\quad,\quad(\forall)x\in{\cal
D}(A),
$$
pour tout $n\in{\bf N}^*$.
Mais, avec le lemme \ref{num30}, il r\'esulte:
$$
{B_\lambda(t)}^{(n)}{T(t)}^{(n)}={T(t)}^{(n)}{B_\lambda(t)}^{(n)}\quad,\quad(\forall)n\in{\bf N}^*.
$$
Donc:
$$
\lambda^ne^{\lambda
T}{B_\lambda(t)}^{(n)}-{B_\lambda(t)}^{(n)}{T(t)}^{(n)}=\lambda^ne^{\lambda
t}{B_\lambda(t)}^{(n)}-{T(t)}^{(n)}{B_\lambda(t)}^{(n)}
$$
et:
$$
{B_\lambda(t)}^{(n)}\left(\lambda^ne^{\lambda
t}I-{T(t)}^{(n)}\right)=\left(\lambda^ne^{\lambda
t}I-{T(t)}^{(n)}\right){B_\lambda(t)}^{(n)}\quad,
$$
pour tout $n\in{\bf N}^*$.
Par suite:
$$
{B_\lambda(t)}^{(n)}Q=Q{B_\lambda(t)}^{(n)}\quad,\quad(\forall)n\in{\bf N}^*
$$
et nous voyons que:
$$
(\lambda I-A){B_\lambda(t)}^{(n)}Qx=x\quad,\quad(\forall)x\in{\cal E}
$$
et:
$$
{B_\lambda(t)}^{(n)}Q(\lambda I-A)x=x\quad,\quad(\forall)x\in{\cal
D}(A)\quad,
$$
d'o\`u on obtient que $\lambda\in\rho(A)$. Nous en d\'eduisons que
$\lambda\in\sigma(A)$ implique $\lambda^ne^{\lambda
t}\in\sigma\left({T(t)}^{(n)}\right)$ pour tout $n\in{\bf N}^*$. Par cons\'equent:
$$
\left\{\left.\lambda^ne^{\lambda
t}\right|\lambda\in\sigma(A)\right\}\subset\sigma\left({T(t)}^{(n)}\right)\quad,
$$
ou bien:
$$
\left\{\left.\left(e^{\lambda
t}\right)^{(n)}\right|\lambda\in\sigma(A)\right\}\subset\sigma\left({T(t)}^{(n)}\right)
$$
et finalement:
$$
\left(e^{t\sigma(A)}\right)^{(n)}\subset\sigma\left({T(t)}^{(n)}\right)\quad,
$$
pour tout $n\in{\bf N}^*$ et tout $t>0$.\fin
\vspace{2cm}

    \section{$C_0$-semi-groupes analytiques}

\hspace{1cm}Par la suite nous \'etudions la possibilit\'e d'\'etendre
l'intervalle $]0,\infty)$ \`a une r\'egion du plan complexe, sans
abandonner les propri\'et\'es de $C_0$-semi-groupe. Nous d\'esignerons par $\Delta$ l'ensemble:
$$
\{z\in{\bf C}|Re\:z>0\mbox{ \em et }\varphi_1<\arg
z<\varphi_2\:,\:\varphi_1<0<\varphi_2\}
$$

\begin{definitie}
On appelle $C_0$-semi-groupe analytique une famille $\left\{T(z)\right\}_{z\in\Delta}\subset{\cal
B(E)}$ v\'erifiant les propri\'et\'es suivantes:\\
i) $T(0)=I$;\\
ii) $T(z_1+z_2)=T(z_1)T(z_2)$, $(\forall)z_1,z_2\in\Delta$;\\
iii) $\lim_{z\rightarrow 0}T(z)x=x$, $(\forall)x\in{\cal E}$,
$z\in\Delta$;\\
iv) l'application:
$$
\Delta\ni z\longmapsto T(z)\in{\cal B(E)}
$$
est analytique dans le secteur $\Delta$.
\end{definitie}

Comme la multiplication par $e^{\omega t}$ n'a aucun effet sur la
possibilit\'e ou l'impossibilit\'e d'extension \`a un semi-groupe
analytique, il est suffit de consid\'erer seulement les
$C_0$-semi-groupes uniform\'ement born\'es. Le th\'eor\`eme suivant donne
une caract\'erisation pour les $C_0$-semi-groupes analytiques
uniform\'ement born\'es.

\begin{teo}
Soient $\left\{T(t)\right\}_{t\geq 0}\in{\cal SG}(M,0)$ et $A$ son
g\'en\'erateur infinit\'esimal tel que $0\in\rho(A)$. Les affirmations suivantes sont \'equivalentes:\\
i) il existe $\delta>0$ tel que $\left\{T(t)\right\}_{t\geq 0}$ peut \^etre \'etendu \`a un
semi-groupe analytique dans le secteur:
$$
\Delta_\delta=\{z\in{\bf C}|Re\:z>0\mbox{ \em et }|\arg z|<\delta\}\quad,\quad\delta>0
$$
et $\left\{T(z)\right\}_{z\in{\Delta_{\delta'}}}$ est uniform\'ement
born\'e dans tout sous secteur ferm\'e
$\Delta_{\delta'}\subset\Delta_\delta$, o\`u
$\delta'\in]0,\delta[$;\\
ii) il existe une constante $C>0$ telle que pour tout $\gamma>0$ et tout $\eta\neq 0$ on ait:
$$
\left\|R(\gamma+i\eta;A)\right\|\leq\frac{C}{|\eta|}\quad;
$$
iii) il existe $\delta\in\left]0,\frac{\pi}{2}\right[$ et $K>1$ tel
que:
$$
\rho(A)\supset\Sigma_\delta=\left\{\lambda\in{\bf
C}\left||\arg\lambda|<\frac{\pi}{2}+\delta\right.\right\}\cup\{0\}
$$
et:
$$
\left\|R(\lambda;A)\right\|\leq\frac{K}{|\lambda|}\quad,\quad(\forall)\lambda\in\Sigma_\delta-\{0\};
$$
iv) l'application:
$$
]0,\infty)\ni t\longmapsto T(t)\in{\cal B(E)}
$$
est diff\'erentiable et il existe une constante $L>0$ tel que:
$$
\|AT(t)\|\leq\frac{L}{t}\quad,\quad(\forall)t>0.
$$
\end{teo}
\dem
$i)\Longrightarrow ii)$ Soit $\delta\in\left]0,\frac{\pi}{2}\right[$. Si $\delta'\in]0,\delta[$, il existe $C'>0$ tel
que:
$$
\|T(z)\|\leq C'
$$
pour tout $z\in\Delta_{\delta'}$. Comme l'application:
$$
\Delta_\delta\ni z\longmapsto T(z)\in{\cal B(E)}
$$
est analytique dans le sous secteur
$\Delta_{\delta'}$, avec le th\'eor\`eme de Cauchy (\cite[pag. 225]{dunford-schwartz}), on voit que
$$
\int\limits_{\Gamma}\!T(z)\:dz=0\quad,
$$
quel que soit le contour de Jordan lisse et ferm\'e $\Gamma\subset\Delta_{\delta'}$. Par cons\'equent, dans l'int\'egrale
$$
R(\gamma+i\eta;A)x=\int\limits_{0}^{\infty}\!e^{-(\gamma+i\eta)t}T(t)x\:dt\quad,\quad\gamma>0,
$$
on peut changer le chemin d'int\'egration par:
$$
\Gamma_\theta=\{r(\cos\theta+i\sin\theta)|\:0<r<\infty,|\theta|\leq\delta'\}\quad.
$$
Si $\eta>0$, alors pour le chemin
$$
\Gamma_{-\delta'}=\{r(\cos\delta'-i\sin\delta')|\:0<r<\infty\}\quad,
$$
nous obtenons:
\begin{eqnarray*}
&
&\|R(\gamma+i\eta;A)x\|=\left\|\int\limits_{0}^{\infty}\!e^{-(\gamma+i\eta)t}T(t)x\:dt\right\|=\\
&=&\left\|\:\int\limits_{\Gamma_{-\delta'}}\!e^{-(\gamma+i\eta)r(\cos\delta'-i\sin\delta')}T(r(\cos\delta'-i\sin\delta'))x\:d(r(\cos\delta'-i\sin\delta'))\right\|\leq\\
&\leq&\int\limits_{0}^{\infty}\left|e^{-(\gamma+i\eta)r(\cos\delta'-i\sin\delta')}\right|\left\|T(r(\cos\delta'-i\sin\delta')\right\|\|x\||\cos\delta'-i\sin\delta'|\:dr\leq\\
&\leq&C'\|x\|\int\limits_{0}^{\infty}\!e^{-r(\gamma\cos\delta'+\eta\sin\delta')}\:dr=
\frac{C'\|x\|}{\gamma\cos\delta'+\eta\sin\delta'}\leq\frac{C'}{\eta\sin\delta'}\|x\|\quad.
\end{eqnarray*}
Si nous notons:
$$
C=\frac{C'}{\sin\delta'}\quad,
$$
alors on obtient:
$$
\|R(\gamma+i\eta;A)\|\leq\frac{C}{\eta}\quad.
$$
Soit maintenant $\eta<0$. Alors pour le chemin
$$
\Gamma_{\delta'}=\{r(\cos\delta'+i\sin\delta')|\:0<r<\infty\}\quad,
$$
nous avons:
\begin{eqnarray*}
&
&\|R(\gamma+i\eta;A)x\|=\left\|\int\limits_{0}^{\infty}\!e^{-(\gamma+i\eta)t}T(t)x\:dt\right\|=\\
&=&\left\|\int\limits_{\Gamma_{\delta'}}\!e^{-(\gamma+i\eta)r(\cos\delta'+i\sin\delta')}T(r(\cos\delta'+i\sin\delta'))x\:d(r(\cos\delta'+i\sin\delta'))\right\|\leq\\
&\leq&\int\limits_{0}^{\infty}\left|e^{-(\gamma+i\eta)r(\cos\delta'+i\sin\delta')}\right|\left\|T(r(\cos\delta'+i\sin\delta')\right\|\|x\||\cos\delta'+i\sin\delta'|\:dr\leq\\
&\leq&C'\|x\|\int\limits_{0}^{\infty}\!e^{-r(\gamma\cos\delta'-\eta\sin\delta')}\:dr=
\frac{C'\|x\|}{\gamma\cos\delta'-\eta\sin\delta'}\leq\frac{C'}{-\eta\sin\delta'}\|x\|=\frac{C}{-\eta}\|x\|\quad.
\end{eqnarray*}
Par cons\'equent:
$$
\|R(\gamma+i\eta;A)\|\leq\frac{C}{-\eta}\quad.
$$
Finalement on voit que:
$$
\|R(\gamma+i\eta;A)\|\leq\frac{C}{|\eta|}\quad.
$$
$ii)\Longrightarrow iii)$ Comme $\left\{T(t)\right\}_{t\geq 0}\in{\cal
SG}(M,0)$, avec le th\'eor\`eme de Hille-Yosida on voit que:
$$
\|R(\lambda;A)\|\leq\frac{M}{\mbox{\em Re}\lambda}
$$
pour tout $\lambda\in\Lambda_0$. Compte tenu du (ii), il existe $C>0$
tel que:
$$
\|R(\lambda;A)\|\leq\frac{C}{|\mbox{\em Im}\lambda|}
$$
pour tout $\lambda\in\Lambda_0$ avec $\mbox{\em Im}\lambda\neq 0$. Compte tenu des in\'egalit\'es:
$$
\mbox{\em Re}\lambda\|R(\lambda;A)\|\leq M
$$
et:
$$
|\mbox{\em Im}\lambda|\|R(\lambda;A)\|\leq C\quad,
$$
nous obtenons:
$$
\mbox{\em Re}^2\lambda\|R(\lambda;A)\|^2\leq M^2
$$
et:
$$
\mbox{\em Im}^2\lambda\|R(\lambda;A)\|^2\leq C^2\quad,
$$
d'o\`u il r\'esulte que:
$$
\left(\mbox{\em Re}^2\lambda+\mbox{\em Im}^2\lambda\right)\|R(\lambda;A)\|^2\leq M^2+C^2\quad.
$$
Par cons\'equent, il existe une constante $K_1=\sqrt{M^2+C^2}>1$ tel que:
$$
\|R(\lambda;A)\|\leq\frac{K_1}{|\lambda|}\quad,\quad(\forall)\lambda\in\Lambda_0.
$$
Consid\'erons $\lambda\in{\bf C}$ avec $\mbox{\em
Re}\lambda\leq 0$. Soit $\gamma>0$ suffisamment petit et $\eta\neq
0$. Comme l'application $R(\:.\:;A)$ est
analytique sur $\rho(A)$, pour tout $\lambda=\mbox{\em
Re}\lambda+i\eta\in\rho(A)$ avec $\mbox{\em Re}\lambda\leq 0$, nous
obtenons:
$$
R(\lambda;A)=\sum\limits_{n=0}^{\infty}\frac{\left(\mbox{\em Re}\lambda-\gamma\right)^n}{n!}{R(\gamma+i\eta;A)}^{(n)}\quad.
$$
Avec la proposition \ref{num14} (iii), on voit que:
$$
{R(\gamma+i\eta;A)}^{(n)}=(-1)^nn!{R(\gamma+i\eta;A)}^{n+1}\quad.
$$
Alors il vient:
$$
R(\lambda;A)=\sum\limits_{n=0}^{\infty}(-1)^n(\mbox{\em
Re}\lambda-\gamma)^n{R(\gamma+i\eta;A)}^{n+1}
$$
et cette s\'erie est uniform\'ement convergente pour:
$$
\left\|R(\gamma+i\eta;A)\right\||\mbox{\em Re}\lambda-\gamma|\leq\alpha<1\quad.
$$
Compte tenu de la propri\'et\'e (ii), elle est convergente pour la topologie
de la norme si $\lambda=\mbox{\em Re}\lambda+i\eta\in\rho(A)$ est tel que sa partie r\'eelle v\'erifie $\mbox{\em Re}\lambda\leq 0$ et
$$
\left|\mbox{\em
Re}\lambda-\gamma\right|\leq\frac{\alpha|\eta|}{C}<\frac{|\eta|}{C}\quad.
$$
C'est-\`a-dire qu'il existe un voisinage ${\cal V}_\varepsilon$,
$\varepsilon=\frac{\alpha|\eta|}{C}$, de $\gamma+i\eta\in\Lambda_0$
contenu dans $\rho(A)$ lorsque $\gamma>0$ est suffisamment petit. Dans
ce voisinage ${\cal V}_\varepsilon$, il existe $\lambda\in{\bf C}$
tel que $\mbox{\em Re}\lambda\leq 0$ et $\lambda\in\rho(A)$. Si nous
d\'efinissons $\delta\in\left]0,\frac{\pi}{2}\right[$ tel que:
$$
\tan\delta=\frac{|\mbox{\em Re}\lambda|}{|\mbox{\em
Im}\lambda|}=\frac{|\mbox{\em
Re}\lambda|}{|\eta|}=\frac{\alpha}{C}\quad,\quad\alpha\in]0,1[,
$$
o\`u:
$$
\delta=\arctan\frac{\alpha}{C}\quad,\quad\alpha\in]0,1[,
$$
alors on voit que:
$$
\left\{\zeta\in{\bf
C}\left||\arg\zeta|<\frac{\pi}{2}+\delta\right.\right\}\subset\rho(A)\quad.
$$
Nous d\'esignerons par $\Sigma_\delta$ l'ensemble $\left\{\lambda\in{\bf
C}\left||\arg\lambda|<\frac{\pi}{2}+\delta\right.\right\}\cup\{0\}$, o\`u $\delta\in\left]0,\frac{\pi}{2}\right[$.
Si $\lambda\in\Sigma_\delta-\{0\}$ et $\mbox{\em
Re}\lambda\leq 0$, alors nous avons:
\begin{eqnarray*}
& &\|R(\lambda;A)\|\leq\sum\limits_{n=0}^{\infty}\left\|{R(\gamma+i\eta;A)}^{n+1}\right\||\mbox{\em
Re}\lambda-\gamma|^n\leq\\
&\leq&\sum\limits_{n=0}^{\infty}\frac{C^{n+1}}{|\eta|^{n+1}}\frac{\alpha^n|\eta|^n}{C^n}=\frac{1}{1-\alpha}\frac{C}{|\mbox{\em
Im}\lambda|}\quad.
\end{eqnarray*}
Comme:
$$
\frac{|\mbox{\em Re}\lambda|}{|\mbox{\em
Im}\lambda|}<\frac{1}{C}\quad,
$$
il vient:
$$
\frac{|\mbox{\em Re}\lambda|^2}{|\mbox{\em
Im}\lambda|^2}<\frac{1}{C^2}\quad,
$$
d'o\`u:
$$
\frac{|\mbox{\em Re}\lambda|^2}{|\mbox{\em
Im}^2\lambda|^2}+1<\frac{1}{C^2}+1\quad.
$$
Par cons\'equent
$$
\frac{|\mbox{\em Re}\lambda|^2+|\mbox{\em Im}\lambda|^2}{|\mbox{\em
Im}\lambda|^2}<\frac{1+C^2}{C^2}\quad,
$$
ou
$$
\frac{|\lambda|^2}{|\mbox{\em
Im}\lambda|^2}<\frac{1+C^2}{C^2}\quad,
$$
Il s'ensuit donc que:
$$
\frac{1}{|\mbox{\em Im}\lambda|}<\frac{\sqrt{1+C^2}}{C|\lambda|}\quad.
$$
Par suite:
$$
\|R(\lambda;A)\|\leq\frac{\sqrt{1+C^2}}{(1-\alpha)|\lambda|}\quad,\quad(\forall)\lambda\in\Sigma_\delta-\{0\}
$$
et si nous notons
$$
K_2=\frac{\sqrt{1+C^2}}{1-\alpha}>1\quad,
$$
alors il vient:
$$
\|R(\lambda;A)\|\leq\frac{K_2}{|\lambda|}\quad,\quad(\forall)\lambda\in\Sigma_\delta-\{0\}.
$$
Finalement, on obtient qu'il existe $K>1$ tel que:
$$
\|R(\lambda;A)\|\leq\frac{K}{|\lambda|}\quad,\quad(\forall)\lambda\in\Sigma_\delta-\{0\}.
$$
$iii)\Longrightarrow iv)$ Supposons qu'il existe
$\delta\in\left]0,\frac{\pi}{2}\right[$ et $K>1$ tel que:
$$
\Sigma_\delta=\left\{\lambda\in{\bf
C}\left||\arg\lambda|<\frac{\pi}{2}+\delta\right.\right\}\cup\{0\}
$$
et
$$
\|R(\lambda;A)\|\leq\frac{K}{|\lambda|}\quad,\quad(\forall)\lambda\in\Sigma_\delta-\{0\}.
$$
Compte tenu du th\'eor\`eme \ref{num45}, on voit que l'op\'erateur $A$
est le g\'en\'erateur infinit\'esimal d'un semi-groupe
$\left\{T(t)\right\}_{t\geq 0}$ pour lequel il existe $M>1$ tel que 
$$
\|T(t)\|\leq M\quad,\quad(\forall)t\geq 0.
$$
De plus, pour $\nu\in\left]\frac{\pi}{2},\frac{\pi}{2}+\delta\right[$ on consid\`ere le chemin
$$
\Gamma_\nu=\Gamma_\nu^{(1)}\cup\Gamma_\nu^{(2)}\quad,
$$
o\`u
$$
\Gamma_\nu^{(1)}=\left\{\left.r(\cos\nu-i\sin\nu)\right|0<r<\infty\right\}
$$
et
$$
\Gamma_\nu^{(2)}=\left\{\left.r(\cos\nu+i\sin\nu)\right|0<r<\infty\right\}\quad,
$$
tel que
$$
T(t)=\frac{1}{2\pi i}\int\limits_{\Gamma_\nu}\!e^{zt}R(z;A)\:dz\quad,\quad(\forall)t\geq 0,
$$
l'int\'egrale \'etant uniform\'ement convergente par rapport \`a $t>0$.\\
Soit 
$$
\varphi:[0,\infty)\times\Sigma_\delta\longrightarrow{\cal{B(E)}}\quad,
$$
$$
\varphi(t,z)=e^{zt}R(z;A)\quad.
$$
Il est clair que l'application $\varphi$ est diff\'erentiable par rapport \`a $t>0$ et:
$$
\frac{\partial\varphi(z,t)}{\partial t}=ze^{zt}R(z;A)\quad.
$$
De plus:
$$
\left\|\frac{\partial\varphi(z,t)}{\partial t}\right\|=\left\|ze^{zt}R(z;A)\right\|\leq K\left|e^{zt}\right|\quad,\quad(\forall)t>0.
$$
Avec le th\'eor\`eme de d\'erivation de Lebesgue, on voit que l'application
$$
]0,\infty)\ni t\mapsto T(t)\in{\cal{B(E)}}\quad,
$$
$$
T(t)=\frac{1}{2\pi i}\int\limits_{\Gamma_\nu}\!e^{zt}R(z;A)\:dz\quad,\quad(\forall)t\geq 0,
$$
est diff\'erentiable et on a:
$$
T(t)'=\frac{1}{2\pi i}\int\limits_{\Gamma_\nu}\!z e^{zt}R(z;A)\:dz\quad,\quad(\forall)t\geq 0.
$$
Comme $\nu\in\left]\frac{\pi}{2},\frac{\pi}{2}+\delta\right[$, il
vient $\cos\nu<0$ et compte tenu que:
\begin{eqnarray*}
& &\|T(t)'\|=\left\|\frac{1}{2\pi
i}\int\limits_{\Gamma_\nu^{(1)}}\!z e^{zt}R(z;A)\:dz+\frac{1}{2\pi
i}\int\limits_{\Gamma_\nu^{(2)}}\!z e^{zt}R(z;A)\:dz\right\|\leq\\
&\leq&\frac{1}{2\pi}\int\limits_{0}^{\infty}\!r e^{rt\cos\nu}\|R(z;A)\|\:dr+\frac{1}{2\pi}\int\limits_{0}^{\infty}\!r
e^{rt\cos\nu}\|R(z;A)\|\:dr\leq\\
&\leq&\frac{K}{\pi}\int\limits_{0}^{\infty}\!e^{-rt(-\cos\nu)}\:dr=\frac{K}{\pi}\left(\frac{1}{-t\cos\nu}\right)\quad,\quad(\forall)t>0,
\end{eqnarray*}
on d\'eduit qu'il existe $L=\frac{M}{\pi(-\cos\nu)}>0$ tel que:
$$
\|AT(t)\|=\|T(t)'\|\leq\frac{L}{t}\quad,\quad(\forall)t>0.
$$
$iv)\Longrightarrow i)$ Soit $t_0>0$. Avec la formule de Taylor et compte tenu de la remarque \ref{num46}, on a:
\begin{eqnarray*}
T(t)&=&\sum\limits_{k=0}^{\infty}\frac{(t-t_0)^k}{k!}T^{(k)}(t_0)+\frac{1}{(n-1)!}\int\limits_{t_0}^{t}\!(t-u)^{n-1}T^{(n)}(u)\:du=\\
&=&\sum\limits_{k=0}^{\infty}\frac{(t-t_0)^k}{k!}A^{k}T(t_0)+\frac{1}{(n-1)!}\int\limits_{t_0}^{t}\!(t-u)^{n-1}A^{n}T(u)\:du\quad,
\end{eqnarray*}
pour tout $n\in{\bf N}^{*}$.\\
Compte tenu du (iv) et de la remarque \ref{num46}, on voit que:
\begin{eqnarray*}
& &\left\|\frac{1}{(n-1)!}\int\limits_{t_0}^{t}\!(t-u)^{n-1}A^{n}T(u)\:du\right\|\leq\frac{1}{(n-1)!}\int\limits_{t_0}^{t}\!(t-u)^{n-1}\left\|\left[AT\left((\frac{u}{n}\right)\right]^n\right\|\:du\leq\\
&\leq&\frac{1}{(n-1)!}\int\limits_{t_0}^{t}\!(t-u)^{n-1}\left(\frac{Ln}{u}\right)^n\:du\leq\frac{1}{(n-1)!}\left(\frac{Ln}{t_0}\right)^n\int\limits_{t_0}^{t}\!(t-u)^{n-1}\:du=\\
&=&\frac{1}{(n-1)!}\left(\frac{Ln}{t_0}\right)^n\int\limits_{0}^{t-t_0}\!s^{n-1}\:du=\frac{1}{n!}\left(\frac{Ln}{t_0}\right)^n(t-t_0)^n\quad.
\end{eqnarray*}
Avec la formule de Stirling
$$
n!=n^n\sqrt{2\pi n}e^{-n+\frac{u_n}{12n}}\quad,\quad u_n\in]0,1[,
$$
on obtient
$$
n!e^n\geq n^n\quad.
$$
Par cons\'equent:
$$ 
\left\|\frac{1}{(n-1)!}\int\limits_{t_0}^{t}\!(t-u)^{n-1}A^{n}T(u)\:du\right\|\leq\frac{1}{n!}\frac{(t-t_0)^n}{t_0^n}L^nn!e^n=\left(Le\frac{t-t_0}{t_0}\right)^n\quad,
$$
pour $t\geq t_0>0$ et $n\in{\bf N}^{*}$ suffisamment grand.\\
Il en r\'esulte que la s\'erie de Taylor est convergente vers $T(t)$ si $t-t_0<\frac{t_0}{Le}$ et on a:
$$
T(t)=\sum\limits_{n=0}^{\infty}\frac{(t-t_0)^n}{n!}A^nT(t_0)\quad.
$$
Il s'ensuit donc que pour $z\in{\bf C}$ v\'erifiant 
$$
Re\:z>0\mbox{ et }\left|z-t_0\right|<\frac{t_0}{Le}\quad,
$$
on peut d\'efinir une fonction analytique
$$
T(z)=\sum\limits_{n=0}^{\infty}\frac{(z-t_0)^n}{n!}A^nT(t_0)\quad.
$$
La s\'erie de la partie droite de cette \'egalit\'e est uniform\'ement convergente par rapport \`a $z\in{\bf C}$ v\'erifiant les conditions $Re\:z>0$ et 
$$
\left|z-t_0\right|<\alpha\frac{t_0}{Le}\quad,
$$
o\`u $\alpha\in]0,1[$.\\
Soit $z\in{\bf C}$ tel que $Re\:z=t_0>$. On voit que:
$$
|Im\:z|=|z-t_0|<\frac{Re\:z}{Le}\quad,
$$
d'o\`u:
$$
\frac{|Im\:z|}{|Re\:z|}<\frac{1}{Le}
$$
ou encore
$$
|\arg z|\leq\arctan\left(\frac{1}{Le}\right)\quad.
$$
En prenant
$$\delta=\arctan\left(\frac{1}{Le}\right)\quad,
$$
nous obtenons que l'application
$$
\Delta_\delta\ni z\mapsto T(z)\in{\cal B(E)}
$$
est analytique dans le secteur
$$
\Delta_\delta=\{z\in{\bf C}|\:Re\:z>0\:et\:|\arg z|<\delta\}\quad.
$$
De plus, si nous consid\'erons $z\in{\bf C}$ avec les propri\'et\'es $Re\:z>0$ et $|z-t_0|\leq\alpha\frac{t_0}{Le}$, $\alpha\in]0,1[$, alors nous d\'eduisons que:
\begin{eqnarray*}
& &\|T(z)\|\leq\|T(t_0)\|+\sum\limits_{n=1}^{\infty}\frac{|z-t_0|^n}{n!}\left\|A^nT(t_0)\right\|\leq\\
&\leq&\|T(t_0)\|+\sum\limits_{n=1}^{\infty}\frac{\alpha^n}{n!}\left(\frac{t_0}{Le}\right)^n\left(\frac{Ln}{t_0}\right)^n\leq\\
&\leq&M+\sum\limits_{n=1}^{\infty}\alpha^n=M+\frac{\alpha}{1-\alpha}\quad.
\end{eqnarray*}
Par cons\'equent, si nous notons
$$
\delta^{'}=\arctan\left(\alpha\frac{1}{Le}\right)\quad,\quad\alpha\in]0,1[,
$$
nous voyons que l'application
$$
\Delta_\delta\ni z\mapsto T(z)\in{\cal B(E)}
$$
est uniform\'ement born\'ee dans le sous secteur
$$
\Delta_{\delta^{'}}=\{z\in{\bf C}|\:Re\:z>0\:et\:|\arg z|\leq\delta^{'}\}\subset\Delta_\delta\quad.
$$
Il est \'evident que $T(0)=I$ parce que $\left\{T(t)\right\}_{t\geq 0}\in{\cal SG}(M,0)$. De plus, pour tout $t>0$ et tout $z\in\Delta_\delta$, il r\'esulte que:
\begin{eqnarray*}
& &T(t)T(z)=\sum\limits_{n=0}^{\infty}\frac{(z-t_0)^n}{n!}A^nT(t_0+t)=\\
&=&\sum\limits_{n=0}^{\infty}\frac{[(z+t)-(t_0+t)]^n}{n!}A^nT(t_0+t)=T(t+z)\quad.
\end{eqnarray*}
Alors, pour tous $z_1,z_2\in\Delta_\delta$, nous obtenons:
\begin{eqnarray*}
& &T(z_1)T(z_2)=T(z_1)\sum\limits_{n=0}^{\infty}\frac{(z_2-t_0)^n}{n!}A^nT(t_0)=\\
&=&\sum\limits_{n=0}^{\infty}\frac{(z_2-t_0)^n}{n!}A^nT(z_1)T(t_0)=\sum\limits_{n=0}^{\infty}\frac{(z_2-t_0)^n}{n!}A^nT(z_1+t_0)=\\
&=&\sum\limits_{n=0}^{\infty}\frac{[(z_2+z_1)-(z_1+t_0)]^n}{n!}A^nT(z_1+t_0)=T(z_1+z_2)\quad.
\end{eqnarray*}
Soit 
$$
{\cal E}_0=\bigcup_{0<t<\infty}T(t){\cal E}\quad.
$$
Nous prouvons que cet ensemble est dense dans $\cal E$. Soient $x\in{\cal E}$ et $t_n>0$, $n\in{\bf N}$, tel que $\lim_{n\rightarrow\infty}t_n=0$. Alors pour $x_n=T(t_n)x\in{\cal E}_0$, $n\in{\bf N}$, nous obtenons:
$$
\lim_{n\rightarrow\infty}x_n=\lim_{n\rightarrow\infty}T(t_n)x=x\quad.
$$
Par cons\'equent, $\overline{{\cal E}_0}={\cal E}$.\\
De plus, nous avons vu que $\left\{T(z)\right\}_{z\in\Delta_{\delta^{'}}}$ est uniform\'ement born\'e dans tout sous secteur ferm\'e $\Delta_{\delta^{'}}$. De m\^eme, pour $x\in{\cal E}$ on obtient $T(t)x\in{\cal E}_0$ et:
$$
\lim_{z\rightarrow 0}T(z)T(t)x=\lim_{z\rightarrow 0}T(z+t)=T(t)x\quad.
$$
Compte tenu du th\'eor\`eme de Banach-Steinhaus (\cite[Theorem II.1.11, pag. 52]{dunford-schwartz}), il en r\'esulte que
$$
\lim_{z\rightarrow 0}T(z)x=x\quad,\quad(\forall)x\in{\cal E},\:z\in\Delta_\delta.
$$
Finalement, on voit que $\left\{T(z)\right\}_{z\in\Delta_{\delta}}$ est un $C_0$-semi-groupe analytique qui \'etend le semi-groupe $\left\{T(t)\right\}_{t\geq 0}\in{\cal SG}(M,0)$.\fin
\vspace{2cm}

    \section{$C_0$-semi-groupes de contractions}

\hspace{1cm}Dans la suite nous pr\'esentons quelques probl\`emes concernant la
classe du $C_0$-semi-groupes $\left\{T(t)\right\}_{t\geq 0}$
v\'erifiant la propri\'et\'e $\|T(t)\|\leq 1$, pour tout $t\geq 0$.

\begin{definitie}
On dit que $\left\{T(t)\right\}_{t\geq 0}$ est un $C_0$-semi-groupe de
contractions sur l'espace de Banch ${\cal E}$ si
$\left\{T(t)\right\}_{t\geq 0}\in{\cal SG}(1,0)$.
\end{definitie}

\begin{lema}
Soit $\left\{T(t)\right\}_{t\geq 0}\in{\cal SG}(M,\omega)$. Alors,
l'application:
$$
|||\:.\:|||:{\cal E}\longrightarrow{\bf R}_+
$$
$$
|||x|||=\sup_{t\geq 0}e^{-\omega
t}\|T(t)x\|\quad,\quad(\forall)x\in{\cal E},
$$
est une norme sur ${\cal E}$ \'equivalente avec la norme initiale
$\|\:.\:\|$.
\end{lema}
\dem
Soit $x\in{\cal E}$. Pour tout $t\geq 0$, on a:
$$
e^{-\omega t}\|T(t)x\|\leq e^{-\omega t}\|T(t)\|\|x\|\leq M\|x\|\quad.
$$
En passant \`a la borne sup\'erieure par rapport \`a $t$, on voit que:
$$
|||x|||\leq M\|x\|\quad,\quad(\forall)x\in{\cal E}.
$$
D'autre part, nous avons:
$$
|||x|||=\sup_{t\geq 0}e^{-\omega t}\|T(t)x\|\geq e^{-\omega
0}\|T(0)x\|=\|x\|\quad,\quad(\forall)x\in{\cal E},
$$
d'o\`u il r\'esulte que:
$$
\|x\|\leq|||x|||\leq M\|x\|\quad,\quad(\forall)x\in{\cal E}.
$$
Par cons\'equent les normes $|||\:.\:|||$ et $\|\:.\:\|$ sont
\'equivalentes.\fin

\begin{teo}
Soient $\left\{T(t)\right\}_{t\geq 0}\in{\cal SG}(M,\omega)$, $A$ son
g\'en\'erateur infinit\'esimal et:
$$
S(t)=e^{-\omega t}T(t)\quad,\quad(\forall)t\geq 0.
$$
Alors:\\
i) $\left\{S(t)\right\}_{t\geq 0}\in{\cal SG}(1,0);$\\
ii) le $C_0$-semi-groupe $\left\{S(t)\right\}_{t\geq 0}$ a pour
g\'en\'erateur infinit\'esimal l'op\'erateur $B=A-\omega I$.
\end{teo}
\dem
i) Il est clair que la famille $\left\{S(t)\right\}_{t\geq 0}$ est un
$C_0$-semi-groupe. De plus, pour tout $t\geq 0$, on a:
\begin{eqnarray*}
& &|||S(t)x|||=\sup_{s\geq 0}e^{-\omega s}\|T(s)e^{-\omega
t}T(t)x\|=\\
&=&\sup_{\tau\geq t}e^{-\omega\tau}\|T(\tau)x\|\leq\sup_{\tau\geq
0}e^{-\omega\tau}\|T(\tau)x\|=|||x|||\quad,\quad(\forall)x\in{\cal E},
\end{eqnarray*}
d'o\`u on obtient:
$$
|||S(t)x|||\leq|||x|||\quad,\quad(\forall)x\in{\cal E}\mbox{ et }t\geq
   0.
$$
Il s'ensuit que:
$$
\|S(t)\|\leq 1\quad,\quad(\forall)t\geq 0
$$
et, par cons\'equent, $\left\{S(t)\right\}_{t\geq 0}\in{\cal
SG}(1,0)$.\\
ii) Elle est analogue \`a celle du th\'eor\`eme \ref{num44}.\fin

Pour les $C_0$-semi-groupes de contractions, on peut formuler la version
suivante du th\'eor\`eme de Hille-Yosida.

\begin{teo}\label{num89}
Un op\'erateur lin\'eaire:
$$
A:{\cal D}(A)\subset{\cal E}\longrightarrow{\cal E}
$$
est le g\'en\'erateur infinit\'esimal d'un
semi-groupe $\left\{T(t)\right\}_{t\geq 0}\in{\cal SG}(1,0)$ si et seulement si:\\
i) $A$ est un op\'erateur ferm\'e et $\overline{{\cal D}(A)}={\cal
E}$;\\
ii) $\Lambda_0=\{\lambda\in{\bf
C}\left|\mbox{\em Re}\lambda>0\right.\}
\subset\rho(A)$ et pour
$\lambda\in\Lambda_0$, on a:
$$
\left\|{R(\lambda;A)}^n\right\|\leq\frac{1}{(\mbox{Re}\lambda)^n}\quad,\quad(\forall)n\in{\bf
N}^*.
$$
\end{teo}
\dem
$i)\Longrightarrow ii)$ Comme $\left\{T(t)\right\}_{t\geq 0}\in{\cal SG}(1,0)$, nous avons:
$$
\|T(t)\|\leq1\quad,\quad(\forall)t\geq 0.
$$
Par suite, on peut prendre $M=1$ et $\omega=0$. Avec le th\'eor\`eme de Hille-Yosida, il r\'esulte que:\\
(i) $A$ est un op\'erateur ferm\'e et $\overline{{\cal D}(A)}={\cal E}$;\\
(ii) $\Lambda_0=\{\lambda\in{\bf
C}\left|\mbox{\em Re}\lambda>0\right.\}
\subset\rho(A)$ et pour
$\lambda\in\Lambda_0$, on a:
$$
\left\|{R(\lambda;A)}^n\right\|\leq\frac{1}{(\mbox{Re}\lambda)^n}\quad,\quad(\forall)n\in{\bf
N}^*.
$$
$ii)\Longrightarrow i)$ Soit 
$$
A:{\cal D}(A)\subset{\cal E}\longrightarrow{\cal E}
$$
un op\'erateur lin\'eaire v\'erifiant les propri\'et\'es (i) et (ii) de l'\'enonc\'e. Avec le th\'eor\`eme de Hille-Yosida, il en r\'esulte que $A$ est le g\'en\'erateur infinit\'esimal d'un $C_0$-semi-groupe $\left\{T(t)\right\}_{t\geq 0}$ pour lequel il existe $M=1$ et $\omega=0$ tel que:
$$
\|T(t)\|\leq1\quad,\quad(\forall)t\geq 0.
$$
Donc $A$ est le g\'en\'erateur infinit\'esimal d'un semi-groupe de contractions.\fin

Une autre caract\'erisation tr\`es int\'eressante des
$C_0$-semi-groupes de contractions est donn\'ee par le fameux
th\'eor\`eme de Lumer-Phillips, dans lequel interviennent les
op\'erateurs $m$-dissipatifs.

\begin{definitie}
On appelle op\'erateur $m$-dissipatif un op\'erateur lin\'eaire
$A:{\cal D}(A)\subset{\cal E}\longrightarrow{\cal E}$ v\'erifiant les
propri\'et\'es suivantes:\\
i) $A$ est op\'erateur dissipatif;\\
ii) il existe $\alpha_0>0$ tel que ${\cal I}\mbox{m
}(\alpha_0I-A)={\cal E}$.
\end{definitie}

\begin{teo}[Lumer - Phillips]\label{num53}
Soit $A:{\cal D}(A)\subset{\cal E}\longrightarrow{\cal E}$ un
op\'erateur lin\'eaire tel que $\overline{{\cal D}(A)}={\cal E}$.\\
L'op\'erateur $A$ est le g\'en\'erateur infinit\'esimal d'un
semi-groupe $\left\{T(t)\right\}_{t\geq 0}\in{\cal SG}(1,0)$ si et
seulement si $A$ est un op\'erateur $m$-dissipatif.
\end{teo}
\dem
$\Longrightarrow$ Soit $A:{\cal D}(A)\subset{\cal
E}\longrightarrow{\cal E}$ le g\'en\'erateur infinit\'esimal d'un
semi-groupe $\left\{T(t)\right\}_{t\geq 0}\in{\cal SG}(1,0)$. Si
$x\in{\cal D}(A)$ et $x^*\in{\cal J}(x)$, alors pour tout $t\geq 0$ on
voit que:
\begin{eqnarray*}
& &\mbox{\em Re}\langle T(t)x-x,x^*\rangle=\mbox{\em Re}\langle
T(t)x,x^*\rangle-\mbox{\em Re}\langle x,x^*\rangle\leq\\
&\leq&\left|\langle T(t)x,x^*\rangle\right|-\mbox{\em Re}\langle
x,x^*\rangle\leq\|T(t)x\|\|x^*\|_*-\|x\|^2\leq\\
&\leq&\|x\|^2-\|x\|^2=0\quad,
\end{eqnarray*}
d'o\`u:
$$
\frac{1}{t}\mbox{\em Re}\langle T(t)x-x,x^*\rangle\leq
0\quad,\quad(\forall)t\geq 0.
$$
En passant \`a limite pour $t\searrow 0$, il vient:
$$
\mbox{\em Re}\langle Ax,x^*\rangle\leq 0
$$
pour tout $x\in{\cal E}$ et tout $x^*\in{\cal J}(x)$. Il s'ensuit donc que $A$
est un op\'erateur dissipatif.\\
D'autre part, avec le th\'eor\`eme \ref{num89}  , on voit que $]0,\infty)\subset\rho(A)$. Donc $\alpha I-A\in{\cal GL(E)}$, $(\forall)\alpha\in]0,\infty)$. Par suite, ${\cal
I}\mbox{\em m }(\alpha I-A)={\cal E}$ pour tout $\alpha>0$ et finalement on voit que $A$ est un op\'erateur $m$-dissipatif.\\
ii) Soit $A:{\cal D}(A)\subset{\cal
E}\longrightarrow{\cal E}$ un op\'erateur $m$-dissipatif tel que
$\overline{{\cal D}(A)}={\cal E}$. Alors, il existe
$\alpha_0\in]0,\infty)$ tel que ${\cal I}\mbox{\em m
}(\alpha_0I-A)={\cal E}$. En appliquant la proposition \ref{num35}, on
voit que ${\cal I}\mbox{\em m }(\alpha I-A)={\cal E}$, pour tout
$\alpha\in]0,\infty)$. Il s'ensuit que
$]0,\infty)\subset\rho(A)$. Comme $A$ est un op\'erateur dissipatif,
compte tenu de la proposition \ref{num34}, il vient:
$$
\left\|(\alpha
I-A)x\right\|\geq\alpha\|x\|\quad,\quad(\forall)x\in{\cal D}(A),
$$
d'o\`u:
$$
\left\|R(\alpha;A)\right\|\leq\frac{1}{\alpha}
$$
pour tout $\alpha\in]0,\infty)$. Avec le th\'eor\`eme \ref{num89},
on voit que $A$ est le g\'en\'erateur infinit\'esimal d'un semi-groupe $\left\{T(t)\right\}_{t\geq 0}\in{\cal SG}(1,0)$.\fin
\begin{prop}
Soit $A\in{\cal B(E)}$ tel que $\|A\|\leq 1$. Alors
$\left\{e^{t(A-I)}\right\}_{t\geq 0}$ est un semi-groupe uniform\'ement
continu de contractions.
\end{prop}
\dem
Il est \'evident que $\left\{e^{t(A-I)}\right\}_{t\geq 0}$ est un
semi-groupe uniform\'ement continu. De plus:
$$
\left\|e^{t(A-I)}\right\|\leq\left\|e^{tA}\right\|\:\left\|e^{-t}\right\|\leq
e^{t\|A\|}e^{-t}\leq
1\quad,\quad(\forall)t\geq 0.\fin
$$
\vspace{2cm}

     \section{Notes}

{\footnotesize
Pour les propri\'et\'es des $C_0$-semi-groupes diff\'erentiables nous
avons consult\'e \cite[pag. 51]{pazy1}, \cite[pag. 73]{ahmed} et
\cite[pag. 28]{davies}. Le th\'eor\`eme \ref{num52} se trouve dans
\cite{lemle1}.

Les propri\'et\'es des $C_0$-semi-groupes analytiques uniform\'ement
born\'es se trouvent dans \cite[pag. 60]{pazy1}, \cite[pag. 81]{ahmed}
ou \cite[pag.59]{davies}. Une introduction tr\`es int\'eressante des $C_0$-semi-groupes analytiques, par la construction d'un
calcul fonctionnel ad\'equat, est donn\'e dans
\cite[pag. 121]{clement}.

Le th\'eor\`eme \ref{num53} a \'et\'e montr\'e par Lummer et Phillips
dans \cite{lumer-phillips}.
}

  \chapter{La formule de Lie - Trotter}

     \section{Le cas des semi-groupes uniform\'ement continus}

\hspace{1cm}Dans cette section nous pr\'esentons la formule du produit de Lie-Trotter pour les
semi-groupes uniform\'ement continus.

\begin{teo}\label{num43}
Soient $A$ le g\'en\'erateur infinit\'esimal d'un semi-groupe
uniform\'e-\\ment continu $\left\{T(t)\right\}_{t\geq 0}$ et
$\left(\left\{A_n(t)\right\}_{t\geq 0}\right)_{n\in{\bf
N}}\subset{\cal B(E)}$ tel que:
$$
\lim_{n\rightarrow\infty}\|A_n(t)\|=0\quad,
$$
uniform\'ement par rapport \`a $t$ sur les intervalles compacts de
$[0,\infty)$.
Alors:
$$
T(t)=\lim_{n\rightarrow\infty}\left[I+\frac{t}{n}\left(A+A_n(t)\right)\right]^n\quad,
$$
uniform\'ement par rapport \`a $t$ sur les intervalles compacts de
$[0,\infty)$.
\end{teo}
\dem
Soient $0\leq a<b$. Si nous notons:
$$
V_n(t)=I+\frac{t}{n}\left(A+A_n(t)\right)\quad,
$$
alors, pour tout $0\leq k\leq n$, on a:
\begin{eqnarray*}
&
&\left\|V_n^k(t)\right\|=\left\|\left[I+\frac{t}{n}\left(A+A_n(t)\right)\right]^k\right\|\leq\\
&\leq&\left[1+\frac{t}{n}\left(\|A\|+\|A_n(t)\|\right)\right]^k\leq\left[1+\frac{t}{n}\left(\|A\|+\|A_n(t)\|\right)\right]^n=\\
&=&\sum\limits_{i=0}^{n}\frac{t^i\left(\|A\|+\|A_n(t)\|\right)^i}{n^i}\leq\sum\limits_{i=0}^{n}\frac{t^i\left(\|A\|+\|A_n(t)\|\right)^i}{i!}\leq\\
&\leq&\sum\limits_{i=0}^{\infty}\frac{t^i\left(\|A\|+\|A_n(t)\|\right)^i}{i!}=e^{t\left(\|A\|+\|A_n(t)\|\right)}\leq
M\quad,
\end{eqnarray*}
cette derni\`ere quantit\'e etant uniform\'ement born\'ee par rapport \`a $t\in[a,b]$.
De m\^eme, pour:
$$
U_n(t)=e^{\frac{t}{n}A}\quad,\quad(\forall)t\geq 0,
$$
nous obtenons:
$$
U_n^n(t)=e^{tA}=T(t)\quad,\quad(\forall)t\geq 0
$$
et pour tout $0\leq k\leq n$, nous avons:
$$
\left\|U_n^k(t)\right\|\leq e^{\frac{kt}{n}\|A\|}\leq e^{\frac{nt}{n}\|A\|}\leq
e^{t\|A\|}\leq
N\quad,
$$
uniform\'ement par rapport \`a $t\in[a,b]$. Il vient:
\begin{eqnarray*}
& &V_n^n(t)-U_n^n(t)=V_n^n(t)U_n^0(t)-V_n^{n-1}(t)U_n^1(t)+\\
&+&V_n^{n-1}(t)U_n^1(t)-V_n^{n-2}(t)U_n^2(t)+\\
&+&V_n^{n-2}(t)U_n^2(t)-\cdots-V_n^0(t)U_n^n(t)=\\
&=&\sum\limits_{i=0}^{n-1}\left[V_n^{n-i}(t)U_n^i(t)-V_n^{n-i-1}(t)U_n^{i+1}(t)\right]=\\
&=&\sum\limits_{i=0}^{n-1}V_n^{n-i-1}(t)\left[V_n(t)-U_n(t)\right]U_n^i(t)\quad.
\end{eqnarray*}
Comme:
\begin{eqnarray*}
& &V_n(t)-U_n(t)=I+\frac{t}{n}A+\frac{t}{n}A_n(t)-e^{\frac{t}{n}A}=\\
&=&I+\frac{t}{n}A+\frac{t}{n}A_n(t)-I-\frac{t}{n}A-\frac{1}{2!}\frac{t^2}{n^2}A^2-\cdots\quad,
\end{eqnarray*}
il r\'esulte que:
$$
\left\|V_n(t)-U_n(t)\right\|\leq\frac{t}{n}\|A_n(t)\|+\frac{1}{2!}\frac{t^2}{n^2}\|A\|^2+\cdots\quad.
$$
Par cons\'equent:
\begin{eqnarray*}
&
&\left\|V_n^n(t)-U_n^n(t)\right\|\leq\sum\limits_{i=0}^{n-1}M\left(\frac{t}{n}\|A_n(t)\|+\frac{1}{2!}\frac{t^2}{n^2}\|A\|^2+\cdots\right)N=\\
&=&MN\left(t\|A_n(t)\|+\frac{1}{2!}\frac{t^2}{n}\|A\|^2+\cdots\right)\longrightarrow
0\quad\mbox{si}\quad n\rightarrow\infty,
\end{eqnarray*}
uniform\'ement par rapport \`a $t$ sur les intervalles compacts de
$[0,\infty)$, ce qui ach\`eve la preuve.\fin

\begin{teo}[la formule exponentielle]
Soit $A$ le g\'en\'erateur infinit\'esi-\\mal d'un semi-groupe
uniform\'ement continu $\left\{T(t)\right\}_{t\geq 0}$. Alors:
$$
T(t)=\lim_{n\rightarrow\infty}\left(I+\frac{t}{n}A\right)^n=\lim_{n\rightarrow\infty}\left(I-\frac{t}{n}A\right)^{-n}=\lim_{n\rightarrow\infty}\left[\frac{n}{t}R\left(\frac{n}{t};A\right)\right]^n\quad,
$$
uniform\'ement par rapport \`a $t$ sur les intervalles compacts de
$[0,\infty)$.
\end{teo}
\dem
La premi\`ere \'egalit\'e r\'esulte du th\'eor\`eme \ref{num43}
pour $A_n(t)=0$, quels que soient $n\in{\bf N}$ et $t\geq 0$.\\
Soient $0\leq a<b$. On a:
$$
\left\|\frac{t}{n}A\right\|<1
$$
pour $n$ suffisamment grand et $t\in[a,b]$. Avec le lemme \ref{num2}, il vient:
$$
I-\frac{t}{n}A\in{\cal GL(E)}
$$
et:
$$
\left(I-\frac{t}{n}A\right)^{-1}=\sum\limits_{i=0}^{\infty}\frac{t^iA^i}{n^i}=I+\frac{t}{n}\left(A+A_n(t)\right)\quad,
$$
o\`u:
$$
A_n(t)=\frac{t}{n}A^2+\frac{t^2}{n^2}A^3+\cdots
$$
et:
$$
\lim_{n\rightarrow\infty}\|A_n(t)\|=0\quad,
$$
uniform\'ement par rapport \`a $t\in[a,b]$. Avec le th\'eor\`eme \ref{num43}, on voit que:
$$
T(t)=\lim_{n\rightarrow\infty}\left(I+\frac{t}{n}(A+A_n(t))\right)^n=\lim_{n\rightarrow\infty}\left[I-\frac{t}{n}A\right]^{-n}\quad,
$$
uniform\'ement par rapport \`a $t$ sur les intervalles compacts de
$[0,\infty)$.\\
La troisi\`eme \'egalit\'e en r\'esulte compte tenu que:
$$
\left(I-\frac{t}{n}A\right)^{-n}=\left[\frac{n}{t}\left(\frac{n}{t}I-A\right)^{-1}\right]^n=\left[\frac{n}{t}R\left(\frac{n}{t};A\right)\right]^n\quad.\fin
$$

\begin{teo}[la formule de Lie-Trotter]
Soit $A_1$ le g\'en\'erateur infinit\'e-\\simal du semi-groupe
uniform\'ement continu $\left\{T_1(t)\right\}_{t\geq 0}$ et $A_2$ le g\'en\'erateur infinit\'esimal du semi-groupe
uniform\'ement continu $\left\{T_2(t)\right\}_{t\geq 0}$, alors l'op\'erateur:
$$
A:{\cal E}\longrightarrow{\cal E}\quad,
$$
$$
Ax=A_1x+A_2x
$$
est le g\'en\'erateur infinit\'esimal d'un semi-groupe
$\left\{T(t)\right\}_{t\geq 0}$ uniform\'ement
continu, tel que:
$$
T(t)=\lim_{n\rightarrow\infty}\left[T_1\left(\frac{t}{n}\right)T_2\left(\frac{t}{n}\right)\right]^n\quad,
$$
uniform\'ement par rapport \`a $t$ sur les intervalles compacts de
$[0,\infty)$.
\end{teo}
\dem
Nous avons successivement:
\begin{eqnarray*}
&
&T_1\left(\frac{t}{n}\right)T_2\left(\frac{t}{n}\right)=e^{\frac{t}{n}A_1}e^{\frac{t}{n}A_2}=\\
&=&\left[I+\frac{t}{n}A_1+\frac{1}{2!}\frac{t^2}{n^2}A_1^2+\frac{1}{3!}\frac{t^3}{n^3}A_1^3+\cdots\right]\left[I+\frac{t}{n}A_2+\frac{1}{2!}\frac{t^2}{n^2}A_2^2+\frac{1}{3!}\frac{t^3}{n^3}A_2^3+\cdots\right]=\\
&=&I+\frac{t}{n}(A_1+A_2)+\frac{t^2}{n^2}\left(\frac{1}{2!}A_1^2+A_1A_2+\frac{1}{2!}A_2^2\right)+\cdots=\\
&=&I+\frac{t}{n}\left[(A_1+A_2)+\frac{t}{n}\left(\frac{1}{2!}A_1^2+A_1A_2+\frac{1}{2!}A_2^2\right)+\cdots\right]=\\
&=&I+\frac{t}{n}\left[A+A_n(t)\right]\quad,
\end{eqnarray*}
o\`u l'op\'erateur:
$$
A_n(t)=\frac{t}{n}\left(\frac{1}{2!}A_1^2+A_1A_2+\frac{1}{2!}A_2^2\right)+\frac{t^2}{n^2}\left(\frac{1}{3!}A_1^3+\frac{1}{2!}A_1^2A_2+\frac{1}{2!}A_1A_2^2+\frac{1}{3!}A_2^3\right)+\cdots
$$
a la propri\'et\'e:
\begin{eqnarray*}
\|A_n(t)\|&\leq&\frac{t}{n}\left\|\frac{1}{2!}A_1^2+A_1A_2+\frac{1}{2!}A_2^2\right\|+\\
&+&\frac{t^2}{n^2}\left\|\frac{1}{3!}A_1^3+\frac{1}{2!}A_1^2A_2+\frac{1}{2!}A_1A_2^2+\frac{1}{3!}A_2^3\right\|+\cdots\longrightarrow 0\quad\mbox{ si }n\rightarrow\infty\quad,
\end{eqnarray*}
uniform\'ement par rapport \`a $t$ sur les intervalles compacts de
$[0,\infty)$. Avec le th\'eor\`eme \ref{num43}, on voit que:
$$
T(t)=\lim_{n\rightarrow\infty}\left[T_1\left(\frac{t}{n}\right)T_2\left(\frac{t}{n}\right)\right]^n\quad,
$$
uniform\'ement par rapport \`a $t$ sur les intervalles compacts de
$[0,\infty)$.\fin

\begin{obs}
\em
Si $A,B\in{\cal B(E)}$, alors on a:
$$
e^{t(A+B)}=\lim_{n\rightarrow\infty}\left(e^{\frac{t}{n}A}e^{\frac{t}{n}B}\right)^n\quad,
$$
uniform\'ement par rapport \`a $t$ sur les intervalles compacts de
$[0,\infty)$.
\end{obs}
\vspace{2cm}

    \section{Propri\'et\'es de convergence des $C_0$-semi-groupes}

\hspace{1cm}Dans cette section on introduit la topologie de la r\'esolvante
sur l'ensemble ${\cal GI}({\cal E})$ des g\'en\'erateurs
infinit\'esimaux et on montre le th\'eor\`eme de Trotter - Kato.

\begin{definitie}
On dit que la suite $\left(A_n\right)_{n\in{\bf N}^*}\subset{\cal
GI}({\cal E})$ est convergente vers $A\in{\cal GI}({\cal E})$ pour la
topologie forte de la r\'esolvante si pour tout
$\lambda\in\bigcap\limits_{n\in{\bf N}^*}\rho(A_n)\cap\rho(A)$, on a:
$$
R(\lambda;A_n)x\longrightarrow
R(\lambda;A)x\quad,\quad\mbox{si}\quad n\rightarrow\infty\;,\;(\forall)x\in{\cal E}.
$$
De m\^eme, on dit que la suite $\left(A_n\right)_{n\in{\bf N}^*}\subset{\cal
GI}({\cal E})$ est convergente
vers $A\in{\cal GI}({\cal E})$ pour la topologie uniforme de la
r\'esolvante si pour tout $\lambda\in\bigcap\limits_{n\in{\bf N}^*}\rho(A_n)\cap\rho(A)$, on a:
$$
\left\|R(\lambda;A_n)-R(\lambda;A)\right\|\longrightarrow
0\quad,\quad\mbox{si}\quad n\rightarrow\infty.
$$
\end{definitie}

Par la suite, nous supposerons que ${\cal GI}({\cal E})$ est dot\'e de la
topologie forte de la r\'esolvante.

\begin{lema}\label{num18}
Soient $\left\{T(t)\right\}_{t\geq 0}\;,\;\left\{S(t)\right\}_{t\geq
0}\in{{\cal SG}(M,\omega)}$ et $A$, respectivement $B$, leur
g\'en\'erateurs infinit\'esimaux. Alors pour tout
$\lambda\in\Lambda_\omega$ et tout $x\in{\cal E}$ on a l'\'egalit\'e:
$$
R(\lambda;B)\left[T(t)-S(t)\right]R(\lambda;A)x=\int\limits_{0}^{t}\!S(t-s)\left[R(\lambda;A)-R(\lambda;B)\right]T(s)x\:ds
$$
quel que soit $t\geq 0$.
\end{lema}
\dem
Soient $x\in{\cal E}$ et $\lambda\in\Lambda_\omega$. Alors
$R(\lambda;A)x\in{\cal D}(A)$ et $R(\lambda;B)x\in{\cal D}(B)$. L'application:
$$
[0,t]\ni s\longrightarrow S(t-s)R(\lambda;B)T(s)R(\lambda;A)x\in{\cal
E}
$$
est diff\'erentiable et pour $s\in[0,t]$ et $x\in{\cal E}$ nous avons:
\begin{eqnarray*}
& &\frac{d}{ds}S(t-s)\left[R(\lambda;B)T(s)R(\lambda;A)x\right]=\\
&=&S(t-s)(-B)R(\lambda;B)T(s)R(\lambda;A)x+S(t-s)R(\lambda;B)T(s)AR(\lambda;A)x=\\
&=&S(t-s)(\lambda I-B-\lambda I)R(\lambda;B)T(s)R(\lambda;A)x+\\
&+&S(t-s)R(\lambda;B)T(s)(-\lambda I+A+\lambda I)R(\lambda;A)x=\\
&=&S(t-s)T(s)R(\lambda;A)x-\lambda
S(t-s)R(\lambda;B)T(s)R(\lambda;A)x-\\
&-&S(t-s)R(\lambda;B)T(s)x+\lambda S(t-s)R(\lambda;B)T(s)R(\lambda;A)x=\\
&=&S(t-s)\left[T(s)R(\lambda;A)x-R(\lambda;B)T(s)x\right]=\\
&=&S(t-s)\left[R(\lambda;A)-R(\lambda;B)\right]T(s)x
\end{eqnarray*}
puisque la r\'esolvante $R(\lambda;A)$ commute avec $T(t)$,
$(\forall)t\geq 0$. Par cons\'equent:
$$
S(t-s)R(\lambda;B)T(s)R(\lambda;A)x|_0^t=\int\limits_{0}^{t}\!S(t-s)\left[R(\lambda;A)-R(\lambda;B)\right]T(s)x\:ds\quad,
$$
ou encore:
\begin{eqnarray*}
& &R(\lambda;B)T(t)R(\lambda;A)x-S(t)R(\lambda;B)R(\lambda;A)x=\\
&=&\int\limits_{0}^{t}\!S(t-s)\left[R(\lambda;A)-R(\lambda;B)\right]T(s)x\:ds\quad,\quad(\forall)x\in{\cal
E}.
\end{eqnarray*}
Comme $S(t)R(\lambda;B)=R(\lambda;B)S(t)$ pour tout $t\geq 0$, on obtient finalement:
$$
R(\lambda;B)\left[T(t)-S(t)\right]R(\lambda;A)x=\int\limits_{0}^{t}\!S(t-s)\left[R(\lambda;A)-R(\lambda;B)\right]T(s)x\:ds
$$
pour tout $x\in{\cal E}$.\fin

Le th\'eor\`eme suivant pr\'esente une tr\`es jolie correspondance entre les $C_0$-semi-groupes d'op\'erateurs lin\'eaires born\'es et
leur g\'en\'erateurs infinit\'esimaux.

\begin{teo}\label{num20}
Soient $\left(\left\{T_n(t)\right\}_{t\geq 0}\right)_{n\in{\bf
N}^*}\subset{{\cal SG}(M,\omega)}$ ayant pour g\'en\'erateurs
infinit\'esimaux les op\'erateurs $(A_n)_{n\in{\bf N}^*}\subset{{\cal GI}(M,\omega)}$
et \semi ayant pour g\'en\'erateur infinit\'esimal l'op\'erateur $A\in{{\cal
GI}(M,\omega)}$.\\ Les affirmations suivantes sont \'equivalentes:\\
i) $A_n\longrightarrow A$, si $n\rightarrow\infty$, pour la
topologie forte de la r\'esolvante;\\
ii) pour tout $t_0\in]0,\infty)$ on a l'\'egalit\'e:
$$
\lim_{n\rightarrow\infty}\left\{\sup_{t\in[0,t_0]}\left\|T_n(t)x-T(t)x\right\|\right\}=0\quad,\quad(\forall)x\in{\cal
E}.
$$
\end{teo}
\dem
$i)\Longrightarrow ii)$
Supposons que $A_n\longrightarrow A$, si $n\rightarrow\infty$, pour la
topologie forte de la r\'esolvante. Pour tout
$\lambda\in\Lambda_\omega$, nous avons:
$$
R(\lambda;A_n)x\longrightarrow
R(\lambda;A)x\quad,\quad\mbox{si }n\rightarrow\infty\;,\;(\forall)x\in{\cal E}.
$$
Soient $t_0\in]0,\infty)$, $x\in{\cal E}$ et
$\lambda\in\Lambda_\omega$ arbitrairement fix\'ees. Puisque la
r\'esolvante commute avec le semi-groupe associ\'e, il r\'esulte que:
\begin{eqnarray*}
&
&\left[T_n(t)-T(t)\right]R(\lambda;A)x=T_n(t)\left[R(\lambda;A)-R(\lambda;A_n)\right]x+\\
&+&R(\lambda;A_n)\left[T_n(t)-T(t)\right]x+\left[R(\lambda;A_n)-R(\lambda;A)\right]T(t)x.
\end{eqnarray*}
Montrons que cette expression tend vers zero si
$n\rightarrow\infty$.\\
Comme $\left(\left\{T_n(t)\right\}_{t\geq 0}\right)_{n\in{\bf
N}^*}\subset{{\cal SG}(M,\omega)}$, il est clair que:
$$
\left\|T_n(t)\right\|\leq Me^{\omega t_0}\quad,\quad(\forall)t\in[0,t_0].
$$
Compte tennu de (i), on voit que le premier terme converge vers zero si $n\rightarrow\infty$,
uniform\'ement par rapport \`a $t\in[0,t_0]$.\\
De m\^eme, la continuit\'e de l'application $t\mapsto T(t)x$ sur l'intervalle
compact $[0,t_0]$, conduit au fait que l'ensemble
$\left\{T(t)x\left|t\in[0,t_0]\right.\right\}$ est compact, comme l'image d'un compact par une fonction continue. On en
d\'eduit facilement que le troisi\`eme terme est fortement convergent vers zero lorsque $n\rightarrow\infty$ et cette convergence est uniforme par rapport \`a $t\in[0,t_0]$.\\
Pour le deuxi\`eme terme, compte tenu du lemme \ref{num18}, on a:
$$
R(\lambda;A_n)\left[T(t)-T_n(t)\right]R(\lambda;A)x=\int\limits_{0}^{t}\!T_n(t-s)\left[R(\lambda;A)-R(\lambda;A_n)\right]T(s)x\:ds,
$$
pour tout $t\in[0,t_0]$.
Si pour $s\in[0,t_0]$, on pose
$$
f_{t,n}(s)x=T_n(t-s)\left[R(\lambda;A)-R(\lambda;A_n)\right]T(s)x\quad,\quad
0\leq s\leq t\leq t_0,
$$
alors on voit que:
\begin{eqnarray*}
\left\|f_{t,n}(s)x\right\|&=&\left\|T_n(t-s)\left[R(\lambda;A)-R(\lambda;A_n)\right]T(s)x\right\|\leq\\
&\leq&\left\|T_n(t-s)\right\|\left\|R(\lambda;A)-R(\lambda;A_n)\right\|\left\|T(s)\right\|\|x\|\leq\\
&\leq&Me^{\omega(t-s)}\left(\left\|R(\lambda;A)\right\|+\left\|R(\lambda;A_n)\right\|\right)Me^{\omega s}\|x\|\leq\\
&\leq&Me^{\omega t}\left(\frac{M}{Re\:\lambda-\omega}+\frac{M}{Re\:\lambda-\omega}\right)Me^{\omega t}\|x\|=\\
&=&\frac{2M^3}{Re\:\lambda-\omega}e^{2\omega t}\|x\|\quad.
\end{eqnarray*}
De plus, compte tenu des in\'egalit\'es:
\begin{eqnarray*}
\left\|f_{t,n}(s)x\right\|&=&\left\|T_n(t-s)\left[R(\lambda;A)-R(\lambda;A_n)\right]T(s)x\right\|\leq\\
&\leq&\left\|T_n(t-s)\right\|\left\|R(\lambda;A)-R(\lambda;A_n)\right\|\left\|T(s)\right\|\|x\|\leq\\
&\leq&Me^{\omega(t-s)}\left\|R(\lambda;A)-R(\lambda;A_n)\right\|Me^{\omega s}\|x\|\leq\\
&\leq&M^2e^{\omega t}\left\|R(\lambda;A)-R(\lambda;A_n)\right\|\|x\|\quad,
\end{eqnarray*}
nous obtenons
$$
\lim_{n\rightarrow\infty}\left\|f_{t,n}(s)x\right\|=0\quad,
$$
quels que soient $s\in[0,t]$ et $t\in[0,t_0]$.
Avec le th\'eor\`eme de la convergence domin\'ee de Lebesgue (\cite[Theorem III.3.7,
pag. 124]{dunford-schwartz}) il r\'esulte que:
$$
\lim_{n\rightarrow\infty}\int\limits_{0}^{t}\left\|f_{t,n}(s)x\right\|\:ds=\int\limits_{0}^{t}\lim_{n\rightarrow\infty}\left\|f_{t,n}(s)x\right\|\:ds\quad,
$$
pour tout $t\in[0,t_0]$. Il s'ensuit donc que:
$$
\lim_{n\rightarrow\infty}\left\|R(\lambda;A_n)\left[T(t)-T_n(t)\right]R(\lambda;A)x\right\|=0\quad,\quad(\forall)
x\in{\cal E},
$$
uniform\'ement par rapport \`a $t\in[0,t_0]$. Si
nous notons $y=R(\lambda;A)x\in{\cal D}(A)$, on voit que:
$$
\lim_{n\rightarrow\infty}\left\|R(\lambda;A_n)\left[T(t)-T_n(t)\right]y\right\|=0\quad,\quad(\forall)y\in{\cal
D}(A),
$$
uniform\'ement par rapport \`a $t\in[0,t_0]$. Par cons\'equent, si
$x\in{\cal D}(A)$, le deuxi\`eme terme tend vers zero pour
$n\rightarrow\infty$, uniform\'ement par rapport \`a $t\in[0,t_0]$.\\
Il s'ensuit que:
$$
\lim_{n\rightarrow\infty}\left\{\sup_{t\in[0,t_0]}\left\|\left[T_n(t)-T(t)\right]R(\lambda;A)x\right\|\right\}=0\quad,\quad(\forall)x\in{\cal
D}(A),
$$
d'o\`u il r\'esulte imm\'ediatement:
$$
\lim_{n\rightarrow\infty}\left\{\sup_{t\in[0,t_0]}\left\|\left[T_n(t)-T(t)\right]y\right\|\right\}=0\quad,\quad(\forall)y\in
R(\lambda;A){\cal
D}(A).
$$
Comme ${R(\lambda;A){\cal D}(A)}={\cal D}\left(A^2\right)$, compte tenu du th\'eor\`eme \ref{num37} on voit que $\overline{R(\lambda;A){\cal D}(A)}={\cal E}$. Nous obtenons finalement:
$$
\lim_{n\rightarrow\infty}\left\{\sup_{t\in[0,t_0]}\left\|T_n(t)x-T(t)x\right\|\right\}=0\quad,\quad(\forall)x\in{\cal
E}.
$$
$ii)\Longrightarrow i)$ En appliquant le th\'eor\`eme \ref{num13}, nous obtenons pour
$\lambda\in\Lambda_\omega$:
$$
\left[R(\lambda;A_n)-R(\lambda;A)\right]x=\int\limits_{0}^{\infty}\!\!e^{-\lambda
t}\left[T_n(t)-T(t)\right]x\:dt\quad,\quad(\forall)x\in{\cal E},
$$
d'o\`u il r\'esulte:
$$
\left\|\left[R(\lambda;A_n)-R(\lambda;A)\right]x\right\|\leq\int\limits_{0}^{\infty}\!\!e^{-\scriptstyle{Re}\lambda
t}\left\|\left[T_n(t)-T(t)\right]x\right\|\:dt\quad,\quad(\forall)x\in{\cal
E}.
$$
Mais:
$$
\left\|\left[T_n(t)-T(t)\right]x\right\|\leq 2Me^{\omega t}\|x\|
$$
quels que soient $x\in{\cal E}$, $t\geq 0$ et $n\in{\bf N}^*$. Dans ce
cas, en posant:
$$
f_n(t)=e^{\lambda
t}\left[T_n(t)-T(t)\right]x\quad,\quad(\forall)t\geq 0
$$
on voit que:
$$
\left\|f_n(t)\right\|\leq
2Me^{-(\scriptstyle{Re}\lambda-\omega)t}\|x\|\quad,\quad(\forall)t\geq
0.
$$
De plus, compte tenu de l'in\'egalit\'e:
$$
\left\|f_n(t)\right\|\leq e^{-\scriptstyle{Re}\lambda
t}\left\|\left[T_n(t)-T(t)\right]x\right\|\quad,
$$
nous obtenons:
$$
\lim_{n\rightarrow\infty}\left\|f_n(t)\right\|=0\quad,\quad(\forall)t\geq
0.
$$
Avec le th\'eor\`eme de la convergence domin\'ee de Lebesgue
(\cite[Theorem III.3.7, pag. 124]{dunford-schwartz}), il vient:
$$
\lim_{n\rightarrow\infty}\left\|\left[R(\lambda;A_n)-R(\lambda;A)\right]x\right\|=0
$$
pour tout $x\in{\cal E}$ et tout $\lambda\in\Lambda_\omega$. Donc $A_n\longrightarrow A$, si
$n\rightarrow\infty$, pour la topologie forte de la
r\'esolvante.\fin

Une version int\'eressante du th\'eor\`eme \ref{num20} est le
th\'eor\`eme suivant.

\begin{teo}\label{num31}
Soient $\left(\left\{T_\alpha(t)\right\}_{t\geq 0}\right)_{\alpha>0}\subset{{\cal SG}(M,\omega)}$ ayant pour g\'en\'erateurs
infinit\'esimaux les op\'erateurs $(A_\alpha)_{\alpha>0}\subset{{\cal GI}(M,\omega)}$
et \semi ayant pour g\'en\'erateur infinit\'esimal l'op\'erateur $A\in{{\cal
GI}(M,\omega)}$.\\
Les affirmations suivantes sont \'equivalentes:\\
i) pour tout $x\in{\cal D}(A)$, il existe $x_\alpha\in{\cal
D}(A_\alpha)$ tel que: 
$$
\lim_{\alpha\searrow 0}x_\alpha=x
$$
et:
$$
\lim_{\alpha\searrow 0}A_\alpha x_\alpha=Ax;
$$
ii) pour tout $\lambda\in\Lambda_\omega$, on a:
$$
\lim_{\alpha\searrow
0}R(\lambda;A_\alpha)x=R(\lambda;A)x\quad,\quad(\forall)x\in{\cal E};
$$
iii) pour tout $t_0\in]0,\infty)$, nous avons:
$$
\lim_{\alpha\searrow 0}\left\{\sup_{t\in[0,t_0]}\left\|T_\alpha(t)x-T(t)x\right\|\right\}=0\quad,\quad(\forall)x\in{\cal
E}.
$$
\end{teo}
\dem
$i)\Longrightarrow ii)$ Soient $\lambda\in\Lambda_\omega$ et $x\in{\cal D}(A)$. Alors il existe $x_\alpha\in{\cal D}(A_\alpha)$, $\alpha>0$ tel que
$$
\lim_{\alpha\searrow 0}x_\alpha=x
$$
et:
$$
\lim_{\alpha\searrow 0}A_\alpha x_\alpha=Ax\quad.
$$
Nous d\'efinissons 
$$
y=(\lambda I-A)x\in(\lambda I-A){\cal D}(A)
$$
et
$$
y_\alpha=(\lambda
I-A_\alpha)x_\alpha\in(\lambda I-A_\alpha){\cal D}(A_\alpha)\quad,\quad\alpha>0.
$$
Il r\'esulte que $x=R(\lambda;A)y$ et
$x_\alpha=R(\lambda;A_\alpha)y_\alpha$. 
Compte tenu des \'egalit\'es du (i), nous obtenons:
$$
\lim_{\alpha\searrow 0}R(\lambda;A_\alpha)y_\alpha=R(\lambda;A)y
$$
et:
$$
\lim_{\alpha\searrow 0}A_\alpha R(\lambda;A_\alpha)y_\alpha=AR(\lambda;A)y\quad.
$$
On voit que cette derni\`ere \'egalit\'e devient:
$$
\lim_{\alpha\searrow 0}(\lambda I-\lambda I+A_\alpha)R(\lambda;A_\alpha)y_\alpha=(\lambda I-\lambda I+A)R(\lambda;A)y
$$
ou bien
\begin{eqnarray*}
& &\lim_{\alpha\searrow 0}\lambda R(\lambda;A_\alpha)y_\alpha-\lim_{\alpha\searrow 0}(\lambda I-A_\alpha)R(\lambda;A_\alpha)y_\alpha=\\
&=&\lambda R (\lambda;A)y-(\lambda I-A)R(\lambda;A)y\quad.
\end{eqnarray*}
Il vient:
$$
\lambda R(\lambda;A)y-\lim_{\alpha\searrow 0}y_\alpha=\lambda R(\lambda;A)y-y\quad,
$$
d'o\`u:
$$
\lim_{\alpha\searrow 0}y_\alpha=y\quad.
$$
D'autre part, pour tout $\alpha>0$ on a:
$$
\left\|R(\lambda;A_\alpha)\right\|\leq\frac{M}{\mbox{\em
Re}\lambda-\omega}
$$
et pour $y\in(\lambda I-A){\cal D}(A)$ on voit que:
$$
R(\lambda;A_\alpha)y=R(\lambda;A_\alpha)(y-y_\alpha+y_\alpha)=R(\lambda;A_\alpha)(y-y_\alpha)+R(\lambda;A_\alpha)y_\alpha\quad.
$$
Par suite:
\begin{eqnarray*}
& &\left\|R(\lambda;A_\alpha)y-R(\lambda;A)y\right\|\leq\\
&\leq&\left\|R(\lambda;A_\alpha)(y-y_\alpha)\right\|+\left\|R(\lambda;A_\alpha)y_\alpha-R(\lambda;A)y\right\|\leq\\
&\leq&\frac{M}{\mbox{\em
Re}\lambda-\omega}\|y-y_\alpha\|+\left\|R(\lambda;A_\alpha)y_\alpha-R(\lambda;A)y\right\|\quad,
\end{eqnarray*}
d'o\`u il vient:
$$
\lim_{\alpha\searrow 0}R(\lambda;A_\alpha)y=R(\lambda;A)y\quad,
$$
pour tout $y\in(\lambda I-A){\cal D}(A)$. Comme $\overline{(\lambda I-A){\cal
D}(A)}={\cal E}$, on voit que:
$$
\lim_{\alpha\searrow
0}R(\lambda;A_\alpha)x=R(\lambda;A)x\quad,\quad(\forall)x\in{\cal E}.
$$
$ii)\Longrightarrow i)$ Soient $\lambda\in\Lambda_\omega$ et $x\in{\cal
E}$ tel que:
$$
\lim_{\alpha\searrow
0}R(\lambda;A_\alpha)x=R(\lambda;A)x\quad.
$$
Si nous d\'efinissons:
$$
y_\alpha=R(\lambda;A_\alpha)x\in{\cal D}(A_\alpha)
$$
et:
$$
y=R(\lambda;A)x\in{\cal D}(A)\quad,
$$
nous obtenons:
$$
\lim_{\alpha\searrow 0}y_\alpha=y\quad.
$$
De plus:
\begin{eqnarray*}
& &\lim_{\alpha\searrow 0}A_\alpha y_\alpha=\lim_{\alpha\searrow
0}A_\alpha R(\lambda;A_\alpha)x=\lim_{\alpha\searrow 0}[\lambda
R(\lambda;A_\alpha)x-x]=\\
&=&\lambda R(\lambda;A)x-x=AR(\lambda;A)x=Ay\quad.
\end{eqnarray*}
$ii)\Longleftrightarrow iii)$ Cette \'equivalence s'obtient avec une
preuve analogue \`a celle du th\'eor\`eme \ref{num20}.\fin

\begin{cor}\label{num32}
Soient $\left(\left\{T_\alpha(t)\right\}_{t\geq 0}\right)_{\alpha>0}\subset{{\cal SG}(M,\omega)}$ ayant pour g\'en\'erateurs
infinit\'esimaux les op\'erateurs $(A_\alpha)_{\alpha>0}\subset{{\cal GI}(M,\omega)}$
et \semi ayant pour g\'en\'erateur infinit\'esimal l'op\'erateur $A\in{{\cal
GI}(M,\omega)}$. Supposons que pour tout $x\in{\cal D}(A)$, il existe
$\delta>0$ tel que pour tout $\alpha\in]0,\delta[$ on ait $x\in{\cal
D}(A_\alpha)$ et $\lim_{\alpha\searrow 0}A_\alpha x=Ax$. Alors, pour
tout $t_0\in]0,\infty)$ nous avons:
$$
\lim_{\alpha\searrow 0}\left\{\sup_{t\in[0,t_0]}\left\|T_\alpha(t)x-T(t)x\right\|\right\}=0\quad,\quad(\forall)x\in{\cal
E}.
$$
\end{cor}
\dem
Dans le th\'eor\`eme \ref{num31}, nous pouvons prendre $x_\alpha=x$, $(\forall)\alpha\in]0,\delta[$.\fin

Le th\'eor\`eme suivant montre que sous certaines conditions, ${\cal
GI}(M,\omega)$ est une sous-classe ferm\'ee dans ${\cal GI}({\cal E})$.

\begin{teo}\label{num21}
Soient $(A_n)_{n\in{\bf N}^*}\subset{{\cal GI}(M,\omega)}$ et
$\lambda_0\in\Lambda_\omega$ tel que:\\
i) $\left(R(\lambda_0;A_n)\right)_{n\in{\bf N}^*}$ est fortement convergente vers
$R_{\lambda_0}\in{\cal B(E)}$;\\
ii) $\overline{{\cal I}\mbox{m }R_{\lambda_0}}={\cal E}$.\\
Alors il existe un unique op\'erateur $A\in{\cal GI}(M,\omega)$ tel
que $R_{\lambda_0}=R(\lambda_0;A)$.
\end{teo}
\dem
Nous notons:
$$
{\cal S}=\left\{\lambda\in\Lambda_\omega\left|\left(R(\lambda;A_n)\right)_{n\in{\bf N}^*}\mbox{\em est fortement
convergente}\right.\right\}.
$$
Montrons que ${\cal S}=\Lambda_\omega$.\\
Prouvons que ${\cal S}$ est ensemble ouvert dans
$\Lambda_\omega$.
Soit $\mu\in{\cal S}$. Pour tout $n\in{\bf N}^*$, l'application:
$$
\rho(A_n)\ni\lambda\longmapsto R(\lambda;A_n)\in{\cal B(E)}
$$
est analytique et nous avons:
\begin{eqnarray*}
R(\lambda;A_n)&=&\sum\limits_{k=0}^{\infty}\frac{(\lambda-\mu)^k}{k!}\frac{d^k}{d\mu^k}R(\mu;A_n)=\\
&=&\sum\limits_{k=0}^{\infty}\frac{(\lambda-\mu)^k}{k!}(-1)^kk!{R(\mu;A_n)}^{k+1}=\\
&=&\sum\limits_{k=0}^{\infty}(\mu-\lambda)^k{R(\mu;A_n)}^{k+1}\quad.
\end{eqnarray*}
Comme $A_n\in{\cal GI}(M,\omega)$ implique:
$$
\left\|{R(\mu;A_n)}^k\right\|\leq\frac{M}{(\mbox{\em
Re}\mu-\omega)^k}\quad,\quad(\forall)k\in{\bf N}^*,
$$
on voit que:
$$\left\|R(\lambda;A_n)\right\|\leq\sum\limits_{k=0}^{\infty}|\mu-\lambda|^k\left\|R(\mu;A)^{k+1}\right\|\leq\frac{M}{\mbox{\em
Re}\mu-\omega}\sum\limits_{k=0}^{\infty}\left(\frac{|\mu-\lambda|}{\mbox{\em
Re}\mu-\omega}\right)^k\quad.
$$
La s\'erie de la partie droite de cette in\'egalit\'e est convergente
sur l'ensemble:
$$
{\cal V}=\left\{\lambda\in\Lambda_\omega\left|\:|\mu-\lambda|(\mbox{\em Re}\mu-\omega)^{-1}<1\right.\right\}.
$$
Il en r\'esulte que la s\'erie:
$$
R(\lambda;A_n)=\sum\limits_{k=0}^{\infty}(\mu-\lambda)^k{R(\mu;A_n)}^{k+1}
$$
est uniform\'ement convergente sur les compacts
$$
{\cal V_\nu}=\left\{\lambda\in\Lambda_\omega\left|\:|\mu-\lambda|(\mbox{\em Re}\mu-\omega)^{-1}\leq\nu<1\right.\right\}\subset{\cal V}\quad.
$$
Comme
$$
\left\|R(\lambda;A)\right\|\leq\frac{M}{Re\:\mu-\omega}\sum\limits_{k=0}^{\infty}\nu^k\quad,
$$
on voit que la suite $\left(R(\lambda;A_n)\right)_{n\in{\bf N}^*}$ est fortement convergente pour tout $\lambda\in{\cal V}_\nu$. Donc il existe un voisinage de $\mu$ contenu
dans ${\cal S}$. Par cons\'equent ${\cal S}$ est ensemble ouvert dans
$\Lambda_\omega$.\\
Maintenant, nous allons montrer que ${\cal S}$ est un ensemble relativement ferm\'e dans $\Lambda_\omega$. Soient $\left(\lambda_m\right)_{m\in{\bf N}}\subset{\cal S}$ et $\lambda\in\Lambda_\omega$ tel que
$$
\lambda=\lim_{m\rightarrow\infty}\lambda_m\quad.
$$
Pour tout $\nu\in]0,1[$, il existe $\lambda_{m,\nu}\in{\cal S}$ tel que:
$$
\left|\lambda_{m,\nu}-\lambda\right|\left(Re\:\lambda_{m,\nu}-\omega\right)^{-1}\leq\nu<1\quad.
$$
Compte tenu de la premi\`ere partie de la preuve, on voit que la s\'erie
$$
R(\lambda;A_n)=\sum\limits_{k=0}^{\infty}\left(\lambda_{m,\nu}-\lambda\right)^kR(\lambda_{m,\nu};A_n)^{k+1}
$$
est uniform\'ement convergente et que la suite $(R(\lambda;A_n))_n\in{\bf N}^*$ est fortement convergente. Par cons\'equent, $\lambda\in{\cal S}$ et $\cal S$ est un ensemble relativement ferm\'e dans $\Lambda_\omega$. Comme $\lambda_0\in{\cal S}$, nous voyons que ${\cal S}=\Lambda_\omega$ par connexit\'e.\\
Pour $\lambda\in\Lambda_\omega$, d\'efinissons l'op\'erateur
$R_\lambda\in{\cal B(E)}$ par:
$$
R_\lambda
x=\lim_{n\rightarrow\infty}R(\lambda;A_n)x\quad,\quad(\forall)x\in{\cal
E}.
$$
Soient $\lambda\:,\:\mu\in\Lambda_\omega$ arbitraires. On a:
\begin{eqnarray*}
\left(R_\lambda-R_\mu\right)x&=&\lim_{n\rightarrow\infty}\left[R(\lambda;A_n)-R(\mu;A_n)\right]x=\\
&=&\lim_{n\rightarrow\infty}(\mu-\lambda)R(\lambda;A_n)R(\mu;A_n)x=\\
&=&(\mu-\lambda)R_\lambda R_\mu x\quad,\quad(\forall)x\in{\cal E}.
\end{eqnarray*}
Par cons\'equent $R_\lambda$ est une pseudo-r\'esolvante, quel que soit
$\lambda\in\Lambda_\omega$. Comme il existe $\lambda_0\in\Lambda_\omega$
tel que $\overline{{\cal I}\mbox{\em m }R_{\lambda_0}}={\cal E}$,
compte tenu du th\'eor\`eme \ref{num19} (ii), on d\'eduit que
$\overline{{\cal I}\mbox{\em m }R_\lambda}={\cal E}$, quel que soit
$\lambda\in\Lambda_\omega$.
Avec l'in\'egalit\'e:
$$
\left\|{R(\lambda;A_n)}^m\right\|\leq\frac{M}{(\mbox{\em
Re}\lambda-\omega)^m}\quad,\quad(\forall)\lambda\in\Lambda_\omega\mbox{ et }m\in{\bf N}^*,
$$
on voit que pour tout compact $K\subset\Lambda_\omega$, il existe $M_K>0$ tel que
$$
\sup_{\lambda\in K}\|R(\lambda;A_n)\|\leq M_K\quad,
$$
quel que soit $n\in{\bf N}^*$. Avec le lemme de Montel (\cite[pag. 220]{gaspar-suciu-3}), on d\'eduit qu'il existe une sous-suite $\left(R(\lambda;A_{n_k})\right)_{k\in{\bf N}^*}$ telle que
$$
R_\lambda x=\lim_{k\rightarrow\infty}R(\lambda;A_{n_k})x\quad,\quad(\forall)x\in{\cal E},
$$
uniform\'ement par rapport \`a $\lambda$ sur les compacts de $\Lambda_\omega$. Comme $R(\lambda;A_{n_k})$ est un op\'erateur injectif pour tout $k\in{\bf N}^*$, avec le th\'eor\`eme de Hurwitz (\cite[pag. 193]{gaspar-suciu-3}) nous obtenons que $R_\lambda$ est un op\'erateur injectif, donc ${\cal K}\mbox{\em er
}R(\lambda;A_n)=\{0\}$. En appliquant le th\'eor\`eme
\ref{num19} (iii), on voit que pour tout $\lambda\in\Lambda_\omega$, il existe un op\'erateur
lin\'eaire $A:{\cal D}(A)\longrightarrow{\cal E}$,  $A=\lambda
I-R_\lambda^{-1}$ ferm\'e et d\'efini sur un sous espace dense tel que
$R_\lambda=R(\lambda;A)$, $(\forall)\lambda\in\Lambda_\omega$. De plus:
$$
\left\|{R(\lambda;A)}^m\right\|\leq\frac{M}{\left(\mbox{\em
Re}\lambda-\omega\right)^m}\quad.
$$
et le th\'eor\`eme de Hille - Yosida implique alors que $A\in{\cal
GI}(M,\omega)$.\fin

Maintenant, nous avons toutes les conditions pour formuler un autre
r\'esultat important concernant les $C_0$-semi-groupes.

\begin{teo}[Trotter - Kato]\label{num54}
Soit $\left(\left\{T_n(t)\right\}_{t\geq 0}\right)_{n\in{\bf
N}^*}\subset{\cal SG}(M,\omega)$ ayant pour g\'en\'erateurs
infinit\'esimaux les op\'erateurs $\left(A_n\right)_{n\in{\bf N}^*}\subset{\cal
GI}(M,\omega)$.\\
S'il existe $\lambda_0\in\Lambda_\omega$ tel que:\\
i) $\left(R(\lambda_0;A_n)\right)_{n\in{\bf N}^*}$ est fortement convergente vers
$R_{\lambda_0}\in{\cal B(E)}$;\\
ii) $\overline{{\cal I}\mbox{m }R_{\lambda_0}}={\cal E}$,\\
alors il existe un unique op\'erateur $A\in{\cal GI}(M,\omega)$ tel
que $R_\lambda=R(\lambda;A)$, $(\forall)\lambda\in\Lambda_\omega$.
De plus, si $\left\{T(t)\right\}_{t\geq 0}$ est le $C_0$-semi-groupe
engendr\'e par $A$, alors pour tout $t_0\in]0,\infty)$ on a:
$$
\lim_{n\rightarrow\infty}\left\{\sup_{t\in[0,t_0]}\left\|T_n(t)x-T(t)x\right\|\right\}=0\quad,\quad(\forall)x\in{\cal
E}.
$$
\end{teo}
\dem
Les affirmations du th\'eor\`eme r\'esultent du th\'eor\`eme \ref{num20} et du th\'eor\`eme
\ref{num21}.\fin
\vspace{2cm}

    \section{Formule de Lie - Trotter pour les $C_0$-semi-groupes}

\hspace{1cm}Dans la suite, nous montrons le th\'eor\`eme de repr\'esentation
g\'en\'erale, la formule exponentielle et la formule de Lie-Trotter pour les semi-groupes fortement
continus. Nous commen\c cons par un r\'esultat technique.

\begin{lema}\label{num33}
Soient $T\in{\cal B(E)}$ et $M,N\geq 1$ tel que:
$$
\left\|T^k\right\|\leq MN^k\quad,\quad(\forall)k\in{\bf N}^*.
$$
Alors, pour tout $n\in{\bf N}$, nous avons:
$$\left\|e^{n(T-I)}x-T^nx\right\|\leq
MN^{n-1}e^{(N-1)n}\sqrt{n^2(N-1)^2+nN}\|Tx-x\|
$$
pour tout $x\in{\cal E}$.
\end{lema}
\dem
Soient $k,n\in{\bf N}$ tel que $k\geq n$. Alors, nous avons:
\begin{eqnarray*}
&
&\left\|T^kx-T^nx\right\|=\left\|\sum\limits_{i=n}^{k-1}\left(T^{i+1}x-T^ix\right)\right\|\leq\\
&\leq&\sum\limits_{i=n}^{k-1}\left\|T^i\right\|\|Tx-x\|\leq \|Tx-x\|\sum\limits_{i=n}^{k-1}MN^i\leq\\
&\leq&M\|Tx-x\|\sum\limits_{i=n}^{k-1}N^{k-1}=(k-n)MN^{k-1}\|Tx-x\|\leq\\
&\leq&|k-n|MN^{n+k-1}\|Tx-x\|\quad,\quad(\forall)x\in{\cal
E}.
\end{eqnarray*}
Compte tenu de la sym\'etrie, il est clair que cette in\'egalit\'e
reste valable si nous consid\'erons $n>k$. Par suite, on voit que:
$$
\left\|T^kx-T^nx\right\|\leq|k-n|MN^{n+k-1}\|Tx-x\|\quad,\quad(\forall)x\in{\cal
E}\mbox{ et }n,k\in{\bf N}.
$$
Si $t\geq 0$ et $n\in{\bf N}$, alors nous avons:
\begin{eqnarray*}
& &\left\|e^{t(T-I)}x-T^nx\right\|=\left\|e^{-t}\sum\limits_{k=0}^{\infty}\frac{t^k}{k!}\left(T^kx-T^nx\right)\right\|\leq\\
&\leq&e^{-t}\sum\limits_{k=0}^{\infty}\frac{t^k}{k!}\left\|T^kx-T^nx\right\|\leq
MN^{n-1}e^{-t}\|Tx-x\|\sum\limits_{k=0}^{\infty}\frac{(tN)^k}{k!}|k-n|\quad.
\end{eqnarray*}
Avec l'in\'egalit\'e de Cauchy-Schwartz, il vient:
\begin{eqnarray*}
& &\sum\limits_{k=0}^{\infty}\frac{(tN)^k}{k!}|k-n|=\sum\limits_{k=0}^{\infty}\left(\sqrt{\frac{(tN)^k}{k!}}\right)\left(\sqrt{\frac{(tN)^k}{k!}}|k-n|\right)\leq\\
&\leq&\sqrt{\sum\limits_{k=0}^{\infty}\frac{(tN)^k}{k!}}\sqrt{\sum\limits_{k=0}^{\infty}\frac{(tN)^k}{k!}(k-n)^2}=e^{tN}\sqrt{(n-Nt)^2+Nt}\quad.
\end{eqnarray*}
Il s'ensuit que:
$$\left\|e^{t(T-I)}x-T^nx\right\|\leq
MN^{n-1}e^{(N-1)t}\sqrt{(n-Nt)^2+Nt}\|Tx-x\|
$$
quel que soit $x\in{\cal E}$. Finalement, en prenant $t=n$, nous obtenons l'in\'egalit\'e consid\'er\'ee dans l'\'enonc\'e.\fin

\begin{teo}[de repr\'esentation g\'en\'erale]
Soit $\left\{F(t)\right\}_{t\geq 0}\subset{\cal B(E)}$ une famille
d'op\'erateurs lin\'eaires born\'es avec $F(0)=I$. Supposons qu'il
existe $\omega\geq 0$ et $M\geq 1$ tel que:
$$
\left\|{F(t)}^k\right\|\leq Me^{k\omega t}\quad,\quad(\forall)k\in{\bf
N}^*,
$$
pour tout $t\geq 0$.\\
Si $A$ est le g\'en\'erateur infinit\'esimal d'un $C_0$-semi-groupe
$\left\{T(t)\right\}_{t\geq 0}$ tel que:
$$
\lim_{t\searrow 0}\frac{F(t)x-x}{t}=Ax\quad,\quad(\forall)x\in{\cal
D}(A),
$$
alors nous avons:
$$
T(t)x=\lim_{n\rightarrow\infty}\left[F\left(\frac{t}{n}\right)\right]^nx\quad,\quad(\forall)x\in{\cal
E},
$$
uniform\'ement par rapport \`a $t$ sur les intervalles compacts de
$[0,\infty)$.
\end{teo}
\dem
Soient $0\leq a<b$. Pour $t\in[a,b]$, d\'efinissons:
$$
A_n=\frac{n}{t}\left[F\left(\frac{t}{n}\right)-I\right]\quad,\quad(\forall)n\in{\bf
N}^*.
$$
Il est clair que $A_n\in{\cal B(E)}$, $(\forall)n\in{\bf N}^*$,
d'o\`u il r\'esulte que pour tout $n\in{\bf N}^*$, $A_n$ est le
g\'en\'erateur infinit\'esimal du semi-groupe uniform\'ement continu
$\left\{e^{tA_n}\right\}_{t\geq 0}$ et, de plus, nous avons:
$$
\lim_{n\rightarrow\infty}A_nx=Ax\quad,\quad(\forall)x\in{\cal D}(A).
$$
Avec le corollaire \ref{num32}, nous voyons que:
$$
\lim_{n\rightarrow\infty}e^{tA_n}x=T(t)x\quad,\quad(\forall)x\in{\cal
E},
$$
uniform\'ement par rapport \`a $t\in[a,b]$. Compte tenu du lemme \ref{num33}, si $x\in{\cal D}(A)$,
il vient:
\begin{eqnarray*}
&
&\left\|e^{tA_n}x-\left[F\left(\frac{t}{n}\right)\right]^nx\right\|=\left\|e^{n\left[F\left(\frac{t}{n}\right)-I\right]}x-\left[F\left(\frac{t}{n}\right)\right]^nx\right\|\leq\\
&\leq&Me^{\omega\frac{t}{n}(n-1)}e^{\left(e^{\omega\frac{t}{n}}-1\right)n}\sqrt{n^2\left(e^{\omega\frac{t}{n}}-1\right)^2+ne^{\omega\frac{t}{n}}}\left\|F\left(\frac{t}{n}\right)x-x\right\|=\\
&=&Me^{\omega\frac{t}{n}(n-1)+\left(e^{\omega\frac{t}{n}}-1\right)n}\sqrt{n^2\left(e^{\omega\frac{t}{n}}-1\right)^2+ne^{\omega\frac{t}{n}}}\frac{t}{n}\left\|\frac{F\left(\frac{t}{n}\right)x-x}{\frac{t}{n}}\right\|=\\
&=&Me^{\omega\frac{t}{n}(n-1)+\frac{e^{\omega\frac{t}{n}}-1}{\frac{t}{n}}t}\sqrt{t^2\left(e^{\omega\frac{t}{n}}-1\right)^2+\frac{t^2}{n}e^{\omega\frac{t}{n}}}\left\|\frac{F\left(\frac{t}{n}\right)x-x}{\frac{t}{n}}\right\|\quad,
\end{eqnarray*}
d'o\`u:
$$
\left\|e^{tA_n}x-\left[F\left(\frac{t}{n}\right)\right]^nx\right\|\longrightarrow
0\quad\mbox{si}\quad n\rightarrow\infty,
$$
pour tout $x\in{\cal D}(A)$, uniform\'ement par rapport \`a $t\in[a,b]$. De plus, pour tout
$x\in{\cal D}(A)$, nous avons:
\begin{eqnarray*}
&
&\left\|T(t)x-\left[F\left(\frac{t}{n}\right)\right]^nx\right\|\leq\\
&\leq&\left\|T(t)x-e^{tA_n}x\right\|+\left\|e^{tA_n}x-\left[F\left(\frac{t}{n}\right)\right]^nx\right\|\longrightarrow
0\quad\mbox{si}\quad n\rightarrow\infty,
\end{eqnarray*}
d'o\`u l'on d\'eduit que:
$$
T(t)x=\lim_{n\rightarrow\infty}\left[F\left(\frac{t}{n}\right)\right]^nx\quad,\quad(\forall)x\in{\cal
D}(A),
$$
uniform\'ement par rapport \`a $t$ sur les intervalles compacts de
$[0,\infty)$.\\
Comme $\overline{{\cal D}(A)}={\cal E}$ et
$\left\|F\left(\frac{t}{n}\right)^n\right\|\leq Me^{\omega t}$, on voit que:
$$
T(t)x=\lim_{n\rightarrow\infty}\left[F\left(\frac{t}{n}\right)\right]^nx\quad,\quad(\forall)x\in{\cal
E},
$$
uniform\'ement par rapport \`a $t$ sur les intervalles compacts de
$[0,\infty)$.\fin

\begin{teo}[la formule exponentielle]
Soient \semi et $A$ son g\'en\'erateur infinit\'esimal. Alors:
$$
T(t)x=\lim_{n\rightarrow\infty}\left(I-\frac{t}{n}A\right)^{-n}x=\lim_{n\rightarrow\infty}\left[\frac{n}{t}R\left(\frac{n}{t};A\right)\right]^nx\quad,\quad(\forall)x\in{\cal
E},
$$
uniform\'ement par rapport \`a $t$ sur les intevalles compacts de $[0,\infty)$.
\end{teo}
\dem
Pour $A\in{\cal GI}(M,\omega)$ et
$t\in\left]0,\frac{1}{\omega}\right[$, nous d\'efinissons:
$$
F(t)=(I-tA)^{-1}=\frac{1}{t}R\left(\frac{1}{t};A\right)\quad.
$$
Compte tenu du th\'eor\`eme de Hille-Yosida, on voit que:
$$
\left\|F(t)^k\right\|=\left\|\left[\frac{1}{t}R\left(\frac{1}{t};A\right)\right]^k\right\|\leq\left(\frac{1}{t}\right)^k\frac{M}{\left(\frac{1}{t}-\omega\right)^k}=\frac{M}{(1-\omega
t)^k}\quad.
$$
Comme:
$$
\frac{1}{(1-\omega t)^k}\leq e^{\frac{k\omega t}{1-\omega t}}\quad,
$$
il vient:
$$
\left\|F(t)^k\right\|\leq Me^{2k\omega t}
$$
pour $t\in\left]0,\frac{1}{2\omega}\right]$. D'autre part, avec le
lemme \ref{num16} nous obtenons:
\begin{eqnarray*}
& &\lim_{t\searrow 0}\frac{F(t)x-x}{t}=\lim_{t\searrow
0}\frac{1}{t}\left[F(t)-I\right]x=\lim_{t\searrow 0}AF(t)x=\\
&=&\lim_{t\searrow
0}A\left[\frac{1}{t}R\left(\frac{1}{t};A\right)x\right]=Ax\quad,\quad(\forall)x\in{\cal
D}(A).
\end{eqnarray*}
Compte tenu du th\'eor\`eme de repr\'esentation g\'en\'erale, on voit
que:
\begin{eqnarray*}
& &T(t)x=\lim_{n\rightarrow\infty}\left[F\left(\frac{t}{n}\right)\right]^nx=\lim_{n\rightarrow\infty}\left(I-\frac{t}{n}A\right)^{-n}x=\\
&=&\lim_{n\rightarrow\infty}\left[\frac{n}{t}R\left(\frac{n}{t};A\right)\right]^nx\quad,\quad(\forall)x\in{\cal
E},
\end{eqnarray*}
uniform\'ement par rapport \`a $t$ sur les intervalles compactes de $[0,\infty)$.\fin

\begin{teo}[la formule de Lie-Trotter]\label{num48}
Soient $A_1\in{\cal GI}(M_1,\omega_1)$ le g\'en\'erateur
infinit\'esimal du semi-groupe $\left\{T_1(t)\right\}_{t\geq
0}\in{\cal SG}(M_1,\omega_1)$, respectivement  $A_2\in{\cal GI}(M_2,\omega_2)$ le g\'en\'erateur
infinit\'esimal du semi-groupe $\left\{T_2(t)\right\}_{t\geq
0}\in{\cal SG}(M_2,\omega_2)$.
Supposons qu'il existe $\omega\geq 0$ et $M\geq 1$ tel que:
$$
\left\|\left[T_1(t)T_2(t)\right]^k\right\|\leq Me^{k\omega
t}\quad,\quad(\forall)k\in{\bf N}^*.
$$
Si l'op\'erateur
$$
A:{\cal D}(A)\subset{\cal E}\longrightarrow{\cal E},
$$
d\'efini par:
$$
Ax=A_1x+A_2x\quad,\quad(\forall)x\in{\cal D}(A)={\cal D}(A_1)\cap{\cal
D}(A_2),
$$
est le g\'en\'erateur infinit\'esimal d'un $C_0$-semi-groupe
$\left\{T(t)\right\}_{t\geq 0}$, alors:
$$
T(t)x=\lim_{n\rightarrow\infty}\left[T_1\left(\frac{t}{n}\right)T_2\left(\frac{t}{n}\right)\right]^nx\quad,\quad(\forall)x\in{\cal
E},
$$
uniform\'ement par rapport \`a $t$ sur les intervalles compacts de
$[0,\infty)$.
\end{teo}
\dem
Soit:
$$
F:[0,\infty)\longrightarrow{\cal B(E)}
$$
$$
F(t)=T_1(t)T_2(t)\quad,\quad(\forall)t\geq 0.
$$
Il est \'evident que $F(0)=I$. De plus, pour $x\in{\cal D}(A)$, nous
avons:
\begin{eqnarray*}
& &\lim_{t\searrow 0}\frac{F(t)x-x}{t}=\lim_{t\searrow
0}\frac{T_1(t)T_2(t)x-x}{t}=\\
&=&\lim_{t\searrow 0}\frac{T_1(t)T_2(t)x-T_1(t)x}{t}+\lim_{t\searrow
0}\frac{T_1(t)x-x}{t}=\\
&=&\lim_{t\searrow 0}T_1(t)\frac{T_2(t)x-x}{t}+A_1x=A_1x+A_2x=Ax\quad.
\end{eqnarray*}
Avec le th\'eor\`eme de repr\'esentation g\'en\'erale, on voit que:
$$
T(t)x=\lim_{n\rightarrow\infty}\left[F\left(\frac{t}{n}\right)\right]^nx=\lim_{n\rightarrow\infty}\left[T_1\left(\frac{t}{n}\right)T_2\left(\frac{t}{n}\right)\right]^nx\quad,\quad(\forall)x\in{\cal
E},
$$
uniform\'ement par rapport \`a $t$ sur les intervalles compacts de
$[0,\infty)$.\fin

\begin{obs}
\em
Si $A\in{\cal GI}(M,\omega)$, compte tenu de la formule
exponentielle, on peut d\'efinir:
$$
e^{tA}x=\lim_{n\rightarrow\infty}\left(I-\frac{t}{n}A\right)^{-n}x=\lim_{n\rightarrow\infty}\left[\frac{n}{t}R\left(\frac{n}{t};A\right)\right]^nx\quad,\quad(\forall)x\in{\cal
E},
$$
uniform\'ement par rapport \`a $t$ sur les intervalles compacts de
$[0,\infty)$. Avec cette notation, dans les hypoth\`eses du
th\'eor\`eme \ref{num48}, nous obtenons pour la formule de Lie-Trotter
l'expression:
$$
e^{tA}x=\lim_{n\rightarrow\infty}\left[e^{\frac{t}{n}A_1}e^{\frac{t}{n}A_2}\right]^nx\quad,\quad(\forall)x\in{\cal
E},
$$
uniform\'ement par rapport \`a $t$ sur les intervalles compacts de $[0,\infty)$.
\end{obs}
\vspace{2cm}

	\section{Notes}

{\footnotesize
Pour les r\'esultats de la section 4.1 on peut consulter
\cite[pag. 35]{kato3}.

Les propri\'et\'es de convergence pour les $C_0$-semi-groupes ont
\'et\'e \'etudi\'ees  par Trotter dans \cite{trotter1}. Pour les
th\'eor\`emes \ref{num20}, \ref{num21}, \ref{num54} on peut consulter
\cite[pag. 84]{pazy1} ou \cite[pag. 131]{ahmed} et pour le th\'eor\`eme
\ref{num31} nous avons utilis\'e \cite[pag. 80]{davies}.

Le th\'eor\`eme \ref{num48} a \'et\'e montr\'e par Trotter dans
\cite{trotter2} et a \'et\'e \'etudi\'e par Chernoff dans
\cite{chernoff1}. Les r\'esultats que nous avons pr\'esent\'es se
trouvent dans \cite[pag. 89]{pazy1}. Dans \cite[pag. 90]{davies}, on
peut trouver ces probl\`emes pour les $C_0$-semi-groupes de contractions.
}

\end{document}